\documentclass[12pt,amsmath,amssym,latexsymb]{article}
\usepackage{amssymb,latexsym}
\usepackage{longtable,tabularx}
\usepackage{makeidx}
\usepackage{graphicx}
\usepackage{psfrag}
 \usepackage{url}
 \usepackage{framed}
\usepackage[T1]{fontenc}
\input xypic
\usepackage{amsmath}
\usepackage{amssymb,verbatim}
\usepackage[usenames]{color}

\def\wtext{}
\def\ptext{}
\def\ftext{}
\def\qtext{}

\def\ktext{}
\def\fftext{}
\def\dtext{}
\def\dtextRB{}
\def\dtextRBX{}

\def\zztext{}

\def\ztext{}

\newcommand{\novtekst}{}
\def\Rtext{}
\def\Mtext{}
\def\stext{}
\def\rtext{}

\def\RobCfortythree{{\ktext C_{3}}}
\def\RobCninety{{\rtext C_{2}}}
\def\RobCeightyfour{{\ktext C_{1}}}

\def\RobCzero{{\rtext C_{4}}}
\def\RobCzeroone{{\rtext C_{5}}}

\def\RobCone{{\rtext C_{6}}}
\def\RobCthree{{\rtext C_{7}}}
\def\RobCfive{{\rtext C_{8}}}
\def\RobCsix{{\rtext C_{9}}}
\def\RobCseven{{\rtext C_{10}}}

\def\MattiCseven{{\rtext C_{11}}}
\def\RobCthirty{{\rtext C_{12}}}
\def\RobCeleven{{\rtext C_{13}}}

 \def\slavaConetheta{{\stext C_{14}}}
 \def\SlavaCnonumber{{\stext C_{15}}}
  \def\Slavacnonumber{c_{1}}

\def\Slavacseven{c_3}

\def\SlavaCthree{{\Rtext C_{16}}}
\def\SlavaCone{{\Rtext C_{17}}}
\def\SlavaCfour{{\Rtext C_{18}}}

\def\SlavaCsixagain{{\stext  C_{20}}}

\def\Slavacsevenagain{{\stext  c_{5}}}

\def\SlavaCnine{{\stext  C_{21}}}

\def\Slavacten{{\stext  C_{7}}}
\def\SlavaCeleven{{\stext  C_{23}}}
\def\SlavaCnineteen{{\stext  C_{24}}}

\def\SlavaCtwelve{{\stext  C_{25}}}
\def\SlavaCthirteen{{\stext  C_{26}}}

\def\SlavaCfourteen{{\stext  C_{27}}}
\def\SlavaCthreeagain{{\stext  C_{28}}}

\def\RobCthirtyfive{{\stext  C_{29}}}
\def\RobCfortyfive{{\stext  C_{30}}}
\def\Cfortyfour{{\stext  C_{31}}}
\def\RobCthirtysix{{\stext  C_{32}}}
\def\RobCfortyseven{{\stext  C_{33}}}

\def\RobCthirtynine{{\rtext  C_{34}}}
\def\RobCeightyone{{\rtext  C_{35}}}
\def\RobCeightytwo{{\rtext  C_{36}}}
\def\RobCfortyeight{{\rtext  C_{37}}}
\def\RobCfortynine{{\rtext  C_{38}}}
\def\RobCeigthy{{\rtext  C_{39}}}
\def\RobCforty{{\rtext  C_{40}}}
\def\RobCfortyone{{\rtext  C_{41}}}

\def\CCfortytwo{{\rtext  C_{42}}}
\def\CCfortythree{{\rtext  C_{43}}}
\def\CCfortyfour{{\rtext  C_{44}}}

\def\origjnull{{\stext j_0}}
\def\jnull{{\stext j_1}} 
\def\janul{{\rtext j_1}} 

 \def\jbasis{{\stext j_0}}

\def\ctext{}
\begin {document}



\def\mtext{}

\def\derho{\gamma}

\def\vektor {\underline }
\def\linetilde{\widetilde}

\def\tilde{\widetilde}
\def \U {\mathcal U}
\def \E {\mathcal E}
\def \bfo {\begin {eqnarray*} }
\def \efo {\end {eqnarray*} }
\def \ba {\begin {eqnarray*} }
\def \ea {\end {eqnarray*} }
\def \beq {\begin {eqnarray}}
\def \eeq {\end {eqnarray}}
\newtheorem {lemma}{Lemma}
\newtheorem {proposition}{Proposition}
\newtheorem {theorem}{Theorem}
\newtheorem {definition}{Definition}

\newtheorem {corollary}{Corollary}

\newtheorem {Remark}{Remark}

\def \ga {{\gamma}}
\def \L {{\cal L}}

\def \F {{\cal F}}

\def \W {{\cal W}}

\def \Box {${\square}$}

\def \Z {{ \mathbb {Z}}}
\def \C {{ \mathbb {C}}}
\def \D {{\cal D}}
\def \R {{ \mathbb {R}}}

\def \H2s {H^{s+1}_0(\partial M\times [0,T/2])}

\def \supp {\hbox{supp }}
\def \diam {\hbox{diam}\,}
\def \dist {\hbox{dist}}

\def\bra{\langle}
\def\cet{\rangle}

\def \e {\varepsilon}

\def \a {\alpha}

\def \pa0 {\partial _0}

\def \p {\partial}

\def\e{\varepsilon}

\def\M{{M}}



\def\tilde{\widetilde}

\def\D{{\cal D}}

\def \mbeq {\begin {eqnarray}}
\def \meeq {\end {eqnarray}}

\def\L{{\cal L}}
\def\la{\lambda}


\def \bfo {\begin {displaymath} }
\def \efo {\end {displaymath} }

\def \beq {\begin {eqnarray}}
\def \eeq {\end {eqnarray}}
\def \ba {\begin {eqnarray*}}
\def \ea {\end  {eqnarray*}}

\def \F {{\cal F}}

\def \H2s {H^{s+1}_0(\partial \M\times [0,T/2])}

\def \supp {\hbox{supp }}
\def \diam {\hbox{diam }}
\def \dist {\hbox{dist}}

\def\bra{\langle}
\def\cet{\rangle}
\def \bD {\bf D}

\def \e {\varepsilon}

\def \a {\alpha}

\def \pa0 {\partial _0}

\def \p {\partial}


\def \vol {\hbox{vol\,}}

\def \bM {{\bf M}}
\def \diam {\hbox {diam}}
\def \inj {\hbox {inj}}


\def \tilde{\widetilde}

\title{Reconstruction and stability in Gel'fand's inverse interior spectral problem }

\author{Roberta Bosi, Yaroslav Kurylev, and Matti Lassas}

\date{December 10, 2019}
\maketitle

\def \bfo {\begin {displaymath} }
\def \efo {\end {displaymath} }
\def \beq {\begin {eqnarray}}
\def \eeq {\end {eqnarray}}

\def \ba {\begin {eqnarray*}}
\def \ea {\end {eqnarray*}}

\def \D {{\cal D}}\def \T {{\cal T}}

\def \W {{\cal W}}

\def \H{{\cal H}}

\def\supp{\hbox{supp }}
\def\diam{\hbox{diam }}
\def\dist{\hbox{dist}}

\def\e{\varepsilon}

\def\F{ {\cal F}}

\def\exp{\,\text{exp}\,}

\def\exp{\hbox{exp}}
 \def\bra{\langle}
\def\cet{\rangle}
\def\p{\partial}

\def\C{{\cal C}}

\def\bra{\langle}
\def\cet{\rangle}

\def\p{\partial }

\def\e{\varepsilon}

\def\a{\alpha}

\begin{abstract}
Assume that $M$ is a compact Riemannian manifold of bounded geometry given by restrictions on its diameter,  Ricci curvature and injectivity radius. Assume we are given, with some error, the first eigenvalues of the Laplacian $\Delta_g$ on $M$ as well as the corresponding eigenfunctions restricted on an open set in $M$. We then construct a stable  approximation to the manifold $(M,g)$. Namely, we construct a metric space and a Riemannian manifold which differ, in a proper sense, just a little from $M$  when the above data are given with a small error. We give an explicit {\dtext $\log\log$-type}  stability estimate on how the constructed manifold and the metric on it depend on the errors in the given data. Moreover a similar stability estimate is derived for the Gel'fand's inverse problem. The proof is based on methods from geometric convergence, a quantitative stability estimate for the unique continuation and a new version of the geometric Boundary Control method.
\end{abstract}

\section{Introduction}

\subsection{Inverse interior spectral data and classes of manifolds}
{\mtext Let $(M,g,p)$ be a pointed   compact Riemannian manifold,
that is, $(M,g)$ is a compact Riemannian manifold without boundary and $p\in M$  is a point on $M$.
Let $\Delta _g$ be the Laplace operator on $(\M,g)$, 
with $0=\la_0 <\la_1\leq \la_2\leq\dots$ being its eigenvalues and
$ \varphi_j$, $j=0,1,2,\dots$ being the complete sequence of
$L^2(M)$-orthonormal
eigenfunctions satisfying $-\Delta_g \varphi_j=\lambda_j \varphi_j$ on $M$.

\begin{definition}\label{def:1} Let $(M, g, p)$ be an $n$ dimensional compact pointed
manifold with $n\ge 2$. Let $r_0>0$. Then
\medskip

\noindent
(i) The pair, consisting of the ball  $(B(p, r_0),g|_{B(p,r_0)})$ on the  Riemannian manifold $M$
and the sequence
  $ \{(\la_j,\, \varphi _j|_{B(p,r_0)});\ j=0, 1,2,\dots\}$  of eigenvalues and eigenfunctions,
   is called the interior spectral data (ISD) of $(\M, g, p)$.
\medskip

\noindent
(ii)
The pair, consisting of the ball  $(B(p, r_0),g|_{B(p, r_0)})$ and a finite collection  $ \{(\la_j,\, \varphi _j|_{B(p, r_0)}),$ $  j=0, 1,2,\dots,J\}$ of the $J+1$  first eigenvalues and eigenfunctions,
 is called the finite interior spectral data (FISD) of $(\M, g, p)$.
\end{definition}

\begin{figure}
\begin {picture}(200,100)(-20,0)
 \put(100,0) {\includegraphics[height=3.0cm]{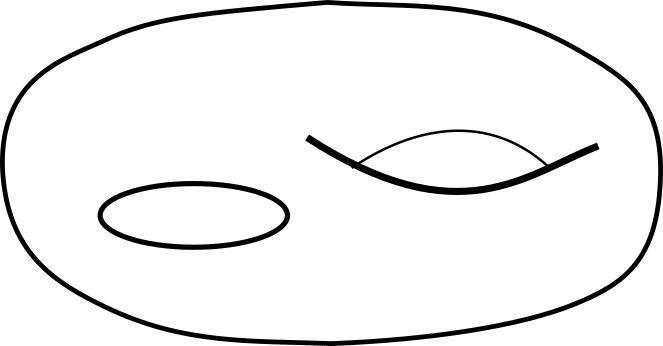} }
\put(143,29){$B$}
\put(90,67){$M$}
\end {picture}
\caption{
\label{fig_1}
{\it The inverse problem is  to reconstruct, in a stable way, the topology, the differentiable structure, and the metric of an unknown Riemannian manifold $(M,g)$, when one is given
an open ball $B=B(p,r_0)\subset M$, the eigenvalues $\lambda_j$ and the restrictions $\varphi_j|_B$ of the eigenfunctions in the ball $B$. We also study the problem of reconstructing an approximation of $(M,g)$ when only finitely many eigenvalues and eigenfunctions are given with errors.} }
\end{figure}

{\it The interior Gel'fand inverse spectral problem} is that of the reconstruction of $(M, g)$ from
its ISD. It was solved in \cite{KrKL}, \cite{KLY}. In this paper we consider the problem 
of an approximate reconstruction of $(\M,g)$ when we
know only its FISD,
namely,  the first eigenvalues,
$\la_j < \delta^{-1}$
  with some small $\delta \in (0,1)$
and the corresponding eigenfunctions of $\varphi_j|_{B(p, r_0)}$.
 Furthermore, we assume that we know all these objects with some error. However, due to the well-known ill-posedness of inverse problems, to achieve this 
goal one needs to assume that the manifold to be approximately reconstructed should lie in a properly 
chosen class of manifolds.  In this paper we concentrate on an appropriate Gromov's class of 
pointed manifolds.

{\dtextRB
Next we define a class of manifolds satisfying geometric bounds, in terms of the constants 
$R, D, i_0$, and $n$, and the radius $r_0$. Those constants have to be consider as global parameters in all calculations.
} 

\begin{definition}\label{def:2}
(Riemannian manifolds of bounded geometry). For any $n\in \Z_+$ and $R>0,$  $D>0$, $i_0>0$,
${\bM}_{\dtext n}:={\bM}_{\dtext n}(R, D, i_0 )$ consists of $n$-dimensional pointed 
compact Riemannian manifolds $(\M,g,p)$
such that
\beq \label{26.1} \nonumber
& & i)\, \sum_{j=0}^3 \|\nabla^j \hbox{Ric}(\M,g)\|_{L^\infty(M,g)}\leq R, 
\\
& & ii)\, \diam (\M,g) \leq D, \\ \nonumber
& & iii) \,\inj (\M,g) \geq i_0.
\eeq
Here $\hbox{Ric}(\M, g)=\hbox{Ric}^M_{jk}$ stands for the Ricci curvature of $M$, 
$\diam(M, g)$ for the diameter of $M$,
and $\inj (\M,g)$  for the injectivity radius of $(M,g)$. At last,
 $\nabla$ stands for the covariant derivative on $(M,g)$. {\dtext 

The norm of $\nabla^j \hbox{Ric}(\M,g)$
is computed using the metric $g$, e.g. $\|\nabla \hbox{Ric}^M\|=(g^{ii'}g^{jj'}g^{kk'}(\nabla_i \hbox{Ric}^M_{jk})(\nabla_{i'} \hbox{Ric}^M_{j'k'}))^{1/2}$.}

\end{definition}

{\ptext We recall that 
a pointed 
compact Riemannian manifold $(\M,g,p)$ consists of a manifold $M$, its Riemannian metric $g$,
and an arbitrary point $p\in M$. This definition is used as we specify the point $p$ near 
which the values of the eigenfunctions are measured.}

\begin{figure}
\begin {picture}(300,100)(35,0)

 \put(20,0){\includegraphics[height=5.0cm]{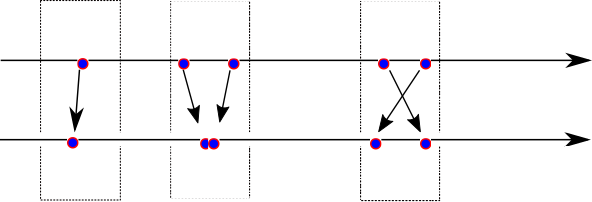}}

\put(420,112){$\R$}
\put(420,58){$\R$}

\put(73,112){$\mu^1_1$}
\put(73,23){$\mu^2_1$}
\put(146,112){$\mu^1_2$}
\put(153,23){$\mu^2_2$}
\put(180,112){$\mu^1_3$}
\put(168,23){$\mu^2_3$}

\put(286,112){$\mu^1_4$}
\put(280,23){$\mu^2_4$}

\put(312,112){$\mu^1_5$}
\put(314,23){$\mu^2_5$}

\end {picture}
\caption{
\label{fig_3A}
{\it Clusters of eigenvalues of operators: Consider a selfadjoint operator $A_1$ having a compact inverse and its perturbation $A_2=A_1+B$, where the operator norm of  the selfadjoint operator $B$ is small.
Then the eigenvalues $(\mu^1_j)_{j\in \mathbb N}$ of the operator $A_1$ and the eigenvalues $(\mu^2_j)_{j\in \mathbb N}$ of $A_2$  are $\delta$-close with some small $\delta$. Note that the eigenvalues change continuously in small perturbations, but the eigenvalues may change order, and several eigenvalues can move together forming an eigenvalue of a higher multiplicity. However, the eigenvalues can
{\dtextRB be}
grouped together to  clusters contained in separated intervals $[a_p,b_p]\subset \R$. The vector spaces spanned by the eigenvectors in such  clusters change continuously in small perturbations. This motivates Definition  \ref{def:3}. } 
}
\end{figure}
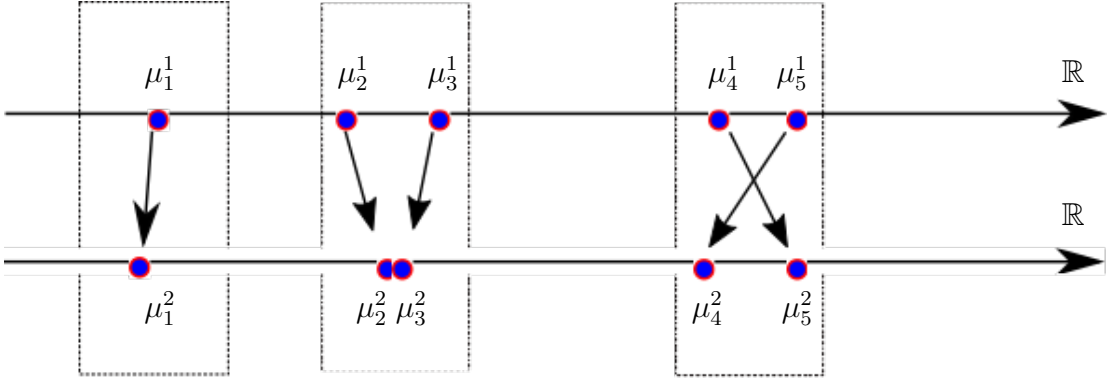

In the future, without loss of generality, we assume 
\beq \label{4.09.2017}
r_0 < \min\bigg(\frac{i_0}{2},\,\frac{\pi}{2 \sqrt{K}},1\bigg).
\eeq
 Here $K$ is the bound for
the sectional curvature on ${\bM}_{\dtext n}$. {\dtext The bound $K$ depends only on $R, D, i_0$, and $n$,
see \eqref{12.06.7}}. This makes it possible to use
in $B(p, r_0)$ the Riemannian normal coordinates which allows us to compare interior spectral data
of different manifolds in ${\bM}_{\dtext n}$.
To formalise the above, let $B(r_0) \subset \R^n$ be an Euclidian ball of radius $r_0$
  and $h$ be 
some Riemannian  coordinates in $B(r_0)$ making it a ball of radius $r_0$ with respect 
to $h$. Let
$\bD $ be a collection of elements {\dtext (Data Sequences)}
\beq\label{elements X}
DS=(\,(B(r_0),h)\,,\, \{(\mu_j,\, \psi_j|_{B(r_0)})\}_{j=0}^\infty\,)
\eeq
where $0=\mu_0 < \mu_1 \le \mu_2 \le \dots, \, \mu_j \to \infty,$ and
$\psi_j\in  L^2(B(r_0), h)$.

\begin{definition}\label{def:3}
 (Interior spectral topology.)
Let $\delta>0$. For $i=1,2,$ consider the collections $DS^i  \in \bD$.

We say that $DS^1$ and $DS^2$
are $\delta$-close
if the following is valid:
\medskip
 \noindent
There are  $P\in \Z_+$  and disjoint intervals
\beq\label{intervals}
I_p = (a_p,b_p) \subset ( -\delta,\ \delta^{-1}+\delta),
\quad p=0, 1,\dots, P,
\eeq
such that

\smallskip
\noindent
i) $b_p-a_p <\delta.$

\smallskip
\noindent
ii) For any $\mu_j^i,\, i=1,2$ with  $|\mu_j^i| < \delta^{-1}$
there is $p$ such that
$\mu_j^i \in  I_p$.

\smallskip
\noindent
iii) {\mtext For $p=0, \, n_0^i=1$. For any $p \ge 1$, the total number $n_p^i$ of elements in sets $\mathcal J^i_p=\{j\in\Z_+;\ \mu_j^i\in I_p\}$
coincide}, i.e.
$ n_p^1 = n_p^2),$ and satisfies {\dtext $n_p^1 = n_p^2\geq 1$.}

\smallskip
\noindent
iv)  {\mtext  There is an orthogonal matrix $O \in O(n)$, 
such that the  metrics $O_* h_1$ and $h_2$  are Lipschitz $\delta$-close
on $B(r_0)$, i.e., for
any {\dtext $x\in B(r_0)$ and} $\xi=(\xi^1, \dots, \xi^n) \in \R^n,\, \xi \neq 0$, {\dtext we have}
\beq \label{1.21.12}
(1+\delta)^{-1} \leq \frac{(O_* h_1)_{jk}(x)\,\xi^j \xi^k}{(h_2)_{jk}(x)\,\xi^j \xi^k} \leq 1+\delta,
\eeq

\smallskip
\noindent
v) For any $p$ there is a unitary matrix
$$
A_p =\left[a^{(p)}_{jk}  \right]_{j,k\in \mathcal J_p} \in U(n_p),
$$
such that
\beq\label{Sobolev estimate for eigenfunctions}
& &
\|A_p\,\cdotp  (O_*  \Psi_p^1)- \Psi_p^{2}\|_{\left(L^2(B(r_0),h_2)\right)^{n_p}} \leq \delta,\\
& & \label{Sobolev estimate for eigenfunctions2}
\|A_p^{-1}\,\cdotp ( (O^{-1})_* \Psi_p^2)- \Psi_p^{1}\|_{\left(L^2(B(r_0),h_1 )\right)^{n_p}} \leq \delta.
\eeq
Here,} $\Psi_p^i$ is the vector-function
$\{\psi_j\}_{j \in \mathcal J^i_p}$.
{\dtext Note that above the number $P$ indicates 
{\dtextRB
how many groups of} eigenvalues are clustered to satisfy conditions i-iv. Moreover, for two sequences $DS^1$ and $DS^2$, the above conditions i-iv may be valid
with several  different values of $P$ and intervals  $I_p,$ $p=1,2,\dots, P$.}
\end{definition}

\begin{Remark} Condition v) can be interpreted as the closedness of the Riesz projectors corresponding to $\Delta_{g_i}$ onto $I_p$.
\end{Remark}

 We note that in a more restricted context of Gelfand's inverse problem for a Schr\"odinger operator with simple spectrum in a domain in $\R^n$
 a similar topology was introduced  in
\cite{AlS}.

\subsection{The main results}

To formulate our result on an approximate reconstruction, we use
the  Gromov-Hausdorff distance.

\begin{definition}\label{def:1b} (GH-topology, see e.g. \cite{Gr}, \cite{BuBuI}).
  Let $(X^i, d^i, p^i),\,i=1,2$ be pointed
compact metric spaces. Then the  pointed Gromov-Hausdorff
distance $d_{GH}(X^1,\,X^2)$ is the infimum of
all $\e>0$  such that there is a metric space
$(Z,d_Z)$ and isometric embeddings
${\bf i}_1:X^1\to Z$ and ${\bf i}_2:X^2\to Z$ which satisfy
\ba
d_H({\bf i}_1(X^1),\,{\bf i}_2(X^2))<\e,\quad d_Z({\bf i}_1(p^1),\,{\bf i}_2(p^2))<\e.
\ea
Here $d_H$ denotes the Hausdorff distance in $Z$, see \cite{BuBuI}.
\end{definition}

The  {\dtext main} result of the paper is:

{\ktext
\begin{theorem}\label{l:10} Let $n\ge 2$, $R, D, i_0$ and $r_0$ satisfying \eqref{4.09.2017}
be given. Then there {\dtextRBX exist}  $\RobCeightyfour>1$ and $\RobCninety>0$,
{depending only on $n$, $R, D, i_0$ and $r_0$,} such that the following is true:
\\
Let $(M^{(1)},g^{(1)},
p^{(1)}),(M^{(2)},g^{(2)},
p^{(2)}) \in {\overline{ {\bf M}_{n}}}$.
Assume that  the
interior spectral data
of the operators $-\Delta_{g^{(i)}}$  on $M^{(i)}$ in the balls $B^{(i)}=B_{M^{(i)}}(p^{(i)},r_0)\subset M^{(i)}$, that is, the collections
\beq\label{e:data 1}
\big(\, (B^{(i)}, g^{(i)})\, ,\,  \{(\lambda_j^{(i)},\, \varphi^{(i)}_j|_{B^{(i)}}) ;\, 
j=0,1,2,\dots\}\big)
\eeq
are $\delta$-close, in the sense of Definition \ref{def:3},
with
{\dtextRBX $0<\delta \le \exp(-e)$. } Then
\beq \label{12.06.3}
d_{GH}((M^{(1)},p^{(1)}),(M^{(2)},p^{(2)}))\leq \RobCeightyfour \big(\ln\big(\ln\frac{1}{\delta}\big)\big)^{-\RobCninety}.
\eeq
\end{theorem}
}

{\dtext The above stability estimate is $\log-\log$  type. It is not known if this type of result is optimal, but the counterexamples of Mandache \cite{Mandache} for equivalent inverse problem show that the stability result can not be better than logarithmic.}

\medskip

The proof of Theorem \ref{l:10} is constructive, and 
is based on the following 
{\ptext result on the reconstruction of the manifold from the data.
{\wtext Below, when we state that a manifold $(M^*,g^*)$ can be constructed from the data, 
we mean 
that there is a sequence of steps, where we solve a finite number of  
quadratic minimization problems in {\zztext 
finite} dimensional spaces, choose elements from finite sets or compute certain explicit functions. 
Indeed, we do the following steps. First, we solve quadratic minimization problems in finite dimensional vector spaces 
(that are equivalent to solving linear equations) to find the finite sequences 
$(d^a_j(\a,i))_{j=0}^{\janul}\in \R^ \janul$, where $(\a,i)$  run over a finite index set, see Theorem \ref{main sliceing with errors}. Second, we use these sequences to compute approximative volumes 
$\,\vol^a(M^*_{ (i)}(\beta))$, of subsets of $M$, where $(i,\beta)$ runs over a finite index set, see Lemma \ref{volumes}.
Third, we choose the set of admissible indexes $\beta$ for which the approximative volumes are larger than a certain threshold value, see Definition \ref{admissible}. The admissible indexes are used in Section \ref{6.1}
and Lemma \ref{1st GH} to define a finite set 
of piecewise constant functions, $R^*_M$, that {\zztext approximate} the collection of the interior distance functions.
Using the finite set $R^*_M$  and a modified version of the construction given in \cite{KaKuLa} we construct 
a finite metric space $(M^*,d^*)$, that approximates the Riemannian manifold $(M,\dist_g)$  in Gromov-Hausdorff sense.}
} 

\begin{proposition}\label{l:10-A} Let $n\ge 2$, $R, D, i_0$ and $r_0$ satisfying \eqref{4.09.2017}
be given. 
Then there exists a constant
$\delta^* = \delta^*(n, R, D, i_0, r_0)$  and {\dtextRBX positive} constants $\RobCninety <1$ and $\RobCfortythree>1,$
{\Mtext depending {\dtext only on} $n$, $R, D, i_0$ and $r_0$,}
such that,
for all $\delta$ with
\beq \label{12.06.4}
0<\delta \le  \delta^*,
\eeq
the following is true:
\\
Assume that $(M,g,
p) \in {\overline{ {\bf M}_{\dtext n}}}$ and we are given a collection
\beq\label{e:data}
\big(\, (B(r_0), g^a)\, ,\,  \{(\mu_j,\, \varphi^a_j|_{B(r_0)}) ;\, 
j=0,1,2,\dots,J\}\big)
\eeq
that  is $\delta$-close, in the sense of Definition \ref{def:3},
 to interior spectral data
of the operator $-\Delta_g$  on $M$.\,
%
%
%
%
%

{\dtext Using the data (\ref{e:data})} we can construct {\dtext a pointed metric space  $(M^*, d^*,p^*)$}
such that
\beq \label{13.10.2017b}  
d_{GH}(M,\,M^*)\leq \e, \quad \hbox{where}
\quad \e=\RobCfortythree \big(\ln\big(\ln\frac{1}{\delta}\big)\big)^{-\RobCninety} .
\eeq
\end{proposition}

{
{\ptext We note that in Proposition \ref{l:10-A} that value of $J$  is not fixed,
but it just has to be so large that every $j$, for which the eigenvalue $\lambda_j$  satisfies
$\lambda_j<\delta^{-1}+\delta$, fulfils the inequality $j\leq J$, see \eqref{intervals}.
}
The relation $J$,  i.e., the number of eigenvalues, and the accuracy parameter  $\delta$
is discussed in Remark \ref{Rem: size of J} below.


\medskip

{\dtext
{\bf A note on the used constants.}
{\ktext In the main part of the  paper,} we will make frequent use of constants 
{\ptext
$c, C, C_1, C_2,$ etc. 
}These constants will depend only on the geometric bounds $n, R, D, i_0, r_0$, see Definition \ref{def:2}, but may change in their value from line to line. The constants that
depend only on the geometric bounds  $n, R, D, i_0, r_0$ will be called `uniform constants'. 
}{\dtextRB
When we define a constant for the first time, we specify whether it is uniform or not and write its further dependencies in parenthesis. For example the constant $C({s,m})$ (or $C_{s,m}$)  depends also on
$s$ and $m$. {\ktext Before Appendix we have collected a  table on the locations where the  constants $C_k$ and $c_k$ are defined. Conventions for constants in the Appendix are explained  in each subsection.}
}

\medskip

Also ${\overline{ {\bf M}_{\dtext n}}}$
is the closure of $ {\bf M}_{\dtext n}$ in the GH topology.
Parameter $K$ above and in Corollary 
\ref{Lip-stability} below is the bound for the sectional
curvature which is uniform, see (\ref{12.06.7}), on ${\overline{ {\bf M}_{\dtext n}}}$.

\begin{Remark} As shown in section 2.1
the class ${\overline{ {\bf M}_{\dtext n}}}$ is compact.
Thus, when checking  condition v) of definition \ref{def:3} it is sufficient
to use the standard $L^2-$ norm on $B(r_0)$.  
\end{Remark}
\noindent Recall that the for {\dtext pointed}
$C^1$-diffeomorphic manifolds $(M_1,p_1)$ and $(M_2,p_2)$  the Lipschitz distance 
is
\beq\label{e:Lip}
d_L((M_1,p_1),(M_2,p_2))& &=\inf_{F:M_1\to M_2}\bigg( \ln(\hbox{Lip}(F))+\ln(\hbox{Lip}(F^{-1}))+\\
& &\nonumber \hspace{2cm} +d_{M_2}(p_2,F(p_1))+d_{M_1}(p_1,F^{-1}(p_2))\bigg)
\eeq
where the infimum is taken over  bi-Lipschitz maps $F:M_1\to M_2$ and   $\hbox{Lip}(F)$ is the Lipschitz-constant of the map $F$, see \cite{GW}.
Inequality (\ref{13.10.2017b}) combined with the
sectional curvature bound (\ref{12.06.7}) and
the solution of the geometric Whitney problem \cite[Thm.\ 1, Cor.\ 1.9]{FIKLN} implies the following stable construction result for the manifold $M$ in the Lipschitz topology. 
\begin{corollary} \label{Lip-stability}
Let ${\dtext (M,g,p)} \in  {\overline{ {\bf M}_{\dtext n}}}$, $\delta>0$, and the metric space $M^*$
be as in 
{\dtextRBX Proposition \ref{l:10-A}. }
Using $M^*$ one can  construct  
a smooth pointed Riemannian manifold ${\dtext (N,g_N,p_N)}$ such that
$ 
|\hbox{Sec}(N)|\le \RobCzero K, \quad
\hbox{inj}(N)\ge \min\{ (\RobCzero K)^{-1/2}, (1-\RobCzeroone K^{1/3}\sigma_0^{2/3})\,i_0\},
$
and $M$ and $N$ are diffeomorphic. Moreover,
\ba
d_L({\dtext (M,p),(N,p_N)})\leq \RobCzero K^{1/3}\sigma_0^{2/3}, \quad \sigma_0=\RobCfortythree \big(\ln\big(\ln\frac{1}{\delta}\big)\big)^{-\RobCninety},
\ea
Here $\hbox{Sec}$ stands for the sectional curvature and $\RobCzero, \RobCzeroone$ are uniform constants.
\end{corollary}


\begin{Remark} Instead of eigenvalues and eigenfunctions 
one can deal with the heat kernels
$H_M(x, y, t)$ of {\Mtext $\partial_t-\Delta_g$,} cf. \cite{BeBeG}, \cite{KasKum}, \cite{KLY}.
Definition \ref{def:3} can be reformulated e.g. as
$\|H_{M^{(1)}}-H_{M^{(2)}} \|_{C(B(p, r_0)^2 \times (\delta, \infty))} < \delta$.
An analog of {\ktext Theorem \ref{l:10-A}} can be obtained. However, we do  not dwell on this
issue in the  paper.
\end{Remark}

To complete this section we recall that stability in the corresponding direct spectral
problem is well-known, see e.g. \cite{Ka}.
 In particular, let $M$ be a compact manifold equipped with metrics $g_\ell,\, \ell=1, 2, \dots,$ and $g_0$. Let  $a,b\not \in \sigma(-\Delta_{g_0})$. Denote by $P_\ell,\, P_0$ the spectral 
 orthoprojectors in 
 $L^2(M, g_\ell),\,L^2(M, g_0)$ on the interval $[a, b]$.
  Then it follows from Theorems  IV.3.16 and VI.5.12 of \cite{Ka} that if
$\|g_\ell-g_0\|_{L^\infty(M)} \to 0$ as $\ell \to \infty$, 
then $\|P_\ell-P_0\|_{L^2(M,g_0)\to L^2(M,g_0)}\to 0 $. This implies
that the ISD of $(M,g_\ell)$ converges to the ISD of $(M,g_0)$.

 \subsection{Earlier results and outline of the paper}

The Gel'fand inverse problem, formulated by I. M. Gel'fand in 50s \cite{Ge}, is the problem of determining 
the coefficients of a second order
elliptic differential operator in a domain $\Omega\subset \R^n$ from the boundary spectral data, that is,
the eigenvalues and the boundary values of the eigenfunction of the operator.
In the geometric Gel'fand inverse problem, a Riemannian manifold with boundary and a metric tensor on it
need  to be constructed from similar data.
 For Neumann boundary value problem for the operator $-\Delta_g$ on manifold $M$,
the boundary spectral data consists of the boundary
 $\p M$, the eigenvalues $\lambda_j$ and the boundary values of the  eigenfunction,
 $\varphi_j|_{\p M}$, $j=1,2,\dots$ 
 The uniqueness of the solution of the Gel'fand inverse problem has been considered in
\cite{Be,BeKu2,Nv1,KaKuLa,NSU}.
\\
To consider the formulation of the stability of the inverse problems, 
let us consider first the Gel'fand inverse on 
 a bounded domain  $\Omega\subset \R^2$ with smooth boundary $\p \Omega$ and a conformally Euclidian metric 
 $g_{jk}(x)=\rho(x)^{-2}\delta_{jk}$. Here, $\rho(x)>0$ is a smooth real valued function. Then the problem has the form
\beq\label{Eq 2}
& &
-\sum_{k=1}^2\rho(x)\left(\frac {\p}{\p x^k}\right)^2\varphi_j(x)-\lambda_j\varphi_j(x)=0,\quad \hbox{in\ }\Omega,
\quad	\partial_\nu \varphi_j|_{\p \Omega}=0.
\eeq
The problem of determining $\rho(x)$ from the boundary spectral data is ill-posed  in sense of Hadamard:
The map from the boundary data
to the coefficient $\rho(x)$ is not continuous so that small change in the data can lead to huge errors in the 
reconstructed function $\rho(x)$.
One way out of this fundamental difficulty 
is to assume a priori higher regularity of coefficients, that is
a widely used trend in inverse problems for isotropic equations,
like (\ref{Eq 2}). This type of results is called 
 conditional
stability results (see e.g.
\cite{Al,AlS,StU}).

For inverse problems for general
metric  this approach bears significant difficulties.
The reason is that the usual $C^k$ norm bounds of coefficients
are not invariant and thus this condition does not suit
 the invariance of the problem  with respect to
 diffeomorphisms. Moreover, if 
the structure of the manifold $M$
is not known a priori, the traditional approach can not be used.
The way to overcome these difficulties
is to impose a priori 
{\dtextRB constraints} in an invariant form and consider
a class of manifolds that satisfy invariant a priori bounds,  for instance
for  curvature, second fundamental form, radii
of injectivity, etc. Under such kind of conditions,  invariant stability results for various inverse problems  
{\dtextRB
have been proven 
}
in \cite{AKKLT,FIKLN,StU}.
In particular, for the Gel'fand inverse problem for manifolds with non-trivial topology, an abstract, i.e., a non-quantitative 
stability 
result  
was proven in \cite{AKKLT}. There, it was shown that the convergence of the boundary spectral data implies
the convergence of the manifolds
with respect to the Gromov-Hausdorff convergence. However, this result was based on compactness arguments
and it did not provide any
{\dtextRB estimate }. In this paper our aim is 
{\dtextRB
to improve} 
this result and to give explicit estimates for an analogous inverse problem.
\\
In this paper we consider a Gel'fand inverse problem for manifolds without boundary. Then, as explained above, instead of assuming that the boundary and the boundary values of 
the eigenfunctions are known we assume that we are given a small open ball $B\subset M$ and
the eigenfunctions $\varphi_j$ are known on this set.
Similar type of formulation of the problem with measurements on open
sets have been considered in \cite{dH1,dH2}.
We show that the Interior Spectral Data (ISD), that is, an open set $B\subset M$, the eigenvalues $\lambda_j$  and the restrictions of the eigenfunctions  $\varphi_j|_B$ 
determine the whole manifold $(M,g)$ in stable way. Also, we quantify this stability by giving explicit
inequalities under a priori assumptions on the geometry of $M$.
 We emphasise that we assume that the eigenfunctions are known only on an open subset $B$ of $M$ that may
{\dtextRB be chosen 
}
to be arbitrarily small but still e.g.\ the topology of $M$ is determined in a stable way. We note that this paper is a slightly extended and polished version of our preprint in Arxiv, published on Feb. 25, 2017. We note that in spectral geometry one has studied similar stability problems where the heat kernel are known on the whole
 manifold, \cite{BeBeG,KasKum,KasKum2}. This data is equivalent to knowing the eigenvalues and the eigenfunctions and the eigenfunctions on the whole manifold.
\\
Outline of the paper: Ch. 2 introduces the geometric set-up.
Ch. 3 formulates the stability of the unique continuation for the solution of the wave equation together with Corollary \ref{main}  for its spatial projection $v$.
Ch. 4 presents Thorems \ref{main sliceing} and \ref{main sliceing with errors} proving the construction of the approximate Fourier coefficients of $\chi_{\Omega}v$ in the case of respectively exact and approximate FISD. 
Ch. 5 shows the related approximate interior distance functions. Ch. 6 collects all the previous inequalities to prove Theorem \ref{l:10}
and Proposition \ref{l:10-A}.

\section{Geometric preliminaries}

\subsection{Properties of the manifolds of bounded geometry }

Here we list some results on 
 the class
${\bM}_{\dtext n}(R, D, i_0)$,
These results can be found in
or immediately follow from
\cite{And,Ch} with further improvements in \cite{AKKLT}.
Namely, the class ${\bM}_{\dtext n}$ is precompact in GH-topology. Its closure,
$\overline{\bM_{\dtext n}}$ consists of pointed Riemannian manifolds $(\M, g,p)$
with $g \in C^5_*(M)$ which satisfy (\ref{26.1}). Here and later $^*$ indicates the Zygmund space.
\\
We define the norm of the space $C^k(M)$ invariantly by
\beq \label{norms}
\|f\|_{C^k(M)}:=\sum_{j=0}^k  \max_{x\in M} \|\nabla^j f(x)\|_g,
\eeq
where {\dtext the norm is computed using the metric $g$.} Next, for $k \in \Z_+, \, \beta \in (0, 1],$ {\dtext we use the the Zygmund spaces}
$$C^{k+\beta}_*(M)=[C^{k_1}(M),C^{k_2}(M)]_\theta, \quad 
k+\beta=\theta k_1+(1-\theta) k_2\in \R_+, \quad
\theta\in (0,1).
$$
Here $[\cdot, \cdot]_\theta$ stands for the interpolation, see e.g \cite{BL}.
Note that, for $\beta \in (0, 1),$ the H\"older spaces 
{\dtextRB fulfill} $C^{k,\beta}(M)=C^{k+\beta}_*(M)$.
\\
To achieve the $C^k_*-$smoothness of $g$, one needs some special coordinates, e.g. harmonic coordinates. For any 
number $Q>1$, that we below choose to be $Q=2$, there is a constant ${\dtext r^{(har)}}$ {\ctext depending only on 
 $n, R, D, i_0, r_0$ and $Q$,} such that, for any
${\dtext (M,g,p)} \in   \overline{\bM_{n}},\, q \in M$, there are {\dtext $Q-$harmonic} coordinates in $B(q, {\dtext r^{(har)}})$, that we denote by $Y:B(q, {\dtext r^{(har)}})\to \R^n$, $U_q=Y(B(q, {\dtext r^{(har)}}))$
that we denote by $y$.
{\dtext For $Q=2$, in these coordinates  the metric tensor $g^{(har)}_{jk}(x)=(Y^{-1})^g$  satisfies
\beq \label{26.2A}
& & 2^{-1}I \leq (g^{(har)}_{jk}(y))_{j,k=1}^n \leq 2 I,  \quad \hbox{for }y\in U_q=Y(B(q, {\dtext r^{(har)}}))\\
& &
\|g^{(har)}_{jk}\|_{C^5_*(\overline U_q)} \leq C^{(har)},\nonumber
\eeq
with some uniform constant $C^{(har)}$, see \cite{And,Ch} and \cite{AKKLT}.
We note that the existence of the harmonic radius $r^{(har)}$
and the constant $C^{(har)}$ for which (\ref{26.2A}) holds for all ${\dtext (M,g,p)} \in   \overline{\bM_{n}},$ $ q \in M$,
is based on compactness results, and therefore the dependency of $r^{(har)}$
and $C^{(har)}$ on $n, R, D, i_0, r_0$ is not explicit.}

Sometimes, with a slight abuse of notation we identify $y\in U_q$ with the corresponding point in $Y(y)\in B(q, {\dtext r^{(har)}})$.}
\\
The inequality  (\ref{26.2A}) immediately implies that {\dtext the sectional curvature $\hbox{Sec}$
and the Riemannian curvature tensor $R_M$ satisfies
\beq \label{12.06.7}
 |\hbox{Sec}(M)|\leq K,\quad \|R_M\|\leq K,\quad \|\nabla R_M\|\leq K,
\eeq
where} $K$  is a uniform constant.

For the sake of simplicity, we will work with H\"older rather then Zygmund
spaces.
It follows from
\cite{And,Ch}, with the terminology described in \cite[Sec.\ 10.3.2] {Pe},  that
 when  
${\dtext (M_k,g_k,p_k) \to (M,g,p)}$ in the GH topology on $ \overline{\bM_{\dtext n}}$,
then, for all $\beta\in (0,1),$ there are $C^{5,\beta}$-smooth diffeomorphism $F_k:M_k\to M$  
{\dtextRB such that}
\beq \label{26.4B harm}
F_*(g_k) \to g \quad \hbox{in $C^{4, \beta}(M)$, as $k\to \infty$}.
\eeq
Thus,  for any $\e>0, \beta<1$, there is
$\sigma=\sigma(\e, \beta)$ {\dtext such that we have the following: For all $M_1,M_2\in   \overline{\bM_{\dtext n}}$
such that $d_{GH}(M^1, M^2) <\sigma$,} there is a
diffeomorphism $F: M^1 \to M^2$ and
\beq \label{26.4}
\|g^h_1-F_*(g^h_2)\|_{C^{4, \beta}(M^i)} <\e,\quad i=1, 2,
\eeq
 cf.
 \cite[Sec.\ 10.3.2]{Pe}.
Returning to (\ref{26.4B harm}), for large $k$, $M_k$ and $M$ are diffeomorphic,  so that it
is possible to use results from \cite{Ka}, see the end of sec. 1.2. This implies stability of the 
direct problem
in the GH topology on $\overline{\bM_{\dtext n}}$.
\\
{\dtext Note that we can solve the ordinary differential equations that define the geodesics  in the harmonic coordinates.
Then it follows from (\ref{26.2A}) that
there is a uniform constant $\RobCone>1$, such that for  any ball $B(x,r)\subset M$,  where $(M,g,p) \in \overline{\bM_{\dtext n}}$, we have}
\beq \label{12.06.13}
\frac 1{\RobCone} r^n \le \hbox{vol}(B(x, r)) \le \RobCone r^n, \quad 0 \leq r \leq D.
\eeq
{\dtext  Thus,  the volume of balls having radius $i_0/2$ is bounded below by a uniform constant $v_0$. 
Furthermore,
by \cite{GW}, the class of Riemannian manifolds $(M,g)$ that satisfy (\ref{12.06.7}) 
and conditions $\diam (\M,g) \leq D$  and $\vol (\M,g) \geq v_0$ are pre-compact with respect to 
the Lipschitz distance $d_L((M_1,p_1),(M_2,p_2))$, see \eqref{e:Lip} and the closure of this class consists of $C^\infty$-smooth manifolds
with $C^{1,\alpha}$-metric. This implies that there is a uniform constant $C^{(Lip)}$ such that for all 
${\dtext (M_1,g_1,p_1)},{\dtext (M_2,g_2,p_2)}  \in   \overline{\bM_{n}}$ we have
\beq\label{uniform Lip}
d_L((M_1,p_1),(M_2,p_2))\leq C^{(Lip)}.
\eeq
Moreover, by \cite{Katsuda}, we have that for any $\e>0$ there is 
}{\dtextRB
$\zeta(\e)>0$,
such that for $(M_1,g_1,p_1),(M_2,g_2,p_2)  \in   \overline{\bM_{n}}$ we have
\beq\label{eq: Katsuda}
\hbox{if $d_{GH}((M_1,p_1),(M_2,p_2))<\zeta(\e)$ then $d_L((M_1,p_1),(M_2,p_2))<\e.$}
\eeq
}
We turn now to the spectral properties on ${\dtext (M,g,p)} \in  \overline{\bM_{\dtext n}}$.
{\dtext By \cite{Davies}, the inequality \eqref{uniform Lip} implies 
{\dtextRB that} the $j$-th eigenvalue
$\lambda_j(M_i,g_i)$ of the Laplacian on the manifold $(M_i,g_I)$ 
{\dtextRB
satisfies
}
$$
e^{-(n+2)C^{(Lip)}/2}\lambda_j(M_1,g_1)
\leq \lambda_j(M_2,g_2)\leq e^{(n+2)C^{(Lip)}/2}\lambda_j(M_1,g_1)
$$
for all 
${\dtext (M_1,g_1,p_1)},{\dtext (M_2,g_2,p_2)} \in   \overline{\bM_{n}}$.
%
{\dtextRB
Since
} 
the eigenvalues of the manifold  $(M_1,g_1)$  satisfy the Weyl's asymptotics 
$ \la_j(M_1)= c_{{}_{M_1}} j^{2/n}(1+o(1))$ as $j\to \infty$,} then there exists  a uniform constant $\RobCthree>1$ such that, 
\beq \label{26.5} 
\frac 1{\RobCthree} j^{2/n} \leq \la_j(M) \leq \RobCthree j^{2/n}, \quad j\in \Z_+ \hbox{ for all }
{\dtext (M,g,p)} \in  \overline{\bM_{\dtext n}(R, D, i_0)}.
\eeq 
Note  that (\ref{26.5}) is valid under a  weaker assumption
that $\hbox{Ric}(M)$ is bounded from below, see
\cite{BeBeG}.

{\Mtext 
\begin{Remark}\label{Rem: size of J}
Assume that the collection of
  $g^a|_{B_e(r_0)}$ and $((\lambda_j^{a},\varphi_j^{a} |_{B_e( r_0)}))_{j=0}^J$ 
 is $\delta$-close to the FISD
  $g|_{B_e(r_0)}$ and $(\lambda_j,\varphi_j |_{B_e( r_0)})_{j=0}^J$ of the manifold ${\dtext (M,g,p)}\in {\bM_{\dtext n}}$.
Then all intervals $I_p= (a_p,b_p)$, $p=0,1,\dots,P$   in (\ref {intervals}) satisfy $b_j\leq   \delta^{-1}+\delta$,
and thus the index $j$ of any eigenvalue $\lambda_j$ that is in some of these intervals satisfies
by (\ref{26.5}) the inequality $ \RobCthree^{-1} j^{2/n} \leq  \delta^{-1}+\delta\leq  2\delta^{-1}$.
On the other hand, if $j < ( \RobCthree^{-1}  \delta^{-1})^{n/2}$, then
$\lambda_j<\delta^{-1}.$
 Thus, without loss of generality, 
we can always  {\novtekst assume that the}
 value of  $J$ in 
Proposition \ref{l:10-A} 
  satisfies
\beq\label{bound for J}
( \RobCthree^{-1}  \delta^{-1})^{n/2}\leq J\leq   (2\RobCthree\delta^{-1})^{n/2}.
\eeq
\end{Remark}

\begin{Remark}\label{Rem: j and k ratio}
Below we will assume that $\delta < (3\RobCthree)^{-1}$.
Then for $j\geq 1$ we have $\lambda_j\geq \RobCthree^{-1}$ and 
$\lambda_j>3\delta$.
Next, assume that $\lambda_j$  and $\lambda_k$ with $k>j\geq 1$ belong in the same interval $I_p= (a_p,b_p)$
with $b_p-a_b<\delta$.
Since $\lambda_j\ge \RobCthree^{-1} >3\delta$, we have $a_p>2\delta$ so that
$b_p<2a_p$.
Then by (\ref{26.5}) we have
$$
\RobCthree^{-1} k^{2/n}\leq \lambda_k\le 
b_p\leq 2a_p\leq 2\lambda_j\leq 2\RobCthree j^{2/n},
$$
implying
\beq\label{j and k ratio}
j< k\leq (2^{1/2}\RobCthree)^{n} j.
\eeq
\end{Remark}}
\noindent
Next, instead of harmonic coordinates, we can use coordinates made of the eigenfunctions
 $\varphi_j$. It turns out, cf. \cite{BeKu,AKKLT}, that in a neighbourhood of any $x \in M$
 there are $\varphi_{j(1; x}), \dots, \varphi_{j(n; x)}$ which form $C^6_*$-smooth coordinates.
 Moreover, by the compactness arguments, there are uniform {\dtext constants $r$ and $C$ so that 
 these coordinates are well defined in any ball $B(x, r)\subset M$, where $(M,g,p) \in \overline{\bM_{\dtext n}}$, and
the metric tensor $g$ in these coordinates} satisfies (\ref{26.2A}).
There is also a uniform number $N \in \Z_+,$ such that  {\novtekst we can take} $j(\ell; x) \leq N,\, \ell=1, \dots, n.$
\\
Next, using  $((\la_j, \varphi_j))_{j=0}^\infty$,
we introduce the Sobolev spaces
$H^s(M), \, s \in \R$,
\beq \label{26.6}
f(x) =\sum_{j=0}^\infty f_j \varphi_j(x) \in H^s(M)\quad \hbox{iff}\quad\|f\|^2_{H^s}:=\sum_{j=0}^\infty \bra \la_j \cet^s |f_j|^2 <\infty,
\eeq
where $\bra \la \cet = (1+\la^2)^{1/2}$.

\subsection{Distance coordinates}
Recall that there are harmonic coordinates in $B(x,{\dtext r^{(har)}})$ ball near any $x \in M
 \in \overline{\bM_{\dtext n}}$,  see (\ref{26.2A}). In the Proposition below we use
such coordinates as background coordinates near $x$.

{\dtext Below, we say that a subset $Y\subset X$ is a $\tau$-net in 
{the metric space $X$
if
the  union of the balls $B_X(y,\tau)$, $y\in Y$, contains} the whole space $X$.
Also, we say that $Z\subset X$ is $\tau$-separated, if for all $z_1,z_2\in Z$,
$z_1\not =z_2$  we have $d_X(z_1,z_2)\ge \tau$. Observe that if $Z\subset
X$  is a maximal $\tau$-separated subset of $X$ (maximal in the sense that any
other  $\tau$-separated subset of $X$ that contains $Z$  has to be equal to $Z$)
{\dtextRB 
, then it is} a $\tau$-net in $X$.}

\begin{figure}
\begin {picture}(200,100)(40,0)
\put(100,0) {\includegraphics[height=5.0cm]{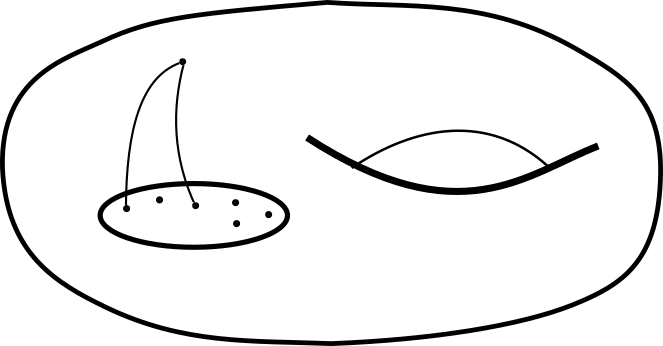} }
\put(173,49){$z_{j_\ell}$}
\put(183,110){$y$}
\end {picture}
\caption{
\label{fig_3B}
{\it There are $L-1$ points $z_1,z_2,\dots,z_{L-1}\in B(p,r_0/4)$ such that for any $y\in M$ there are  $n$ points $z_{j_1(y)},\dots,z_{j_n(y)}$ so that 
the distance functions $X^\ell(x)=d_M(x,z_{j_\ell(y)})$, $\ell=1,2,\dots,n$ define local smooth coordinates $x\mapsto (X^\ell(x))_{\ell=1}^n$ in a neighbourhood $U_y\subset M$ of the point  $y$.
}
}
\end{figure}

%

\begin{proposition} \label{uniform_covering} $ $\\
There are uniform {\dtext constants} $\tau_0,\rho_0<\min\{{\dtext r^{(har)}}/4, r_0/128\}$ and uniform constants
 {\dtext $L\in \Z_+$  and} $\RobCfive, \RobCsix,\,\RobCseven>0$ 
{\ztext depending only on $n$, $R, D, i_0$ and $r_0$,}
 such that,
for any  ${\dtext (M,g,p)}  \in \overline{\bM_{\dtext n}(R, D, i_0)}$ 
the following holds true:
There is a $\tau_0$-net in $B(p, r_0/4)$ with at most $L-1$ points.
 Let  $\{z_1, \dots, z_{L-1}\}\subset  B(p, r_0/4) $  be 
 {\ptext an arbitrary collection of points that is a} ${\dtext \tau_0}$-net in $ B(p, r_0/4)$.
Then,
\medskip

(i) For all $x \in M$, there are $n$  points $z_{j(i)}\in Z$, 
$j(i)=j(i; x)$, $i=1,2,\dots,n$  such that
the map $X:B(x,\rho_0)\to \R^n$,
\beq \label{5.09.2017}
X:y =(y^1, \dots, y^n)\mapsto (d_M(y, z_{j(1)}),d_M(y, z_{j(2)}), \dots, d_M(y, z_{j(n)})),
\eeq
for coordinates where $X:B(x,\rho_0)\to X(B(x,\rho_0))$ is a Lipschitz-smooth diffeomorphism 
and
 \beq\label{tau coordinates final}
 \|DX\|_{L^\infty(B(x,\rho_0))}+ \|DX^{-1}\|_{L^\infty(X(B(x,\rho_0))))}\leq \RobCfive.
 \eeq
where {\dtext the norms are computed using the metric $g$ on $M$ and the Euclidean norm in $\R^n$.} 
Moreover, $z_{j(i)}$ can be chosen so that $d(x, z_{j(i)}) > r_0/{\dtext 16}$ and
the metric tensor $(g_{ij})_{i,j=1}^n=X_*g$ in these coordinates satisfies
\beq \label{distance_co}
& &\RobCfive^{-1} I \le (g_{ij}(z))_{i,j=1}^n\le \RobCfive I, \,\,z\in X(B(x, \rho_0)).
\eeq

(ii) The map $H:M\to \R^{L-1}$, defined by $H^L(x)=(d_M(x,z_j))_{j=1}^{L-1}$, satisfies
  \beq\label{tau coordinates}
\frac 1{\RobCsix}\leq \frac {|H^L(x)-H^L(y)|}{d(x,y)}\leq \RobCsix,\quad
\hbox{for\,all}\,\, x,y\in M,\,x\not =y.
 \eeq
 {\novtekst Here as a norm in (\ref{tau coordinates}) we can take e.g. the Euclidian norm in
$\R^{L-1}$}.
 
\end{proposition}
%

\noindent
{\bf Proof.}
 Let us first consider one pointed manifold ${\dtext (M,g,p)} \in \overline{\bM_{\dtext n}}$. 
{\dtext  Let us consider 
the extended exponential map
 $$
F:TM\to M\times M,\quad F(x,\xi)= (x,\exp_x(\xi)).$$ 
Inequalities  \eqref {26.1}  and \eqref{12.06.7} imply
that in the set $S=\{(x,\xi)\in {TM}:\ \|\xi\|\leq 2D\}$ the map $F$ is $C^2$-smooth and its norm
 in $C^2(\overline S)$ is bounded by a uniform constant. The proof this is analogous to that of Lemma 2 in
 \cite{KatsudaKuLa}. 
Let
  $x_0 \in M$ and $\gamma_{x_0,\xi_0}([0,s_0])$, $\xi_0\in S_{x_0}M$  be a shortest geodesic
  from $x_0$ to $p$ where $s_0=d(x_0,p)$. When $s_0\ge r_0/2$, choose $s_1=s_0- {\fftext r_0/5}$,
  and when  $s_0<  r_0/2$,  choose $s_1=s_0+ {\fftext r_0/5}$. Then the  point $p_1=\gamma_{x_0,\xi_0}(s_1)$
  satisfies 
  $p_1 \in \p B(p, {\fftext r_0/5}) $,
and $d(x_0, p_1) {\fftext \ge r_0/5} $. As $r_0<i_0$, we see that the geodesic $\gamma_{x_0,\xi_0}([0,s_1+{\fftext r_0/5}])$
is a length minimising  geodesic  between its endpoints. In particular, this implies that 
$\gamma_{x_0,\xi_0}([0,s_1])$ 
{\dtextRB 
continues behind
$p_1$ as a shortest   curve between  its points. 
} 
As in \cite[Lemma 4]{KatsudaKuLa},
(see also \cite{Ivanov} where related results are proven with lower regularity assumptions), we see that \eqref {26.1}  and \eqref{12.06.7} imply
that there is a uniform constant ${\fftext r_*\in (0,\frac {r_0}{100})}$ 
 such that following is true.  Let $\mathbb B=B_{TM}((x_0,s_1\xi_0),r_*)$ be the ball of radius $r_*$ and center $(x_0,s_1\xi_0)$ defined in the tangent bundle 
 $(TM,g)$ using the Sasaki metric. Then for vectors $(x,s\xi)\in \mathbb B$, where  $\xi\in S_xM$, $s>0$, the geodesics $\gamma_{x,\xi}([0,s])$
 are length minimizing curves between their end points. Moreover,
 the exponential map in $\mathbb B$, that is,
 $$
 F:\mathbb B\to F(\mathbb B)
 $$
 is a diffeomorphism and satisfies $\|dF\|\leq C$ in $\mathbb B$, and  $\|(dF)^{-1}\|\leq C$ in  $F(\mathbb B),$  and $F(\mathbb B) {\fftext \subset M\times B(p_1, r_0/100)}$. {\fftext Here,
$M\times B(p_1, r_0/100)
 \subset M\times B(p, r_0/4)$.}
 In particular, $d(x,z)=|F^{-1}(x,z)|$ for $(x,z)\in F(\mathbb B)$.
 Let $\xi_j\in S_{x_0}M$ and $t_j>0$, $j=1,2,\dots,n$ be such that 
 \beq\label{conditions for tj}
 \|\xi_j-\xi_0\|<r_*/s_1,\quad \|\xi_j-\xi_k\|>r_*/(8s_1)\hbox{ for $j\not =k$, and } |t_j-s_1|<r_*/8.
 \eeq
 Then $s_1\xi_j\in \mathbb B$.
 Let $z_j=\exp_{x_0}(t_j\xi_j)\in B(p, r_0/4).$
 We see that  $\nabla d_M(\,\cdotp,z_j)|_{x_0}=-\xi_j$
 and  $\|(dF|_{x_0})^{-1} \|\leq C.$

Inverse function theorem, see e.g. \cite{H1}, and the facts that $\|(dF|_{x_0})^{-1} \|\leq C$
and that $F$  has a uniformly bounded $C^2$-norm in $S$,
 imply for the map  $$H^M_{z_1,\dots,z_n}(x)=(d_M(x,z_j))_{j=1}^n$$
 that there are uniform constants $\rho_*>0$ and $c_*>0$ such that we have
 \beq\label{eq H-map}
 |H^M_{z_1,\dots,z_n}(x)-H^M_{z_1,\dots,z_n}(x')|\geq c_*d_M(x,x'),\quad
 \hbox{for all } x,x'\in B_M(x_0,\rho_*).
 \eeq
 Let us now choose $\xi^0_j\in S_{x_0}M$, $j=1,\dots,n$ 
 such that $s_1\xi_j^0\in \mathbb B$ satisfy
 \ba
 \|\xi^0_j-\xi_0\|<c_*/(2s_1),
 \quad \|\xi^0_j-\xi^0_k\|>c_*/(4s_1)\hbox{  for $j\not =k$.}
 \ea
  Let 
  $z_j^0=\exp_{x_0}(s_1\xi^0_j)$. As 
  $\|dF^{-1}\|\leq C$ in $F(\mathbb B)$, there is a uniform constant 
  $\fftext{\tau_*\in (0,r_0/100)}$
 such that if $\tilde  z_j\in B(p, r_0/4)$, $j=1,\dots,n$ satisfy $d(\tilde  z_j,z_j^0)<\tau_*$ then 
 there are $\tilde  \xi_j\in S_{x_0}M$ and $\tilde t_j>0$ such 
 that $\tilde  z_j=\exp_{x_0}(\tilde  t_j\tilde  \xi_j)$ and
 $$
 |\tilde  \xi_j-\xi_j|<r_*/(8s_1)\quad\hbox{and}\quad |\tilde  t_j-s_1|<r_*/8.
 $$ 
 Then $\tilde  \xi_j$ and $\tilde  t_j$  satisfy \eqref{conditions for tj}.
 Thus,  (\ref{eq H-map}) implies that the map $H^M_{\tilde  z_1,\dots,\tilde  z_n}(x)$
 satisfies \eqref{eq H-map}.
This implies that if $\{\hat z_i\in B_M(p, r_0/4)$, $i=1,2,\dots,i_M\}$ is any $\tau_*-$net in $B_M(p, r_0/4)$ then for all $j=1,2,\dots,n$ there are $i_j\in \{1,2,\dots,i_M\}$ such that $d_M(\hat z_{i_j},z^0_{j})\leq \tau_*.$
Then the above implies that 
 for $H^M_{\hat z_{i_1},\dots,\hat z_{i_n}}(x)$
 satisfies \eqref{eq H-map}.

 {\fftext Observe that above $(x,s\xi)\in \mathbb B$, so that 
  $d_M(x,x_0)<r_*<r_0/100$. Moreover, $d(x_0, p_1) \ge r_0/5$, 
  $z_j^0\in B(p_1, r_0/100)$,
  $d_M(\hat z_{i_j},z^0_{j})\leq \tau_* <r_0/100$
   yield $d_M(x,\hat z_{i_j})
 \ge  r_0/5-3r_0/100
 >r_0/8.$}

Note that above $c_*$, $\tau_*$ and $\rho_*$ are uniform {\fftext constants} and the estimate 
 (\ref{eq H-map}) is valid for some  {\fftext points
 $\hat z_{i_j}$ in any $\tau_*-$net $\hat z_{i}$ in $B(p, {r_0}/4)$, that satisfy $d_M(x,\hat z_{i_j}) >r_0/8$,}
  and any  ${\dtext (M,g,p)} \in \overline{\bM_{\dtext n}}$.

This proves \eqref{tau coordinates final} and (\ref{distance_co}) in claim (i).

Next we consider the claim (ii).
Let us show that there are  $h_1>0$ and $\tau_1>0$ such that for any
 $(M,g,p)\in  \overline{\bM_{\dtext n}}$  and any maximal $\tau_1$-separeted set 
 $\{z_1,\dots,z_{L-1}\}\subset B(g,r_0/4)$ we have
 \beq\label{claim to contradiction}
 \sup_{x,y\in M,\ x\not =y}\bigg(\sup_{j_1,\dots,j_n}\frac {|(d_M(x,z_{j_i}))_{j=1}^n-(d_M(y,z_{j_i}))_{j=1}^n|_{\R^n}}{d_{M}(x,y)}\bigg)\ge h_1,
 \eeq
where the supremum is taken over all $1\leq j_1<j_2<\dots<j_n\leq L-1$.

Assume the opposite. Then for all $k\in \Z_+$ there are $h_k>0$, $(M_k,g_k,p_k)\in  \overline{\bM_{\dtext n}}$ and $1/k$-nets
$\{z_j^k: j=1,2,\dots,L_k\}\subset B(p_k,r_0/4)$ and points $x_k,y_k\in M_k$, $x_k \neq y_k$  
so that $h_k\to 0$ and 
 \beq\label{limiting sequence}
 \sup_{x,y\in M_k,\ x\not =y}\bigg(\sup_{j_1,\dots,j_n}\frac {|(d_{M_k}(x_k,z^k_{j_i}))_{i=1}^n-(d_{M_k}(y_k,z^k_{j_i}))_{i=1}^n|_{\R^n}}{d_{M_k}(x_k,y_k)}\bigg)<h_k.
 \eeq
Using compactness arguments for $\overline{\bM_{\dtext n}}$ and choosing 
a suitable subsequence of the manifolds $(M_k,g_k,p_k)$ we can assume that 
$(M_k,g_k,p_k) \to (M,g,p)$ in the  Lipschitz-topology. Then  there are diffemorphisms
 $F_k:M_k \to M$ such that $F_k(p_k)\to p$ and Lip$(F_k)\to 1$ and Lip$(F_k^{-1})\to 1$.
Moreover, {\dtextRB we} can assume that   $F_k(x_k)\to x$ and $
F_k(y_k) \to y$  in $M$ and, after using the Cantor diagonalization procedure, we can assume that there are limits $\lim_{k\to \infty} F_k(z_j^k)=  z_j$ in $M$, for all 
$j=1,2,  \dots$. 
Next,  using (\ref{26.4}), we see that $d_{M_k}(x_k,y_k) \to d_M(x, y)$,
$d_{M_k}(x_k, z_j^k) \to d_M(x, z_j)$ and $d_{M_k}(y_k, z_j^k) \to d_M(y, z_j)$. Also
$\{z_j\}_{j=1}^\infty$ is dense in $B_M(p, r_0/4)$. Therefore, $d_M(x, z)=d_m(y, z)$ for
all $z \in B_M(p, r_0/4)$.
Then {\dtext  \cite[Lemma 13]{Helin}} (see also \cite[Lemma 3.30]{KaKuLa}), 
implies that $x=y$. 

Let $k$  be so large that $1/k<\tau_*/2$, and
\ba
& &d_M(F_k(p_k),p)<\tau_*/2,\quad \hbox{Lip}(F_k)\leq 2, \quad \hbox{Lip}
(F_k^{-1})\leq 2,
\\
& &
d_{M}(F_k(x_k), x)<\rho_*/4,\quad
d_{M}(F_k(y_k), y)<\rho_*/4,\quad h_k<c_*.
\ea
 As $x=y$, these imply
$d_{M}(F_k(x_k), F_k(y_k))<\rho_*/2$ 
{\dtextRB and hence}
$d_{M_k}(x_k,y_k)<\rho_*$.
As $1/k<\tau_*/2$, 
the points $z^k_{j},$ $j=1,\dots,L-1$ form a $\tau_*$-net in $B_{M_k}(p_k,r_0/4)$.
Then the inequality \eqref{limiting sequence} for $x_k$ and $y_k$, 
with $h_k<c_*$, is in contradiction with 
the fact that   there is a subset of $n$ of the points in  $\tau_*$-net  
 $z^k_{j},$ $j=1,\dots,L_k$ for which 
\eqref{eq H-map} holds. 
{\dtextRB This proves }
\eqref{claim to contradiction} with some uniform constants $\tau_1$  and $h_1$.

 We observe that a maximal $\tau_1$-separated subset in  the ball $B(p,r_0/4)$
has at most $C_*=\vol(B(p,r_0/4))/\vol_{n,R}(B(x,\tau_1))$ points,
where $\vol_{n,R}(B(x,\tau_1))$  is the volume of the ball of radius $\tau_1$
on the $n$-dimensional sphere having  constant curvature $R$. Hence 
we see that the number  of
{\dtextRB points for a} maximal $\tau_1$-separated subset in $B(p,r_0/4)$ is bounded by a uniform constant  $C_*$. Thus we can choose $L$ to be the integer part of $C_*$ and $\tau_0=\min(\tau_*,\tau_1)$, which makes $L$ and $\tau_0$ uniform constants.

As the number $L-1$  of points in the  $\tau_1$-nets we consider is bounded by a uniform constant, we see that
\eqref{tau coordinates} is valid with $\RobCfive=h_1^{-1}+L$.
These prove the claims (i) and (ii).}

\hfill\Box\medskip

The above considerations 
bring about the following result.


\begin{lemma}  \label{ga-separation} {\Mtext There exist uniform constant $ \MattiCseven>0$ and 
{\dtext uniform constant $N_F\in \mathbb Z_+$} {\dtext (that is, 
{\ztext $\MattiCseven$ and the integer $N_F$ depend only on  $n, R, D, i_0, r_0$)}
}
 such that
\\
\noindent (i)  Let $\sigma \in (0,\tau_0]$. Then {\dtextRB any maximal  $\sigma$-separated set
 $x_1, \dots, x_{N(\sigma)}$ in $M$
 is such that the number of its elements fulfills the bound
}
\beq \label{19.07.1}
 N(\sigma) \le \tilde N(\sigma)= \MattiCseven \sigma^{-n}.
\eeq
Moreover, {\dtextRB the balls $B(x_k, 4 \sigma)$  satisfy the finite intersection property with at most  $N_F$ intersections, that is, any point $x\in M$  belongs to at most $N_F$
balls $B(x_k, 4 \sigma)$.}
\\
\noindent (ii) Let  $\sigma \in (0,\tau_0]$. Then {\dtextRB any maximal  $\sigma$-separated set
  $z_1, \dots, z_{N_1(\sigma)}$
 in $B(p, r_0/4)$ 
  is such that the number of its elements fulfills the bound
}
$N_1(\sigma) \le \tilde N(\sigma)$,
and the balls  $B(z_k, 4 \sigma)$  {\dtext satisfy the finite intersection property with at most  $N_F$ intersections.}}
\end{lemma}

\noindent{\bf Proof.}  It remains to prove the finite intersection property. It follows from (\ref{12.06.13})
if we take into the account that
$B(x_k, 4\sigma) \cap B(x_j, 4\sigma)=\emptyset$ if $d(x_k, x_j) \ge 9 \sigma$ and
$B(x_k, \sigma/2) \cap B(x_j, \sigma/2) =\emptyset$.
\hfill\Box\medskip

\section{Wave equation: stability for the unique continuation}

Consider the initial-value problem
for the  wave
equation
\begin{eqnarray}
& &\p_t^2 w -\Delta_g w=0\hbox{ in } M\times \R,  \label{12}\\
& &w|_{t=0}=v, \nonumber\quad
w_t|_{t=0}=0,\nonumber
\end{eqnarray}
on $(M,g,p)  \in \overline{\bM_{\dtext n}(R, D, i_0)}$ and denote its solution by $w=W(v)$.
Our main interest lies in the case when $v \in  \H^s_{\Lambda}(M)$, $\Lambda >0$,
\beq \label{conditions for u}
\H^s_{\Lambda}(M)=\{v\in H^{s}(M):\,
\|v\|_{H^{s}(M)}\leq \Lambda\}
\eeq
and we assume in the following that 
\beq 
\label{eq s} 
\frac 3 2 < s <2
\eeq
 {\Mtext and denote $\H^0_{\Lambda}(M):=\{v\in L^2(M):\ \|v\|_{L^2(M)}\leq {\Lambda}\}$.}
Using the Fourier decomposition we show that, if $v \in H^{s}(M)$, then
\beq\label{problem}
\|w\|_{H^{s}(M\times [-T,T])}\leq 6 \sqrt T \|v\|_{H^{s}(M)} \leq {\novtekst \SlavaCsixagain} \|v\|_{H^{s}(M)}, \quad T <2 D,
\eeq
where
${\novtekst \SlavaCsixagain}= 6 \sqrt{1+D^2}$.

\begin{figure}
\begin {picture}(300,100)(15,0)
 \put(00,0) {\includegraphics[height=5.0cm]{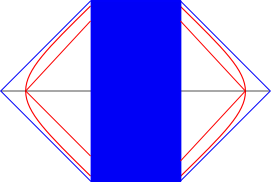} }
\put(38,53){$\Sigma_\gamma$}
\put(-5,53){$\Sigma_{0}$}
\put(83,150){$\Gamma(z,T)$}
\put(53,115){$\D$}

 \put(220,0){\includegraphics[height=5.0cm]{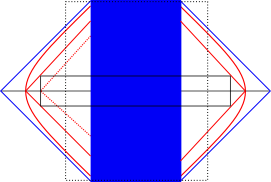}}
\put(227,43){$\Sigma_0$}
\put(273,87){\small $ \Sigma_{2\gamma}$}
\put(267,43){\small $\Sigma_{\gamma}$}
\put(293,150){\quad$\Gamma(z,T)$}
\put(253,150){$\tilde  \Gamma(z,T)$}

\end {picture}
\caption{
\label{fig_3C}
{\it Left: Assume that  the function $u({x,t})$ vanishes in $\Gamma(z,T)$. If  $u({x,t})$ satisfies the wave equation $Pu=0$, then $u({x,t})$  vanishes in the double cone $\Sigma_{0}=\D(z,0, T)$. Theorem \ref{loc_stability} states that if $Pu=f$ is  small, then $u({x,t})$ is small in the  domain $\D=\D(z, \ga, T)$ (that has a red, curved boundary). In the figure, we consider also the double cone $\Sigma_\gamma=\Sigma(z,\gamma,T)\subset \D$.
Right:
 In Corollary \ref{main} 
we assume that  $u({x,t})$ is a solution to  $Pu=0$ and that $u$  is small in the set $\tilde \Gamma(z,T)=B(z, r_0/16+\ga) \times (-T +r_0/16,\, T-r_0/16),$ 
marked with a dotted black boundary.
Note that
$\Gamma(z,T)\subset \tilde \Gamma(z,T)$. 
 Then we apply  
Theorem \ref{loc_stability}, to see that $u$ is small in $ \D=\D(z, \ga, T)$, interpolation, and trace theorem in the black cylindrical set  $K= B(z,T-2\ga)\times (-\gamma,\gamma)$
to see that 
 $u|_{t=0}$ is small in $B(z,T-2\ga)$ that is the intersection of $\{0\}\times M$ and the domain $\Sigma_{2\gamma}= \Sigma(z, 2\ga, T)$, that is shown in the figure with dotted red lines.}
\vspace{-5mm}
}
\end{figure}
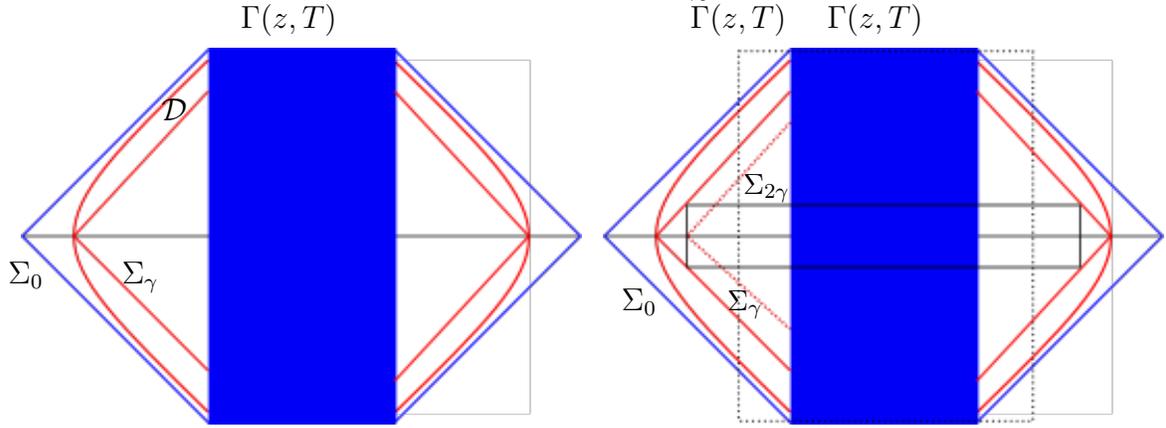

Associated to the wave operator are the double cones of influence. To define these, let $V\subset M$ be open, $T \in \R_+$.
Denote by
$$
\Gamma(V, T):=V\times (-T,T).
$$
Then the double cone of influence is given by
 $$\Sigma(V, T):=\{({x,t}); d(x,V)+|t| < T\}.$$
By Tataru's uniqueness theorem \cite{Ta}, \cite{Ta1}, if $u$ is a solution to  \eqref{12} 
in $M \times (-T, T)$,
which satisfies $u=0$ in $\Gamma(V, T)$, then $u=0$ in $\Sigma(V,  T)$.
However, for our purposes we need an explicit estimate which follows from Theorem 3.3 in \cite{BKL}.
 To formulate the results we introduce, for
\beq \label{22.06.1}
{\ftext \quad 0< \ga \le \frac{r_0}{32}, \quad \frac{r_0}{8} \le T <2D,} 
\eeq
with $r_0$ fulfilling \eqref{4.09.2017} and
$z\in M$, the  domains 
\beq \label{2.3}
& &\Gamma=\Gamma(z, T)=B(z,r_0/16) \times (-T+r_0/16, T-r_0/16),
\\ \nonumber
\D&=&\D(z,\gamma,T)=\{({x,t}):\, (T-d(x,z))^2-t^2 \geq \gamma^2, \,|t| < T-r_0/16\},
\\ \nonumber
& &
\Omega(T)= M \times (-T+r_0/16, \, T-r_0/16).
\eeq
{\ftext Also, let for $b\in \R$,
 \beq \label{2.3b}\\
 \nonumber 
 \Sigma(z,b\gamma,T) = \{(x,t)\in M\times   \R:\ |t|\le T-r_0/16 ,\ |t| \leq  T-b\gamma - d_g(x,z)\}
 \eeq
 be the ``domain of influence'' corresponding to the cylinder $\Gamma(z, T)$.
 Observe that 
 $
\Sigma(z,\gamma,T) \subset \D(z,\gamma,T)  \subset \Sigma(z,0,T).
$
}

{\dtextRB In the following we formulate the stability results for the unique continuation in \cite{BKL}. We note 
that similar results have been obtained by Luc Robbiano in \cite{Robbiano}   with $\theta=1$, but with a  loss
 in the domain of dependence and later by  C. Laurent and M. Leautard
 in \cite{LaLe} with $\theta=1$, but without an explicit calculation
 of the constants in the  domain of dependence.}
 
\begin{theorem}\label{loc_stability}
Let $(M,g,p) \in \overline{\bM_{\dtext n}(R, D, i_0)}$.
Let $P=P(x,D)=\p^2_t-\Delta_g$  be the wave operator associated with  $M$.
Assume that $w({x,t})=0$ for all $({x,t})\in \Gamma$.
Then, for any
{\dtextRB $0 <\theta<1$, there is }
$c_{206}{\ktext(\ga,\theta)} \geq 1$,
{\ptext depending only on $n$, $R, D, i_0,r_0,\theta$ and $\gamma$}
 such that the following stability estimate holds true:
\begin{eqnarray*} \label{2.1}
\|w\|_{L^{2}(\D(z, \ga, T))} \le 
c_{206}{\ktext(\ga,\theta)} \frac{\|w\|_{H^1(\Omega(T))}}
{\Big(\ln \Big(1 + \frac{\|w\|_{H^1(\Omega(T))}}{\|Pw\|_{L^2(\Omega(T))}}\Big)\Big)^\theta}\,,
\end{eqnarray*}
where $c_{206}(\ga,\theta)$ is such that
\beq \label{19.06.1}
c_{206}{\ktext(\ga,\theta)} = c_{205}{\ktext ( \theta)} \exp(\ga ^{-c_{200}}),\quad c_{200}=58(n+1)+1,
\eeq
and  {\ktext $ c_{205}( \theta)\ge 1$ depends on 
$\theta, n, R, D, i_0, r_0$.}
Moreover, for any $0 \leq m \leq 1$,
\begin{eqnarray}\label{m estimate}
\|w\|_{H^{1-m}(\D(z, \ga,T))} \le 
c_{206}{\ktext(\ga,\theta)}^{m} \frac{\|w\|_{H^1((\Omega(T))}}
{\Big(\ln \Big(1 + \frac{\|w\|_{H^1((\Omega(T))}}{\|Pw\|_{L^2((\Omega(T))}}\Big)\Big)^{\theta m}}\,.
\end{eqnarray}
\end{theorem}

\noindent
{\bf Proof.}
Theorem \ref{loc_stability} follows from Theorem 3.3  in \cite{BKL}
with $\ell=r_0/16$ and
$\D(z,\gamma,T) = S(z, r_0/16, T, \ga)$.
Using that $w=0$ in $\Gamma$, the domain $\Lambda$ in the final equation of Theorem 3.3
can be changed into $\D(z, \ga,T)$.
Moreover, for $\theta <1$, the function $f_\theta(a,   b),\,a, b>0,$
\beq\label{f function grows}
f_\theta(a,b)=
a \left(\hbox{ln}(1+\frac{a}{b})\right)^{-\theta},
\eeq
increases when either $a$ or $b$ increases. Thus,
we can change $\|w\|_{H^1(\Omega_1)}$ and $\| P w\|_{L^2(\Omega_1)}$ in Theorem 3.3
to $\|w\|_{H^1(\Omega(T))}$ and $\| P w\|_{L^2(\Omega(T))}$. 
{\Mtext Note that, although the results in \cite{BKL} are formulated for $M \subset \R^n$, 
they can be easily reformulated for an arbitrary compact Riemannian manifold which
possess $C^5$-smooth covering by coordinate systems with $C^4-$smooth metric tensors.
{\ktext To consider parameters  (\ref{19.06.1})  (see the Appendix for details),
we will fix the value of $\theta$  to be} 
\beq\label{eq theta}
\theta = 1/2,
\eeq
 for simplicity. In the general case, we write $c_{205}{\ktext ( \theta)} $ as $\theta-$dependent.
We recall that the constants
in \cite{BKL} (see (3.1)) explicitly depend on parameter {\ktext $ \slavaConetheta >1$} such that
\beq\label{SlavacnonumberA}
{\ktext  \slavaConetheta^{-1}|\xi|^2 \leq g^{jk}(x)\xi_j\xi_k\leq  \slavaConetheta|\xi|^2,
\quad \|g^{jk}(x)\|_{C^4(M)}\leq  \slavaConetheta.}
\eeq
{\dtextRBX
Using harmonic coordinates in balls of radius $r^{(har)}$, this condition is fulfilled  due to \eqref{26.2A}, which also implies $d_g(x,z) \in C^3$. 
}
\hfill\Box\medskip

\noindent
Our main interest will be an estimate for $v(\cdot)=w(0, \cdot)$ in (\ref{12}) in the domain
$B(z, T-2\ga)$. 

\begin{corollary}\label{main} 
Assume \eqref{22.06.1} and let 
{\dtextRB $\theta\in [1/2,1).$} {\Mtext Also, let 
$ {{\Lambda_1}}>0$ and}
$\e_2 \in (0,{\Mtext {{\Lambda_1}}}]$
and
$v \in  \H^s_{{\Mtext {{\Lambda_1}}}}(M)$.
Denote by
$w=W(v)$  the solution to initial-value problem (\ref{12}) and assume that,
\beq\label{small observations}
\|w\|_{L^2(B(z, r_0/16+\ga) \times (-T +r_0/16,\, T-r_0/16))}\leq  \e_2.
\eeq
Then, calling $ \beta=\theta^2/2$ and defining $\e_1
:=\mathcal E_1(\e_2; \theta, \gamma,  {\Mtext {{\Lambda_1}}})$, we get
\beq
\label{19.06.3}
\hspace{-5mm}\|v\|_{L^2(B(z,T-2\ga))} &\leq& \e_1
\eeq
where, for $c_{202}=c_{202}(\theta, \ga)$, {\ktext and $\RobCthirty(\theta)$
depending only on $\theta$  and $n,R,D,i_0,r_0$}
\beq
\label{explicit 1}
 \mathcal E_1(\e_2; \theta, \gamma,  {\Mtext {{\Lambda_1}}}) &=&
 c_{202}\, \frac{{\Mtext {{\Lambda_1}}}}{ \ga^{(2-\theta/2)}
\left(\ln\left[1+ \ga  {\Mtext {{\Lambda_1}}}^{(s-1)/s}\e_2^{-(s-1)/s}  \right]  \right)^{\beta}},\quad
\label{explicit 1B}
\\ \label{30.06.5}
c_{202}&=& \RobCthirty \exp{\left(\ga^{-(c_{200}\, \theta/2)}\right)},
\\
\nonumber
\RobCthirty&=& \RobCthirty(\theta)\, {\Mtext \ge 1}.
\eeq
\end{corollary}

\noindent 
{\bf Proof.} 
Let the cut-off function $\eta(x)\in C^2_0(B(z, r_0/16+\ga/2))$ be equal to
one in $B(z,r_0/16)$ and $\|\eta\|_{C^i(M)}\leq C\gamma^{-i},\, i=0, 1,2$. Then
$w_\eta(x,t)=(1-\eta(x)) w(x,t)$ vanishes in $\Gamma$ and we have
$(\p_t^2- \Delta) w_\eta(x,t)=F,$ where
\beq \label{18.1a}
F(x,t)&=&\big(\Delta_g \eta(x)\big )\,w(x,t)+ 2g(\nabla \eta(x),\nabla_x w(x,t)) \\ \nonumber
&=& \big(\Delta_g \eta(x)\big )\, \left(\tilde \eta (x) w(x,t)\right)+ 2g(\nabla \eta(x),\nabla_x \left( \tilde \eta (x) w(x,t))\right) := F_1+F_2
\eeq
Here $\tilde \eta(x)\in C^2_0(B(z, r_0/16+\ga))$ is equal to
one in $B(z, r_0/16+ \gamma/2)$ and $\|\tilde \eta\|_{C^i(M)}\leq C\gamma^{-i},\, i=0, 1,2$.
Clearly, by hypothesis
$$\|F_1\|_{L^2(M \times (-T +r_0/16,\, T-r_0/16))}\leq  C \ga^{-2} \e_2.$$ 
To estimate
$F_2$, observe that,$
\| \tilde \eta w\|_{H^s(M \times (-T+r_0/16, T-r_0/16))} \le C \ga^{-s}  {\Mtext {{\Lambda_1}}},
$
where we have also used (\ref{problem}). Since
$\|\tilde \eta w\|_{L^2(M \times (-T+r_0/16, T-r_0/16))}  \le  \e_2$,
by interpolation arguments,
we get
\beq 
\| \tilde \eta w\|_{H^1(M \times (-T+r_0/16, T-r_0/16))} \le C \ga^{-1}  {\Mtext {{\Lambda_1}}}^{1/s} \e_2^{1-1/s}
\eeq
Since $\hbox{supp}(\nabla \eta) \cap \hbox{supp}(\nabla \tilde \eta) =\emptyset$, this implies
\ba
\|F_2\|_{L^2(M \times (-T +r_0/16,\, T-r_0/16))}\leq  C \ga^{-1}  {\Mtext {{\Lambda_1}}}^{1/s} \e_2^{1-1/s},\\
\|F\|_{L^2(M \times (-T+r_0/16, T-r_0/16))}\leq C \ga^{-2} {\Mtext {{\Lambda_1}}}^{1/s} \e_2^{1-1/s},
\ea
where we used $\e_2 \le    {\Mtext {{\Lambda_1}}}$. As $s>1$, we have
\ba
\|  w_\eta\|_{H^1(M \times (-T+r_0/16, T-r_0/16))} \le C \ga^{-1} {\Mtext {{\Lambda_1}}}.
\ea
Using growth properties of the function $f_\theta$
of form (\ref{f function grows}), it follows from Theorem \ref{loc_stability} that
\beq \label{18.4a}
\| w_\eta\|_{H^{1-\theta/2}(\D)} \leq
C  c_{206}{\ktext(\ga,\theta)}^{\theta/2}
\frac{\ga^{-1}  {\Mtext {{\Lambda_1}}}}{\left(\ln\left[1+ \ga   {\Mtext {{\Lambda_1}}}^{(s-1)/s} \e_2^{-(s-1)/s}\right]\right)^{\beta}}.
\eeq
Now observe that by the trace-theorem, for any $ \a >1/2\, $ there exists $\RobCeleven=\RobCeleven(\a)$
such that, for $r \geq r_0/16, \,z \in M$:
\beq\label{Sobo1}
 \|w(\,\cdotp,0)\|_{L^2(B(z, r))} \leq
\RobCeleven \,\gamma^{-\a} \|w\|_{H^\a( B(z, r)\times (-\ga,\ga))},\\
\label{Sobo2}
\|w(\,\cdotp,0)\|_{L^2(B(z, T-2 \ga))}  \leq  \RobCeleven\, \gamma^{-\a} \|w\|_{H^\a (\D(z,\gamma,T))}.
\eeq
It follows from  (\ref{Sobo2}) with $\a=1-\theta/2$ and (\ref{18.4a})
that,
\begin{eqnarray}\label{m estimate2}
\|w_\eta(\cdotp,0)\|_{L^2(B(z,T-2\gamma))}
\le C \RobCeleven  \,
\frac{c_{206}{\ktext(\ga,\theta)}^{\theta/2}  {\Mtext {{\Lambda_1}}} }{ \ga^{2-\theta/2}
\left(\ln\left[1+ \ga  {\Mtext {{\Lambda_1}}}^{(s-1)/s}  \e_2^{-(s-1)/s} \right]\right)^{\beta}}.
\end{eqnarray}
Next define $\a=(1-\beta) s +\beta >1/2$. Then by  interpolation,
\ba
\|\eta w  \|_{H^\a( B(z, r)\times (-\ga,\ga))} \le
c_{201}\,\| \eta w\|^{\a/s}_{H^s( B(z, r)\times (-\ga,\ga))} \,\| \eta w\|^{(s-\a)/s}_{L^2(B(z, r)\times (-\ga,\ga))}.
\ea
Using the fact that $\hbox{supp}(\eta) \subset B(z, r_0/16 +\ga)$,
we can apply \eqref{Sobo1} with $r=r_0/16 +\ga$, the previous inequality and (\ref{small observations}), to obtain
\beq \label{19.2c}
\|\eta(\cdot) w(\cdot, 0)  \|_{L^2(B(z,T-2\gamma))}
\le \RobCeleven\ga^{-\a} c_{201}({\novtekst \SlavaCsixagain}  {\Mtext {{\Lambda_1}}})^{\a/s}\epsilon_2^{\beta(s-1)/s}
\nonumber \\
\le 
\RobCeleven \ga^{\beta-\a} c_{201}{\novtekst \SlavaCsixagain}^{\a/s}  {\Mtext {{\Lambda_1}}} 
\left(\ln\left[1+ \ga  {\Mtext {{\Lambda_1}}}^{(s-1)/s}  \e_2^{-(s-1)/s} \right]\right)^{-\beta}
.
\eeq
Here at the last step we
use the fact that
$X \ge \ln(1+X)$ for $X >0$, with $X = \gamma  {\Mtext {{\Lambda_1}}}^{(s-1)/s}  \e_2^{-(s-1)/s}$.
Recall  that $v(x)=w_\eta(x,0)+\eta(x) w(x,0)$.
Comparing (\ref{m estimate2}) and (\ref{19.2c}), we obtain
equation (\ref{explicit 1}).
The  coefficient $c_{202}$ defined in (\ref{30.06.5})  fulfills the inequality
$$c_{202} \ge C \RobCeleven c_{206}{\ktext(\ga,\theta)}^{\theta/2} \ga^{\theta/2-2}
+ \RobCeleven c_{201} {\novtekst \SlavaCsixagain} ^{(1-\beta)+\beta/s} \ga^{(\beta-1)s},
$$
by using (\ref{19.06.1}) and a proper multiplicative coefficient $\RobCthirty$ independent on $\gamma$.
\hfill\Box\medskip

\section{Computation of the projection}

\subsection{Domains of influence} \label{Approximate projections}
{\dtext Let $(M,g,p)\in \overline{\bM_{\dtext n}(R, D, i_0)}$. By Proposition \ref{uniform_covering}, we can choose $L-1$ points $z_j$, $j=1,2,\dots,L-1$ that
form a $\tau_0$-net in $B_M(p,r_0/4)$. Here, $L$ is bounded by a uniform constant.} 

{\dtextRB In  Lemma \ref{ga-separation} we showed that for any $\sigma$ there are 
$N_1(\sigma)$ points, that we enumerate as  $z_{L}, \dots, z_{{\dtext L-1+N_1(\sigma)}}$,
which form a maximal $\sigma$-separated net in $B(p, r_0/4)$ 
and the balls  $B(z_k, 4 \sigma)$, {\dtext $k=L,\dots,{\dtext L-1+N_1(\sigma)}$, satisfy} the  finite intersection property with at most  $N_F$ intersections.
In this section we consider arbitrary $\sigma$, which value will be specified later, and points  $z_\ell, \,\ell=L, \dots, L-1+ N_1(\sigma)$ that satisfy the conditions
of Lemma \ref{ga-separation}. Also, below is $3/2<s<2$.}}

Our {\dtext next goal} is to approximately construct
the values of the distance functions from a variable point $x \in M$ to all points
$z_\ell\in B(p,r_0/4)$, $\ell=1,2, \dots,  {\dtext L-1+N_1}(\sigma)$, defined in Lemma \ref{ga-separation}.
The main step is to
approximately compute the Fourier coefficients of the  functions of form
$\chi_\Omega (x)v(x)$, where $\chi_\Omega(x)$ are the characteristic functions of some special subdomains
$\Omega \subset M$
and $v(x)$ has  a finite Fourier expansion. These subdomains $\Omega$ are defined using distances to
$L$ points
$\{z_1, \dots, z_{L-1}, z_i \}, $ where $i \in \{ L, \dots,  {\dtext L-1+N_1}(\sigma)\}$ is arbitrary.
For  $i \in \{ L,L+1, \dots,  {\dtext L-1+N_1}(\sigma)\}$,
let
$
K_i=\{1,2,\dots,L-1\}\cup \{i\}
$
and define $\mathcal A^{(i)}$ to be the set of those $\a=(\a_\ell)_{\ell=1}^{ {\dtext L-1+N_1}(\sigma)}\in  
\R^{ {\dtext L-1+N_1}(\sigma)}$, such that 
\beq
\nonumber
&&\hbox{$r_0/8\le \a_\ell \le 2D,$ if $\ell  \in K_i,$ 
}\\
\label{06.07.1}
& &\hbox{$\a_\ell=0$, if $\ell \not \in K_i$.}
\eeq
Below, we will assume that
\beq\label{gamma and sigma}
\gamma\leq \sigma.
\eeq
We
}

 {\qtext denote 
\beq \label{2.3 tilde}
& &\hspace{-2cm}\tilde \Gamma(z, T)=B(z,r_0/16+\gamma) \times (-T+r_0/16, T-r_0/16).
\eeq
}
{\novtekst {Next we fix for a while the index $i \in \{ L, \dots,  {\dtext L-1+N_1}(\sigma)\}$.}}
To construct subdomains $\Omega$, we start with
 observation sets 
${\qtext \tilde \Gamma(\a)}$, $\a \in \mathcal A^{(i)}$,
\beq \label{16.4}
{\qtext \tilde \Gamma(\a)=
\bigcup_{\ell \in K_i} \tilde \Gamma(z_\ell,  \a_\ell)}.
\eeq
At last, for $b \in \R$, we define
\beq\label{Ma set}
M(\a,\,b \derho)=\bigcup_{ \ell \in K_i}
 B_M(z_\ell, \a_\ell+b\derho).
\eeq
{\ftext Then the corresponding domains of stable unique continuation are 
\beq \label{15.1}
& &\D(\a)= \bigcup_{ \ell \in K_i}
 \D(z_\ell,  \derho, \a_\ell),\quad \D(\a,b\gamma)= \bigcup_{ \ell \in K_i}
 \D(z_\ell,  \derho, \a_\ell+b\gamma),
 \eeq
 and the corresponding double cones of influences are given by
   \beq  \nonumber
& &\Sigma(\a)= \bigcup_{ \ell \in K_i}
 \Sigma(z_\ell,  \derho, \a_\ell),\quad \Sigma(\a,b\gamma)= \bigcup_{ \ell \in K_i}
 \Sigma(z_\ell,  \derho, \a_\ell+b\gamma),
 \eeq
 }

We have the following volume estimate.

\begin{lemma}\label{y-lemma 1} $a)$ Let $\a \in \mathcal A^{(i)},\, i= L,\dots, {\dtext L-1+N_1}(\sigma)$ and
\beq \label{15.2}
A=A(\a, \ga) = \{x\in M:\ d(x,\p M(\a, \,{\qtext 3\derho}))\leq {\qtext 5\derho}\}.
\eeq
Then, there is a uniform  constant $\SlavaCnonumber>0$,
{\ztext depending only on $n$, $R, D, i_0$ and $r_0$,}
 such
that
$$
\vol(A)\leq \SlavaCnonumber L \derho.
$$
$b)$ Consequently, by defining $b(s)$, for $\frac32 <s <2,$  as 
$b(s) = 1/2, \,n=2, 3$ and $b(s) = s/n, \, n \ge 4$,
we see that there is a uniform constant $\Slavacnonumber{\qtext(s)}$,
{\ztext depending only on  $s$, $n$, $R, D, i_0$ and $r_0$,}
 such that  
\beq \label{Sobolev}
\|\chi_{ B(z_\ell, \a_\ell+{\qtext 8\derho}) \setminus B(z_\ell, \a_\ell-{\qtext 2\derho}} )\,v\|_{L^2(M)}
\le \Slavacnonumber(s)  \ga^{b(s)} \|v\|_{H^s(M)}.
\eeq

\end{lemma}
\medskip
{\bf Proof.} $a)$ Let $ x \in A.$ Then, for some 
$\ell\in K_i$,
\beq \label{61}
x \in B(z_\ell, \a_\ell+{\qtext 8\derho}) \setminus B(z_\ell, \a_\ell-{\qtext  2\derho} ).
\eeq
Since $\|d\,\exp_{z_\ell}|_{\bf v}\|$ is uniformly bounded on $\overline{\bM_{\dtext n}(R, D, i_0)}$ for
${\bf v} \in T_{z_\ell}M,\,|{\bf v}| \le 2D$,
$\vol\left( B(z_\ell, \a_\ell+{\qtext 8\derho}) \setminus B(z_\ell, \a_\ell-{\qtext  2\derho} ) \right)
\le C \derho$, {\novtekst{for all $\ell \in K_i$}}. 

\noindent
$b)$  Similar to part $a)$, we have
$
\vol\left(B(z_\ell, \a_\ell+{\qtext 8\derho}) \setminus B(z_\ell, \a_\ell-{\qtext 2\derho} \right) )\le c \ga.
$
Together with the H\"older inequality and  the Sobolev embedding $H^{s}(M) \to L^q(M),\,
\frac{1}{q} =\frac12 -\frac{s}{n},$ 
(or $ C^0(M)$ for $n= 2, 3$), this implies (\ref{Sobolev}).
Note that $\Slavacnonumber(s)$ is a uniform constant as the embedding can be done in
harmonic coordinates defined in balls with uniform radius.
\hfill\Box\medskip

\subsubsection{Cut-off estimates and finite dimensional projections}\label{subsec 4.2 C}

{\dtextRB Let us apply  Lemma \ref{ga-separation}, with $\ga$ instead of $\sigma$,
to obtain points $x_\ell\in M$, $\ell=1,2, \dots, N(\ga)$
such that the balls $B(x_\ell, 2 \derho),\, \ell=1,2, \dots, N(\ga)$ are a covering of $M$.}

Let $\psi_\ell:M\to \overline \R_+, \,\psi_\ell \in C^6_*(M)$
be {\novtekst in harmonic coordinates} a partition of unity for the covering $B(x_\ell, 2 \derho)$ that satisfy
\beq\label{y-2} \nonumber
\|\psi_\ell\|_{C^{k,\beta}(M)}&\leq& c_{k,\beta} \derho^{-(k+\beta)},
\quad k=0, 1, 2,\quad 0\le \beta  {\dtext <1};
 \\
& &\supp(\psi_\ell) \subset B(x_\ell, 2\ga),\quad
\sum_{\ell=1}^{N(\ga)} \psi_\ell(x)=1.
\eeq
{Below, we use $\Lambda_s\geq 1.$}
\begin{lemma}\label{lemma: tilde u exists} 
{\dtextRB For $\frac32<s<2$ there is}
$\Slavacseven(s)\ge\,{\Mtext 1}$, in (\ref{20.4}) such that, for any
$u\in  \H^s_{ \Lambda_s}(M)$,  $\,i \in \{L, \dots, {\Mtext {\dtext L-1+N_1}}(\sigma)\}$
and $\a \in \mathcal A^{(i)}$, the following holds true:
There 
\beq\label{y-4b}
& &u_\a \in  \H^s_{{\Rtext \frac14}  \SlavaCone(s;\derho) \Lambda_s}(M) \cap {\Mtext \H^0_{\Lambda_s}(M)},\quad
u_\a(x)=0,\,\, \hbox{if }x\in M(\a, {\qtext \derho}), \\ 
& & u_\a(x)=u(x),\,\, \hbox{if }x\in M \setminus M(\a, {\qtext  7\derho}),
\nonumber
\eeq
where
\beq \label{20.4}
  \SlavaCone(s, \derho)=  \Slavacseven(s) \derho^{-s}.
\eeq 
\end{lemma}

\noindent {\bf Proof.}
Define
\beq\label{def of tilde u}
u_\a(x)=
\Psi(x)u(x),\quad \Psi(x)=
\sum_{\hbox{supp}(\psi_\ell)\cap M(\a,{\qtext  3\derho})=\emptyset} \psi_\ell(x).
\eeq
{\dtextRB For a general $w\in H^2(M)$ we have 
the following estimate in Sobolev spaces with $\frac32<s<2$  
\beq\label{eq: interpolation}
\|\Psi w\|_{H^s(M)}\leq C\|\Psi\|_{C^s(M)}\|w\|_{H^s(M)}\leq C\,m\, c_{2,0}\gamma^s\|w\|_{H^s(M)},
\eeq
where $m$ is the number of elements in the set $\{\ell:\ {\hbox{supp}(\psi_\ell)\cap M(\a,{\qtext  3\derho})=\emptyset} \}$ satisfying  $m\leq N(\ga)$. 
}
Thus the existence of  $\Slavacseven(s)$ such that the claim holds follows then
from the finite intersection property
of $B(x_\ell, 2 \derho)$, see Lemma \ref{ga-separation}, and estimates 
(\ref{y-2}).
\hfill\Box\medskip


\subsection{Unique continuation for approximate projections}

Corollary \ref{main}  implies the following result.
{\dtextRB
Note that the notations $\e_2$ and $\mathcal E_2$ are introduced in order to distinguish $\e_2$ from its upper bound $\mathcal E_2$, written as an expression dependent on $\e_1$.
}
{\ptext 
Later, in formula (\ref{new-e0e1}) we set  $\e_1$  to have
a specific value and substitute it in the expression $\mathcal E_2( \frac{\e_1}
{4L};\,\theta, \ga,  \Lambda_s)$  of formula (\ref{19.1 B})
to obtain a specific value for $\e_2$.}

\begin{corollary} \label{E_2} 
Assume that $v$ satisfies
\beq\label{newIC}
\|v\|_{H^s(M)}\leq 
\SlavaCone(s, \derho) \Lambda_s\quad\hbox{and}\quad \|v\|_{L^2(M)}\leq  \Lambda_s, 
\eeq
 with $\SlavaCone(s, \derho)$ defined in \eqref{20.4}, and assume \eqref{22.06.1}.
Let $\e_1 <\Lambda_s$ and $\e_2 \le  \E_2\left( \,  \frac{\e_1}{4L};\theta, \derho, \Lambda_s\right)$    where
\beq \label{19.1}
\mathcal E_2( \frac{\e_1}
{4L};\,\theta, \ga,  \Lambda_s)
\;=\; \frac{\Lambda_s   \ga^{s/(s-1)}}{\left(\exp\left[\Big( 
 \Lambda_s 4L \e_1^{-1}  \gamma^{-(2-\theta/2)}\RobCthirty  
\,\exp( \gamma^{-c_{200}}) \Big)^{1/\beta}\right]\right)^{s/(s-1)}} \quad
\eeq
Let $w= W(v)$ satisfy 
\beq \label{19.2}
\|w\|_{L^2(\tilde \Gamma(z_\ell,  \a_\ell))} \le \e_2
\eeq
on the domain \eqref{2.3 tilde}.
Then, for $ \ell \in K_i$,
\beq \label{19.3}
\|w(0, \cdot)\|_{L^2(B(z_\ell, \a_\ell-{\novtekst 2}\derho))} \le  \frac{\e_1}{4L}, 
\quad \|w(0, \cdot)\|_{L^2(M(\a,-{\novtekst 2}\derho))} \le   {\Mtext \frac 14} \e_1.
\eeq
\end{corollary}


\noindent {\bf Proof.}  
From a small modification of the proof of Corollary  \ref{main} we still can obtain the estimate \eqref{19.06.3} in the following way.
{\dtextRB
The main point is to replace the initial condition $\|v\|_{H^s(M)} \leq \Lambda_1$ with $\eqref{newIC}$.
We then deduce the corresponding estimate for the solution $w=W(v)$ of the wave equation, with $T=\a_\ell$ and $z=z_\ell$,
}
\ba
\|w\|_{H^s(M\times[-T,T])}\leq C\SlavaCone(s, \derho) \Lambda_s\quad\hbox{and}\quad \|w\|_{L^2(M\times[-T,T])}\leq  C\Lambda_s, 
\ea
{\dtextRB Let $\eta$ and  $\widetilde \eta$ be the smooth localizers defined in the proof of Corollary  \ref{main}. 
Calling again $w_{\eta} = (1-\eta(x))w$ and using the definition of $\SlavaCone$ in \eqref{20.4} we get, 
}
\ba
\|\eta w\|_{H^s(M\times (-T+r_0/16,T-r_0/16))} \le C \gamma^{-s} \Lambda_s, \\
\|\widetilde  \eta w\|_{H^s(M\times (-T+r_0/16,T-r_0/16))} \le C \gamma^{-s} \Lambda_s, \\
\| w_{\eta}\|_{H^s(M\times (-T+r_0/16,T-r_0/16))} \le C \gamma^{-s} \Lambda_s, 
\ea
and the intermediate $H^m$ norms follow  by interpolation.  
Here the constant C is
{\dtextRB
dependent of
$c_3(s)$ and independent of $\ga$. Consequently,
\ba
\| \tilde \eta w\|_{H^1(M \times (-T+r_0/16, T-r_0/16))} \le C \ga^{-1}  {\Mtext {{\Lambda_s}}}^{1/s} \e_2^{1-1/s},
\\
\|F\|_{L^2(M \times (-T+r_0/16, T-r_0/16))}\leq C \ga^{-2} {\Mtext {{\Lambda_s}}}^{1/s} \e_2^{1-1/s}.
\ea
Using growth properties of the function $f_\theta$
we get \eqref{18.4a}.
Also \eqref{m estimate2} still holds.
Therefore we obtain \eqref{explicit 1B}, where the new constant $\RobCthirty$ in $\eqref{30.06.5}$ now depends on $c_3(s)$.
}
Next we observe that formula \eqref{explicit 1B} implies that when 
$\e_1={4L}\mathcal E_1(\e_2; \theta, \gamma, \Lambda_s)$, we have
\ba
\e_2 = \frac{\Lambda_s \gamma^{s/(s-1)}}{\left( \exp\left[ \left(  \Lambda_s 4L 
\RobCthirty \e_1^{-1}  \gamma^{-(2-\theta/2)} \exp(\ga^{-(c_{200}\, \theta/2)})
  \right)^{1/\beta} \right] -1 \right)^{s/(s-1)}}  ,
\ea
and $\mathcal E_2$ is defined  by removing $-1$ from the denominator of the expression above, and by replacing $\exp(\ga^{-(c_{200}\, \theta/2)})$ with $\exp(\ga^{-c_{200}})$.
This is done to simplify the calculations of the paper.
The relation \eqref{19.3} follows by imposing on $\e_2$ the $\mathcal E_2$-bound.
\hfill\Box\medskip

Under the conditions of the Corollary and from the growth properties of ${\mathcal E_2}(\e_1)$ it follows that
\beq \label{30.06.7}
\e_2 \le \mathcal E_2( \frac{\e_1}{4L};\,\theta, \ga,  \Lambda_s) \le  \frac{\e_1}{4L}, \quad \e_1 \in (0, \Lambda_s].
\eeq

\subsection{Approximate projections}
{\Mtext
Let $\e_0,\e_1,\e_2$
satisfy  
\beq\label{epsone defined}
\e_0 \le\frac{\Lambda_s}{10},
 \quad
 \e_1=\frac{\e_0^2}{10 \Lambda_s},  \quad
\e_2 =\E_2\left( \,  \frac{\e_1}{4L}, \theta, \derho, \Lambda_s \right).
\eeq

\subsubsection{Finite data with and without errors}
Below we will use several parameters, and for the sake of clarity of presentation,
we have gathered these parameters in this subsection and tell how those will be used.

{\mtext Below, 
we will use $\jnull\in \Z_+$ satisfying
\beq  
\label{24Bnew} \origjnull&\geq & 
\widehat \jbasis({\Mtext \frac{\e_2}{8}};  \derho, \Lambda_s),
\eeq
where
\beq\label{def SlavaCthree(s)}
\widehat \jbasis( \e_*;  \derho, \Lambda_s)=
\SlavaCthree \ga^{-n} \left({\Rtext \frac{\Lambda_s}{\e_*} }\right)^{\frac ns}\ \ \hbox{and}\quad
\SlavaCthree(s)= \Slavacseven(s)^{\frac ns} {\RobCthree^{\frac n2}}
 {(\SlavaCsixagain+1)^{\frac ns}}.
\eeq
We also use 
$\janul\in \Z_+$ satisfying %
\beq \label{janolla}
\origjnull \leq \janul \leq  2^{n/2}(\RobCthree)^{n}  \origjnull
\eeq
Moreover,
we use 
\beq \label{29.5}
\delta &\leq & \widehat \delta_0(\e_2, \ga, \janul, \Lambda_s)=
\Slavacsevenagain\,  \frac 1{\janul}\,  
\frac{\e_2}{\Lambda_s},
\eeq
where 
 $
\Slavacsevenagain=\min(\RobCthree^{-1} , \frac {(1+{\ktext 2\RobCthree}
)^{-1/2}}{100(1+D)^{3/2}L}),$
and $J$ satisfying
\beq\label{bound for J B}
( \RobCthree^{-1}  \delta^{-1})^{n/2}\leq J\leq   (2\RobCthree\delta^{-1})^{n/2},
\eeq
cf. Remark \ref{Rem: size of J}.
Note that (\ref{bound for J B}) implies that $\lambda_J\ge \delta^{-1}$, see 
 Def.\ \ref{def:3} (ii) and (\ref{26.5}).

The use of the above parameters are the following.
We will assume that we are given  the ball $(B_e(r_0), g^a)$ and the pairs
$\{(\lambda_j^a,\, \varphi_j^a|_{B_e(r_0)}) ;\, j=0,1,2,\dots,J\}$. 
We assume
that these data are $\delta$-close to
FISD of some manifold ${\dtext (M,g,p)} \in  \overline \bM_{\dtext n}$,
that is, the ball $(B_e(r_0), g)$ and 
$\{(\lambda_j,\, \varphi_j|_{B_e(r_0)}) ;\, j=0,1,2,\dots,J\},$
where the error size parameter $\delta$ satisfies (\ref{29.5}).
\\
We are going to formulate a minimizaton algorithm that will be used to compute 
volumes of the sets (\ref{Ma set}).
We consider this minimizaton algorithm in the two cases, in the case when we have FISD without errors 
and the case when we have it with  errors.
\\
%
As we have finite data, we need to consider the projection of the solution
of the wave equation to finitely many eigenvectors, and we choose $\origjnull$
so that it is enough to use $\origjnull$ eigenvectors.
This requires that we have the data  $(\lambda_j,\, \varphi_j|_{B_e(r_0)})$ with $j=0,1,2,\dots,\origjnull$.
However, to consider minimization algorithms both for FISD with and without errors, we need to 
increase the amount of data and we will consider
 $(\lambda_j,\, \varphi_j|_{B_e(r_0)})$ with $j=0,1,2,\dots,\janul$,
where $\janul$ is chosen as follows:
In Definition \ref{def:3},
there are intervals $I_p\subset \R,\, p=0,1, \dots, P$ covering the spectrum of $M$ in 
$[0, \delta^{-1}+ \delta]$
each  $I_p$ containing a cluster of
 $n_p$ eigenvalues $\lambda_j$ and approximate eigenvalues $\lambda_j^a$. 
To consider these clusters of eigenvalues, let
$P_0$ be the smallest integer $P_0\leq P$  such that
 \beq\label{pre implications}
\{\lambda_0,\lambda_1,\dots, \lambda_{\origjnull}\}\subset \bigcup_{p=0}^{P_0} I_p
\eeq
and then choose $\janul$ such that  $\origjnull\leq \janul\leq J$ and
 \beq\label{implications}
j\leq \janul\implies \lambda_j\in   \bigcup_{p=0}^{P_0} I_p,\quad
j> \janul\implies \lambda_j\not \in   \bigcup_{p=0}^{P_0} I_p.
\eeq
We note that this happens with some $\janul$  satisfying (\ref{janolla}). 
We also observe that as  $\delta$  satisfies (\ref{29.5})  and $J$ satisfies (\ref{bound for J B}),
{\novtekst{and as $\Lambda_s\geq 1$, $\e_2<1$ and $n\ge 2$,}} we have
\beq\label{estimate for $J$}
J\geq J_0(\delta) =({ \RobCthree^{-1}  \delta^{-1}})^{n/2} \geq \janul.
\eeq

\subsubsection{Minimisation with FISD without errors}
\begin{theorem}\label{main sliceing} 
Let $\e_0,\e_1,\, \e_2$
satisfy (\ref{epsone defined}). 
There is $\derho_0(\e_0;  s, \Lambda_s) $
{\ztext depending only on $\e_0,$ $s,$ $\Lambda_s$, $n$, $R, D, i_0$ and $r_0$,}
with the following properties:
Let $\ga \le \derho_0(\e_0;  s, \Lambda_s) $. 
Assume that
$\origjnull$ satisfies (\ref{24Bnew}) and $\janul$ satisfies (\ref{janolla}), and
\ba 
u(x)=\sum_{j=0}^{\jnull} a_j\varphi_j(x) \in \H^s_{ \Lambda_s}(M),
\ea
Let $\,i \in \{L, \dots,  {\dtext L-1+N_1}(\sigma)\}$ and
$\a \in \mathcal A^{(i)}$.

Moreover, {\dtext assume that we are given} 
\beq\label{data 1}
(B(p, r_0),g|_{B(p, r_0)}),\quad
\left((\lambda_j,\  \varphi_j|_{B(p, r_0)})\right)_{j=0}^{\jnull}\quad\hbox{and}\quad  (a_j)_{j=0}^{\jnull}.
\eeq
{\ptext  The data \eqref{data 1} determine the set $\C^*_{\Rtext m}$ and the function {\wtext  $\L_{\vektor a}:C^*_{\Rtext m}\to \R$,} 
defined in \eqref{29.06.1} and \eqref{minimization E}, 
for which the 
minimizer of $L_{\vektor a}$ in $\C^*_{\Rtext m}$}
%
is {\dtext a sequence} {\wtext $(d_j)_{j=0}^{\jnull}=(d_j(\alpha,i))_{j=0}^{\jnull}\in \R^\jnull$} 
 such that
$
v(x)=\sum_{j=0}^{\jnull}  d_j^{} \varphi_j(x)
$ satisfies
\beq\label{18.1}
v \in  \H^s_{ (2 \SlavaCone(s, \derho) \Lambda_s)}(M) \cap {\Mtext \H^0_{2 \Lambda_s}(M)}, \quad
\|v-\chi_{M(\a,{\qtext -2\derho})}u\|_{L^2(M)}< \e_0,
\eeq

The {\dtextRB above bound  $\derho_0(\e_0; s, \Lambda_s)$ for $\gamma$ }
is defined in \eqref{pre 2 new condition to e1 e2}. 
\end{theorem}

 {\dtextRB Note that the sequence $(d_j)_{j=0}^{\jnull}$ is not unique and that the theorem states the existence of sequences satisfying (\ref{18.1}).}

{\Mtext The next subsections are devoted to the proof of Theorem \ref{main sliceing}.}
In sec. \ref{subsec 4.2 B}, \ref{subsec 4.3} and \ref{subsec 4.4} we keep the index $\,i \in \{L, \dots,  {\dtext L-1+N_1}(\sigma)\}$ fixed not referring to this.

\subsubsection{Finite dimensional projections}\label{subsec 4.2 B}

Next we introduce some special sets of the finite-dimensional functions.
\begin{definition} \label{function_classes}
Let $\vektor b=(b_j)_{j=0}^{\jnull} \in \R^{(\jnull+1)}$ and  $\F^*(\vektor b)$ be its Fourier coimage
\beq \label{extra1}
\F^*(\vektor b)=\sum_{j=0}^{\jnull} b_j\varphi_j \in L^2(M).
\eeq
For $a_1,  a_2>0$  the class of Fourier coefficients $\C_{\jnull,s}(a_1, a_2)$ is defined as
\beq \label{17.1}
\C_{\jnull,s}(a_1,a_2):=\{\vektor b\in \R^{(\jnull+1)}; \,\,
\sum_{j=0}^{\jnull}(1+\lambda_j^2)^{s}|b_j|^2\leq a_1^2,\quad {\sum_{j=0}^{\jnull} |b_j|^2\leq a_2^2}\}.
\eeq
For $w=W(v)$ being the solution to the  problem (\ref{12})
and
$\vektor b\in \R^{(\jnull+1)}$,  we denote
\beq\label{script W}
{\mathcal W}(\vektor b)=W(\F^*(\vektor b)) \in {\dtext C(\R;L^2(M)})
\eeq
and, for any $\e_* >0,\, \a \in \mathcal A^i$, we denote
\beq\label{y-8-B}
& &\C_{\jnull, s}( \e_*;\, a_1,  a_2, \alpha) 
\\ 
  \nonumber
& &=
\{\vektor b\in \C_{\jnull ,s}(a_1, a_2) :\,
\|W(\F^*(\vektor b))\|_{L^2({\qtext  \tilde\Gamma(z_\ell, \a_\ell )})} \le  \e_*, \;
 \forall \ell \in K_i\}.
\eeq
\end{definition}

\begin{lemma}\label{Lem:16.1} (i) Let $v\in  \H^s_{ \,\SlavaCone(s, \,\derho ) \Lambda_s}(M)$ 
and let $ P_{j'}$ be the orthoprojection
$
 P_{j'}v=\sum_{j=0}^{j'}\bra v,\varphi_j\cet_{L^2(M)}\, \varphi_j.
$
Then for any $\a \in \mathcal A^{(i)},\, \e_2>0,$ 
\beq\label{24B}
 \| P_{j'}v-v\|_{L^2(M)}\leq \frac{1}{{\Mtext 8}(\SlavaCsixagain+1)}\, \e_2, \quad
\hbox{if}
\,\,
{\Mtext j' \ge\origjnull\ge \widehat \jbasis({\Rtext \frac{\e_2}{8}};  \derho, \Lambda_s),}
\eeq
see (\ref{problem}) for $\SlavaCsixagain$ and (\ref{def SlavaCthree(s)}) for 
$\widehat \jbasis( \frac{\e_2}{8};  \derho, \Lambda_s).$

\noindent
(ii) Let $u \in  \H^s_{ \Lambda_s}(M)$
and $u_\a$ be  given by (\ref{def of tilde u}).
Let $\jnull$ satisfy (\ref{24Bnew})-(\ref{janolla}).
Then,
\beq
\label{eq: tilde v}
 v_\a=P_{\jnull} u_{\alpha} \in
 \F^*\left(\C_{\jnull, s}(\,  {\Mtext \frac 1{8}}\e_2;  \,{\Mtext \frac 14} \SlavaCone(s;\derho) \Lambda_s,\Lambda_s,  \alpha)\right).
\eeq
\end{lemma}

\noindent
{\bf Proof.}
(i) For $v = \sum_{j=0}^\infty b_j\varphi_j$, we have
$$
\| P_{j'}v-v\|^2_{L^2(M)}= \sum_{j >j'} |b_j|^2
\le 
 |\la_{j'}|^{-s} \SlavaCone(s; \derho)^{2} \Lambda_s^2.
$$
{\novtekst{Here, $\SlavaCone(s; \derho) $ is defined in {(\ref{20.4})   and 
 (\ref{26.5}) with $s=0$}}, and these  impy the estimate  (\ref{24B}).}

\noindent (ii) The  finite propagation speed of waves implies,
due to $u_\a|_{{\qtext M(\a,\derho)}}=0$,  
that
$W( u_\a)|_{\qtext \tilde \Gamma(z_\ell, \a_\ell)}=0$.
By Lemma \ref{lemma: tilde u exists} and (\ref{problem}) 
\begin{eqnarray}\label{e2 error A}
\| W( v_\a)\|_{L^2({\qtext \tilde \Gamma(z_\ell, \a_\ell)})}\leq
{\qtext \| W( v_\a-u_\a)\|_{L^2({\qtext \tilde \Gamma(z_\ell, \a_\ell)})}\leq}
\SlavaCsixagain \|v_\a - u_\a\|_{L^2(M)} \le {\Rtext \frac 18}\e_2.
\end{eqnarray}
Since $\|P_{\jnull}\|_{H^s(M)}=1$  for any $s$, the claim (i) of the lemma with 
$j'=\jnull$, (\ref{e2 error A})
together with  (\ref{y-4b}) prove (\ref{eq: tilde v}).
\hfill\Box\medskip


\begin{figure}
\begin {picture}(200,100)(80,0)
\put(100,0) {\includegraphics[height=5.0cm]{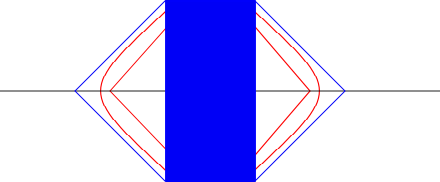}}
\put(153,53){$  \tilde \Sigma$}
\put(90,58){$\{0\}\times M$}
\put(204,53){$  \Sigma_{2\gamma}$}
\put(195,93){$  \D_0$}
\put(243,150){$\tilde \Gamma( \a)$}
\end {picture}
\caption{
\label{fig_4}
{\it In the proof of Theorem \ref{main sliceing}  we solve the minimization problem (\ref{minimization}) where we consider waves $U({x,t})=W(v)$ with initial data $(U(x,0),U_t(x,0))=(v(x),0)$,
$v=\sum_j b_j\varphi_j$, such that $U|_{\tilde \Gamma(\a)}$  is small. By using approximate unique continuation in the domain
$\D_0=\D(\alpha)$, we see that the wave $U$ is small in $   \Sigma_{2\gamma}=\Sigma(\alpha,2\gamma)$ and that $v(x)$
is small in $M(\a, -2\derho)$. In the proof we consider function $u_\alpha(x)$ that
is equal to $u(x)$ in $M \setminus M(\a, {\qtext  7\derho})$. Moreover, $u_\alpha$ 
 vanishes in the 
set $M(\a,\gamma)$ 
and thus the wave 
$U_\alpha=W(u_\alpha)$ produced by  the initial data $(u_\alpha,0)$ 
vanishes in the extended double cone
$ \tilde \Sigma=\Sigma(\alpha,-\gamma)$ and thus in $\tilde \Gamma(\a)$. Using these we see that  $v_\a=P_{\jnull} u_{\alpha}$  is close to the solution  $v^*$ of the minimization problem (\ref{minimization}).
}}
\end{figure}
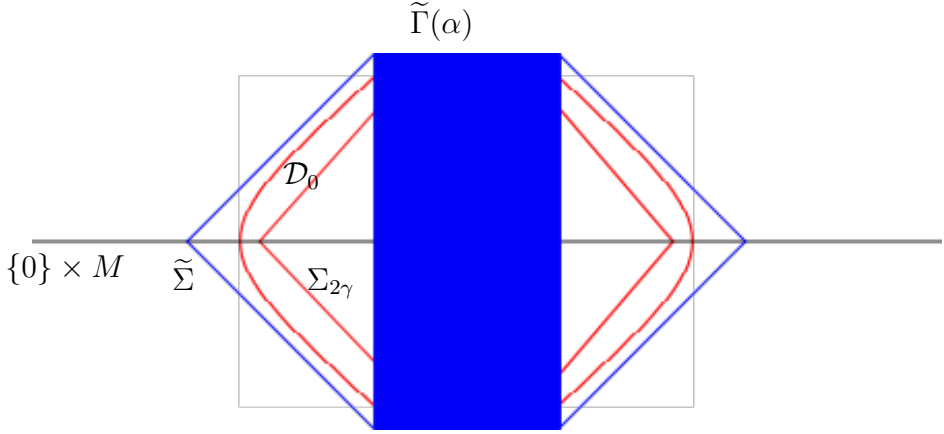

\begin{Remark} \label{rem:norm} The condition
$\|W(\F^*(\vektor b))\|_{L^2({\qtext  \tilde \Gamma(z_\ell, \a_\ell )})} \le  \e_*$, see (\ref{eq: tilde v}),
is equivalent to
\beq\label{Blago analogy2-1}
\big\|  (\sum_{j=0}^{\jnull}  b_j\cos(\sqrt{\lambda_jt})\,
 \varphi_j(x))\big \|_{L^2({\qtext  \tilde \Gamma(z_\ell, \a_\ell)})}\leq
 \e_*, \,\,   \ell \in K_i.
 \eeq
which  can be directly verified if we know
 $\{(\la_j,\,  \varphi_j|_{B(p, r_0)})\}_{j=0}^{\jnull}$. 
\end{Remark}

\subsubsection{Minimisation algorithm}\label{subsec 4.3}

Assume that we are  given $\vektor a=(a_j)_{j=1}^{\jnull} \in \R^{(\jnull+1)}$ {\dtext and denote}
$
u= \F^*(\vektor a) \in  \H^s_{\Lambda_s}(M). 
$
Our next goal is
to use  FISD to find a vector 
$\vektor b \in \C_{\jnull, s}(\SlavaCone(s, \ga)\Lambda_s, \,\Lambda_s)$
such that $\F^*(\vektor b)$ is close  to $\chi_{M(\a)} \F^*(\vektor a)$.
To achieve this goal  we will use a minimisation method.

 {\Rtext }
Let $\e_0,\e_1,\, \e_2$
satisfy (\ref{epsone defined}). 
Let $m\in \{1,2,4\}$ be a parameter we will use below, and
\beq \label{29.06.1}
 \U_{\Rtext m}:=\F^*\left(\C^*\right), \quad \hbox{where}\,\, 
\C^*_{\Rtext m}=\C_{\jnull, s}(  {\Mtext \frac 1{2m}}\e_2;\, {\Mtext \frac 1m}\SlavaCone(s, \ga)\Lambda_s, \Lambda_s, \alpha).
\eeq
\begin{definition} \label{def:minimization}
(i)\;A function $ v \in \U_{\Rtext m}$
is called an $\e_1 $-minimizer of the minimization problem
\beq\label{minimization}
\min_{h\in \U_{\Rtext m}}\L_u(h),\quad \hbox{where }\L_u(h)=
\|h-u\|^2_{L^2(M)},
\eeq
if $ v$ satisfies
\beq\label{minimization2}
\|v-u\|_{L^2(M)}\leq {\Rtext J_{min}(m)}+ 5 \Lambda_s\e_1  ,\quad  {\Rtext J_{min}(m)}:=\inf_{h\in \U_{\Rtext m}}\|h-u\|_{L^2(M)}^2.
\eeq
(ii)\; Equivalently, a vector
$\vektor b=(b_j)_{j=0}^{\jnull}\in \C^*_{\Rtext m}$  is an $\e_1 $-minimizer of the minimization problem
\beq\label{minimization E}
\min_{\vektor c\in \C^*_{\Rtext m}}\L_{\vektor a}(\vektor c),\quad \hbox{where }\L_{\vektor a}(\vektor c)=
\|\vektor c-\vektor a\|^2_{\R^{(\jnull+1)}},
\eeq
if
\beq\label{minimization E2}
\|\vektor b-\vektor a\|_{\R^{(\jnull+1)}}^2\leq {\Rtext J_{min}(m)}+ 5 \Lambda_s \e_1,\quad
{\Rtext J_{min}(m)}:=\inf_{\vektor c\in \C^*_{\Rtext m}}\|\vektor c-\vektor a\|_{\R^{(\jnull+1)}}^2.
\eeq
\end{definition} 
Observe that for $\vektor c \in \C_{\jnull,s}({\Rtext \frac 1m}\SlavaCone(s, \ga)\Lambda_s, \Lambda_s)$
we can check, using Remark \ref{rem:norm} with $\e_*= \e_2{\Mtext /2m}$, that $\vektor c\in \C^*_{\Rtext m}$ and thus
find $\vektor b$ which satisfies (\ref{minimization E2}).

Next we assume that, in addition to  $\e_2$ satisfying (\ref{19.1}), $ \derho$ satisfies 
\beq
\label{pre 2 new condition to e1 e2}
& & \derho\leq \gamma_0={\Mtext \derho_0(\e_0; s, \Lambda_s)=}  \SlavaCthreeagain \left(\frac{\e_1}{\Lambda_s} \right)^{1/(b(s))},
 \quad \hbox{with}\,\,\\
 & &\nonumber \SlavaCthreeagain(s)=\frac{ 1}{( 2 L \Slavacnonumber(s) )^{1/(2b(s))} },  
 \,\, \e_1 = \frac{\e_0^2}{10\Lambda_s},
\eeq
where $b(s)$ and $\Slavacnonumber(s)$ are defined in Lemma \ref{y-lemma 1}, b).
%
%
%
%
%

\begin{lemma}\label{lem: properties of minimizers 1}
Let  
$u\in  \H^s_{\Lambda_s}(M)$, and let
 $\e_0,\e_1,\e_2, \,\jnull,\,\ga$ satisfy (\ref{epsone defined}), (\ref{24Bnew})-(\ref{janolla})
and  (\ref{pre 2 new condition to e1 e2}). 

\noindent
(i) For {\Rtext $m\in \{1,2,4\}$ and} all $h\in  \U_{\Rtext m},$ we have
\beq \label{(H)}
\L_u(h) \geq \|u\|_{L^2(M(\a, -2\derho))}^2-2 \Lambda_s \e_1
+ \|h-u\|_{L^2(M\setminus M(\a, -2\derho))}^2.
\eeq

\noindent
(ii) The function $ v_\a$ defined by (\ref{eq: tilde v}),
(\ref{def of tilde u})
 satisfies
 $ v_\a \in \U_{\Rtext m}$ {\Mtext with $m=4$ and} and
\beq \label{18.11}
\L_u( v_\a)
\le  \|u\|^2_{L^2(M(\a, -2\derho))}+2 \Lambda_s  \e_1 +4 \e_1^2.
\eeq
{\Mtext Note that here $ v_\a \in \U_{\Rtext 4}\subset\U_{\Rtext 2}\subset \U_{\Rtext 1}$.}

\noindent
(iii) For all {\Rtext $m\in \{1,2,4\}$,} the function $v_\a\in \U_{\Rtext m}$ is an $\e_1-$minimiser,
\beq\label{minimal value J}
\L_u(v_\a)  \leq  {\Rtext J_{min}(m)}+ 5\Lambda_s \e_1 .
\eeq

\noindent
{\Rtext (iv) For all $m\in \{1,2,4\}$, we have 
\beq
\bigg |J_{min}(m)- \|u\|^2_{L^2(M(\a, -2\derho))}\bigg|\leq 2 \Lambda_s  \e_1 +4 \e_1^2.
\eeq}
\end{lemma}

\noindent {\bf Proof.}
(i)
We have, for
$h\in \U_{\Rtext m}$,
\ba
& &\|h-u\|_{L^2(M)}^2= \|h-u\|_{L^2(M(\a, -2\derho))}^2
+ \|h-u\|_{L^2(M\setminus M(\a, -2\derho))}^2 \\
&\geq &(\|u\|_{L^2(M(\a, -2\derho))}-\|h\|_{L^2(M(\a, -2\derho))})^2
+ \|h-u\|_{L^2(M\setminus M(\a, -2\derho))}^2.
\ea
Since $h \in  \U_{\Rtext m},$  (\ref{19.1}),
(\ref{19.2}) and (\ref{19.3}) imply that $\|h\|_{L^2(M(\a, -2\derho))} \le \e_1$.
Thus,
\ba
\|h-u\|_{L^2(M)}^2
&\geq & \|u\|_{L^2(M(\a,-2\derho))}^2-2\Lambda_s \e_1 +
\e_1^2
+ \|h-u\|_{L^2(M\setminus M(\a,-2\derho))}^2.\nonumber
\ea
(ii) With $ u_\a, \,v_\a $ defined by (\ref{def of tilde u}) and  (\ref{eq: tilde v}),
$ v_\a \in \U_{\Rtext 4},$  
\beq \label{19.4a}
& &\|u- v_\a\|^2_{L^2(M)} = \|u- v_\a\|^2_{L^2(M(\a, -2\derho))}+
\| u- v_\a\|^2_{L^2(M\setminus M(\a, -2 \derho))} \\ \nonumber
& & \le \|u\|_{L^2(M(\a, -2\derho))}^2 +2  \Lambda_s \e_2+
 \e_2^2
\\ \nonumber
& &\quad +2\|u- u_\a\|^2_{L^2(M\setminus M(\a, -2 \derho))}
+ 2\|  u_\a- v_\a\|^2_{L^2(M\setminus M(\a, -2 \derho))},
\eeq
where we use that $u- v_\a=(u- u_\a) +( u_\a - v_\a)$
and
$\|v_\a\|_{L^2(M(\a, -2 \ga))}^2\le  \e_2^2$, see Lemma \ref{Lem:16.1}.
Observe, that by (\ref{Sobolev}), (\ref{y-4b}) and (\ref{pre 2 new condition to e1 e2}),  
\beq \label{23.1}
\|u-u_\a\|^2_{L^2(M \setminus M(\a, -2\derho))} 
= \|u\|^2_{L^2(M(\a, {\qtext 7\derho}) \setminus M(\a, -2\derho))}
 \nonumber
\le \Slavacnonumber(s) \Lambda_s^2 L^{2}  \derho^{2b(s)} \le \frac12 \e_1^2.
\eeq
where $\Slavacnonumber(s)$ is defined in Lemma \ref{y-lemma 1}, b).
Using (\ref{24B}) and (\ref{y-4b}), we see that 
$ \|  u_\a-v_\a\|^2_{L^2(M\setminus M(\a-2 \derho))} \leq
 \e_2^2.$
Thus, inequality  (\ref{19.4a})  yields that
$$
\L_u(v_\a) \leq \|u\|^2_{L^2(M(\a, -2\derho))}+2 \Lambda_s  \e_2+
3  \e_2^2
+\e_1^2.
$$
 As  $ \e_2 \le \e_1$, see (\ref{30.06.7}), 
we get  (\ref{18.11}).
\\
\noindent
(iii)
The claims (i) and (ii) together with (\ref{epsone defined}) yield that
\ba
\L_u(v_\a)-{\Rtext J_{min}(m)}=\L_u(v_\a) -\min_{h\in  \U_{\Rtext m}}\L_u(h) \qquad\qquad\qquad\qquad\qquad\\
\leq \big( \|u\|^2_{L^2(M(\a, -2\derho))}+2 \Lambda_s \e_1 +
4 \e_1^2 \big)-
\big(\|u\|_{L^2(M(\a, -2\derho))}^2-2 \Lambda_s \e_1 \big)
\leq
5 \Lambda_s \e_1.
\ea
\noindent
(iv)
{\Rtext The claim (iv) follows from (i) and (ii).}
\hfill\Box\medskip

\begin{lemma}\label{lem: properties of minimizers 2}
Let {\Rtext $m\in\{1,2,4\}$ and}
$u\in  \H^s_{\Lambda_s}(M)$, 
 {\Mtext $\e_0,\e_1$, and $\e_2$ satisfy (\ref{19.1}) and  (\ref{epsone defined}),} $\jnull$ satisfies (\ref{24Bnew})-(\ref{janolla})
and $\ga$ satisfies (\ref{pre 2 new condition to e1 e2}).
Let $ v^*= \sum_{j=0}^{\jnull} b_j \varphi_j$ be any  $\e_1 $-minimizer of  the minimization problem (\ref{minimization}), with
$\vektor b \in \C^*_{\Rtext m}$.
Then
\beq\label{eq: v-chi u}
\| v^*-\chi_{(M \setminus M(\a, -2\derho))}u\|^2_{L^2(M)}
\leq \e_0^2.
\eeq
\end{lemma}

\noindent{\bf Proof.}
Since $v_\a$  satisfies  {\novtekst{ by (\ref{18.11})} and } (\ref{minimal value J}),
$$
\| v^*-u\|_{L^2(M)}^2 \le \| v_\a-u\|_{L^2(M)}^2 +5 \Lambda_s \e_1
\\ \nonumber
\leq
\|u\|_{L^2(M(\a, -2\derho))}^2 +7 \Lambda_s \e_1
+ 4 \e_1^2.
$$
Since $ v^*-u$ satisfies (\ref{(H)}), this inequality implies that
\beq \label{important}
 \| v^*-u\|_{L^2(M\setminus M(\a,-2\derho))}^2
\leq
9 \Lambda_s \e_1+ 4 \e_1^2.
\eeq
Since $v^*\in \U$, $w^*=W(v^*)$ satisfies  (\ref{Blago analogy2-1}) with
 $\e^*=\e_2$, where $\e_2$ satisfies (\ref{19.1}) and  (\ref{epsone defined}). 
 It then follows from Corollary \ref{E_2} 
 that
 $$
 \| v^*\|^2_{L^2(M(\a, -2\derho))}
 \le  \e_1^2.
 $$ Due to
(\ref{epsone defined}), this inequality
 together with  (\ref{important}), implies
 (\ref {eq: v-chi u}).
\hfill\Box\medskip\\

\noindent
{\bf Proof of Theorem \ref{main sliceing}}. Assume that
  $\vektor a:=(a_j)_{j=0}^{\jnull}$ satisfies the hypothesis.
\\
First determine $(b_j)_{j=0}^{\jnull}$ so that $v^*=\sum_{j=0}^{\jnull}b_j\varphi_j(x)$ is 
an $\e_1 $-minimizer of (\ref{minimization}), 
$v^* \in \C_{\jnull, s}({\Rtext \frac 1m}\SlavaCone(s, \derho) \Lambda_s, \Lambda_s)$ with $m=1$.
Then, by (\ref{eq: v-chi u}), 
$$
\|\chi_{M(\a, - 2\derho)}u-\sum_{j=0}^{\jnull}(a_j-b_j)\varphi_j
\|_{L^2(M)}< \e_0.
$$
Take $d_j=a_j-b_j$. Then 
$
v(x)= \sum_{j=0}^{\jnull} d_j \varphi_j(x) 
$ satisfies (\ref{18.1}).
\hfill\Box\medskip

\subsubsection{Minimisation with finite interior spectral  data with errors}\label{subsec 4.4}

{\Mtext
In this section we consider an approximate construction  when there 
is a $\delta$-error
 in FISD. {\novtekst
  We assume that we are given the ball $(B_e(r_0), g^a)$ and the pairs
 $(\lambda^a_j,\, \varphi^a_j|_{B_e(r_0)})$ with $j=0,1,2,\dots,J$.}

We assume
that these data are $\delta$-close to
ISD of some manifold ${\dtext (M,g,p)} \in  \overline \bM_{\dtext n}$
in the sense of
Definition \ref{def:3}
with intervals $I_p\subset \R,\, p=0,1, \dots, P$ covering the spectrum of $-\Delta_g$ in 
$[0, \delta^{-1}+ \delta]$.
We will use parameters $\origjnull,\janul\in \Z_+$ and $P_0\in \Z_+$ satisfying 
(\ref{24Bnew})-(\ref{janolla}), (\ref{pre implications}), and  (\ref{implications}).
Note that then  $\origjnull\le \janul\le J$ and that
below we will use  $(\lambda^a_j,\, \varphi^a_j|_{B_e(r_0)})$ with $j=0,1,2,\dots, \janul$.
Denote  $\mathcal J_p=\{j\in\Z_+;\ \lambda^a_j\in I_p\}$ and $n_p$ is the number of elements in
$\mathcal J_p$.

Then, for 
any $p$ there is $A^{p} \in O(n_p)$, $p=1,2,\dots,P_0$ such that, if
$j \in \mathcal J_p$ then 
\beq \label{5.21.12} 
\widetilde\varphi_j= \sum_{k \in \mathcal J_p} A^{p}_{jk} \varphi_k
\eeq
satisfies $\|\widetilde\varphi_j-\varphi_j^a\|_{L^2(M)}<\delta$, where $\varphi_k$ are the eigenfunctions of $\Delta_g$.
Note that $\sum_{p=0}^{P_0} n_p =\janul+1$.
 We use below the matrix $E\in O(\janul+1)$,
\beq \label{5.07.2}
E=[e_{jk} ]_{j, k=0}^{\janul},  \quad
e_{jk}=\bra \linetilde  \varphi_k, \varphi_j\cet_{L^2(M)}
\eeq
and note that $e_{jk}=0$  if $\la_j, \la_k$ do not lie in the same $I_p$.

Let $\vektor b =(b_0, b_1, \dots, b_{\janul}) \in \R^{{\janul}+1}$ then, for
$\vektor{ b}^a=E(\vektor b)$, $\vektor b^a =(b_0^a, b_1^a, \dots, b^a_{\janul})$ we have
\beq\label{tilde rep}
\sum_{j=0}^{{\janul}}   b_j^a \linetilde  \varphi_j(x)=\sum_{j=0}^{{\janul}} b_j  \varphi_j(x).
\eeq

Also, let $\omega_j$  be the center point of the interval $I_p$ containing
 $\lambda_j^a$ so that $|\lambda_j^a-\omega_j|<\delta$.

 The main goal of this section is to prove
 
\begin{theorem}\label{main sliceing with errors} Let 
$\e_0,\e_1,\e_2$  satisfy (\ref{epsone defined}).
{\Rtext Let   $\ga$ satisfy 
\eqref{pre 2 new condition to e1 e2}, 
$\origjnull$ satisfies (\ref{24Bnew}) and $\janul$ satisfy (\ref{janolla}),
$\delta$ satisfies (\ref{29.5}),
 and let $J$ satisfies
\beq \label{08.10.2017}
J_0(\delta)\leq J \leq 2^{n/2}\RobCthree^{n}   J_0(\delta) , \quad \hbox{where }J_0(\delta)=(2\RobCthree\delta)^{-n/2}.
\eeq}
Then  the following is valid:
\medskip

Let
 $z_1,\dots,z_{ {\dtext L-1+N_1}(\sigma)}\in B(p, r_0/4)$ be a $\sigma-$net.
Assume that $g^a|_{B_e(r_0)}$ and $((\lambda_j^{a},\varphi_j^{a} |_{B_e( r_0)}))_{j=0}^{\janul}$
is
 $\delta$-close to  FISD $g|_{B(p,r_0)}$  and
 $((\lambda_j,\varphi_j |_{B(p,r_0)}))_{j=0}^{\janul}$ of  a manifold
  ${\dtext (M,g,p)} \in  \overline{{\bM}_{\dtext n}}$. 
 Also, assume that ${\Rtext {{\tilde{\vektor{a}}}}}= ({\Rtext \tilde a_j})_{j=0}^{\Mtext \janul}$ satisfies
$
\, \sum_{j=0}^{\Mtext \janul} \bra \la_j^a \cet^s  |{\Rtext \tilde a_j}|^2 \le  \Lambda_s^2,
$
and 
$$
\tilde u^a(x)= \widetilde {\F}^*({\Rtext  {{\tilde{\vektor{a}}}}})=\sum_{j=0}^{\Mtext \janul} \tilde a_j  \widetilde\varphi_j(x),\quad\hbox{ for $x\in M.$}$$
Let   $\a\in \mathcal A^{(i)}$.

 {\dtext Assume that we are given 
\beq\label{data 2}
g^a|_{B_e(r_0)},\quad ((\lambda_j^{a},\varphi_j^{a} |_{B_e( r_0)}))_{j=0}^{\janul},\quad
\quad\hbox{and}\quad ({\Rtext \tilde a_j})_{j=0}^{\Mtext \janul}.
\eeq
{\wtext Let $\,i \in \{L, \dots,  {\dtext L-1+N_1}(\sigma)\}$ and
$\a \in \mathcal A^{(i)}$.}

{\ptext  The data \eqref{data 2} determine the set 
 $  \C^{a,*}_2$ and the function {\wtext $\L_{\vektor a}:\C^{a,*}_2\to \R$,} 
defined in \eqref{1.24.12 B} and \eqref{minimi NEW},
for which
the minimizer
of $\L_{\vektor a}$ in  $  \C^{a,*}_2$ is a sequence}}
%
$$
\vektor{  d}^a=\vektor{  d}^a(\a,i)={\wtext ({ d}^a_j(\a,i))_{j=0}^{\Mtext \janul}
\in \R^{\Mtext \janul}
,}
$$ 
such that
$
\widetilde v^a(x)= \linetilde {\F}^*(\vektor{  d}^a)=\sum_{j=0}^{\Mtext \janul}  {d}_j^a  \linetilde{\varphi}_j(x)
$, $x\in M$,
satisfies, cf. (\ref{18.1}),
\beq \label{28.1b}
 \widetilde v^a \in  \H^s_{ 2\SlavaCone(s, \derho) \Lambda_s}(M) \cap {\Mtext \H^0_{ 2\Lambda_s}(M)}, \quad
\|\widetilde v^a-\chi_{M(\a, -2\derho)}\tilde u^a\|_{L^2(M)}< \e_0.
\eeq 
\end{theorem}

\subsubsection{Proof of Theorem \ref{main sliceing with errors}     }\label{subsec 4.4a}

{\novtekst{The rest of this section will be devoted to the proof of Theorem 
\ref{main sliceing with errors}. 
Similar to (\ref{extra1}), 
we introduce }}
%
%
$$
\widetilde \F^*(\vektor b^a)= \sum_{j=0}^{{\Mtext \janul}}   b_j^a  
\,\linetilde  \varphi_j(x),
\quad x \in M;
$$
and
 wave-type functions
\beq\label{31.10.18}
& &w^a(x, t)={\W}^{a}(\vektor {  b}^a):=
   \sum_{j=0}^{{\Mtext \janul}} { b}^a_j \cos(\sqrt{\lambda^{a}_j} \,t)\,
 \varphi^{a}_j(x), \quad x \in B(p, r_0);
\\
\label{Blago analogy2-5}
 & &w(x, t)=\tilde \W\left(\vektor b^a\right)(x, t)=
W\left(\widetilde \F^*(\vektor b^a)\right)(x, t),\quad x \in M;
\\ \label{26.09.2017 cc}
& &\widetilde  w(x, t)=
\sum_{j=1}^{{\Mtext \janul}}  {  b}_j^a \cos \left(\sqrt{\omega_j}\,t \right)\,\linetilde  \varphi_j(x), \quad x \in M;
\\
\label{26.09.2017a}
& &\widetilde  w^a(x, t) =\linetilde{\W}^a(\vektor b^a):= 
\sum_{j=1}^{{\Mtext \janul}}  {  b}_j^a \cos \left(\sqrt{\la_j^a}\,t \right)\,\linetilde  \varphi_j(x), \quad x \in M,
\eeq
where we recall that $W$ is defined by $W(v)=w$ where $w$ satisfies (\ref{12}),
and
\beq \label{1.11.1}
& &\hspace{-15mm}\C^a_{{\Mtext \janul} , s}(a_1, a_2) =\{  \vektor {  b}^a \in \R^{({\Mtext \janul} +1)}: \,\,
 \sum_{j=0}^{\Mtext \janul} \bra \la_j^a \cet^s | {  b}_j^a|^2
\le   a_1^2,\quad 
{\sum_{j=0}^{\Mtext \janul}|{ b}^a_j|^2 \le a_2^2}
\}
\\
\nonumber
& &\hspace{-15mm}\C^a_{{\Mtext \janul}, s}(\e_*;a_1, a_2, \alpha)=
\{\vektor  {  b}^a \in \C^a_{{\Mtext \janul} , s}(a_1, a_2) ; \,
\|\W^a(\vektor  {b}^a) \|_{L^2({\qtext \tilde \Gamma(z_\ell, \a_\ell)})} \le  \e_*,\,  \ell\in K_i \}.
\eeq
We note that (see (\ref{tilde rep}) and (\ref{script W}))
\beq
\widetilde \F^*(\vektor b^a)= \F^*(E^{-1}\vektor b^a),\quad \tilde \W(\vektor b^a)=  \W(E^{-1}\vektor b^a).
\eeq

\begin{lemma}  \label{comparison}
 Let 
 $\vektor { b}^a \in  \C^a_{{\Mtext \janul} , s}(  \SlavaCone(\ga, s) \Lambda_s, \Lambda_s) $.
  If $\delta <1$ satisfies  (\ref{29.5}) 
then, 
\beq\label{Blago analogy 2} \nonumber
& & \|w-w^{a}\|_{L^2(B(p, r_0)\times (-2D, 2D))} \leq \frac1{\Rtext 4} \e_2 .
\eeq
\end{lemma}

\noindent
{\bf Proof.}
 Due 
 to (\ref{Sobolev estimate for eigenfunctions}) and \eqref{26.5}, 
 for $j, k \in \mathcal J_p$,
\beq \label{23.09.2017}
\hspace{-3mm}\big|\sqrt{\lambda^a_j}-\sqrt{\omega_k}\big| \le {\Mtext 2 \sqrt{\SlavaCfour}} \delta, \ \
\|\linetilde  \varphi_j-\varphi_j^a\|_{L^2(B(p, r_0))} \le \delta, \ \
\| \varphi_j^a\|_{L^2(B(p, r_0))} \le 2.
\eeq
{\Mtext
Using this, we obtain for $|t|\le 2D$ the following estimates. First, the Schwartz inequality implies that
\beq \label{28.7.0 M}
\| w^a(\cdot, t)-\linetilde  w^a(\cdot, t)  \|_{L^2({\qtext B(p,r_0)})}& \le&
(\sum_{j=0}^{{\Mtext \janul}} |b^a_j|)\delta
\\
&\le&
 ({\Mtext \janul}+1)^{1/2}( \sum_{j=0}^{{\Mtext \janul}} (b^a_j)^2)^{1/2}\delta
\nonumber \le 2 \Lambda_s ({\Mtext \janul})^{1/2}\delta.
\eeq
Also,
we see that 
\beq \label{28.7.0 MB}
\hspace{-1cm}\|\linetilde  w^a(\cdot, t)- \tilde w(\cdot, t)  \|^2_{L^2({\qtext B(p,r_0)})} &\le&
 \sum_{j=0}^{{\Mtext \janul}} (\cos (\sqrt{\la_j^a}\,t )-\cos (\sqrt{\omega_j}\,t ))^2 (b^a_j)^2\\
\nonumber &\le&   (2D)^2 \SlavaCfour \delta^2
 \Lambda_s^2 =  4D^2 \SlavaCfour \Lambda_s^2\delta^2
  .\eeq
We have
\ba
\widetilde  w(x, t)=
 \sum_{p=0}^P  \sum_{j,k \in \mathcal J_p} b_j^a\cos \left(\sqrt{\omega_k}\,t \right) A^{p}_{jk}\varphi_k(x)
= 
\sum_{k=0}^{{\Mtext \janul}} (\sum_{j \in \mathcal J_p} A^{p}_{jk} b_j^a) \cos \left(\sqrt{\omega_k}\,t \right) \varphi_k(x)
\ea
and
\ba
 w(x, t)=
\sum_{p=0}^P  \sum_{j,k \in \mathcal J_p} b_j^a\cos \left(\sqrt{\la_k}\,t \right) A^{p}_{jk}\varphi_k(x)
= 
\sum_{k=0}^{{\Mtext \janul}} (\sum_{j \in \mathcal J_p} A^{p}_{jk} b_j^a) \cos \left(\sqrt{\la_k}\,t \right) \varphi_k(x),
\ea
and as $A^p$ are orthogonal matrices and $|\sqrt{\la_k}-\sqrt{\omega_k}\big| \le {\Mtext 2 \SlavaCfour^{1/2}} \delta,$  we see similarly to  (\ref{28.7.0 M}) and (\ref{28.7.0 MB})
\beq \label{28.7.0 MM}
\hspace{-1cm}\|\linetilde  w(\cdot, t)- w(\cdot, t)  \|^2_{L^2({\qtext B(p,r_0)})} &\le&  4D^2 \SlavaCfour \Lambda_s^2\delta^2
  .\eeq
Combining the above estimates with 
$\delta<\widehat \delta_0(\e_2, \ga, \janul, \Lambda_s)=  \Slavacsevenagain\,  \frac 1{\janul}\,  
\frac{\e_2}{\Lambda_s} $ and $\Slavacsevenagain=\min(\RobCthree^{-1} , \frac {(1+\SlavaCfour)^{-1/2}}{100(1+D)^{3/2}L}),$
%
we obtain the claim.
}
\hfill\Box\medskip
\\
By Definition \ref{function_classes} we have 
\beq \label{norm-comparison}
& &E\,\C_{{\Mtext \janul}, s}\left(
{\Rtext \frac 12}a_1, a_2  \right) \subset \C^a_{{\Mtext \janul}, s}\left( a_1, a_2  \right)
\subset E\,\C_{{\Mtext \janul}, s}\left(
{\Rtext 2}a_1, a_2  \right).
\eeq
Note that the $\ell^2$-norms of the sequences $(b^a_j)_{j=1}^{\janul}$ do not depend on eigenvalues
and, therefore, the same holds for the exact and approximate data. Also, the $\ell^2$-norms are invariant 
with respect to the operations involving orthogonal matrixes.
{\mtext

Definitions of the sets of sequences in (\ref{y-8-B}) and (\ref{1.11.1}), Lemma \ref{comparison}
and formula (\ref{norm-comparison}) imply that
\beq
\\
\nonumber
\hspace{-1cm}E\C_{{\Mtext \janul}, s}(\e_*-\frac 14 \e_2;\frac 12 a_1, a_2, \alpha)\subset \C^a_{{\Mtext \janul}, s}(\e_*;a_1, a_2, \alpha)
\subset  E\C_{{\Mtext \janul}, s}(\e_*+\frac 14\e_2;2a_1, a_2, \alpha)\hspace{-2cm}
\eeq
Let us use $\e_*=\frac 1{2} \e_2$ and define
\beq \label{1.24.12 B}
 \C^{a,*}_m=  
 \C^a_{\janul, s}\left(\frac 1{m} \e_m, \frac 12  \SlavaCone(\ga, s) \Lambda_s, \Lambda_s, \a\right),\quad m\in\{1,2,4\}
\eeq
Using the notations in (\ref{29.06.1}), we see that
\beq \label{29.06.1B}
E \C^*_{\Rtext 4}\subset  \C^{a,*}_2\subset E\C^*_{\Rtext 1}.
\eeq

Consider the
quadratic function $\L_{\vektor a}:\R^{\janul+1}\to \R$,
\ba
\L_{\vektor a}(\vektor c)=
\|\vektor c-\vektor a\|^2_{\R^{(\janul+1)}},\quad \L_{E\vektor a}(\vektor c)=
\|\vektor c-E\vektor a\|^2_{\R^{(\janul+1)}}.
\ea
%
cf. (\ref{minimization E}). Note that $\L_{\vektor a}(\vektor c)=\L_{E\vektor a}(E\vektor c)$.
Let $\vektor b^{*}\in  \C^{*}_4$ and
$\vektor b^{a,*}\in  \C^{a,*}_2$ be minimizers of $\L_{\vektor {a}}$ and  $\L_{E\vektor {a}}$, respectively, that is
\beq\label{minimi *}
& &\L_{\vektor {a}}(\vektor b^{*})=\min_{\vektor b\in  \C^{*}_4} \L_{\vektor {a}}(\vektor b) =:J_{min}(4),\
\eeq
and
\beq\label{minimi NEW} 
& &\L_{E\vektor {a}}(\vektor b^{a,*})=\min_{\vektor b^a\in  \C^{a,*}_2 } \L_{E\vektor {a}}(\vektor b^a)=:  J_{min}^{a}(2).
\eeq
Note that  we do not anymore consider $\e_1$-minimizers, but the minimizers.
Since
$ \C^{*}$ and $ \C^{a,*}$ are bounded and closed set in
$\R^{\janul +1}$ such  minimizers exist.
When $\e_1<  \Lambda_s/8$,
Lemma
\ref{lem: properties of minimizers 1} (iv) implies 
\beq \label{4 aj 1 minimi}
|J_{min}(4)-J_{min}(1)|\leq 2(2 \Lambda_s  \e_1 +4 \e_1^2)<5\Lambda_s  \e_1 
\eeq
{\novtekst Using  (\ref{29.06.1B}), (\ref{4 aj 1 minimi}),}
 and the fact that $E$ is an isometry, we see that
$$
{\novtekst J_{min}(1)\leq J_{min}^{a}(2)\leq J_{min}(4),}
\quad \hbox{and}\quad J_{min}^{a}(2)
\leq J_{min}(1)+5\Lambda_s  \e_1.
$$
These implies that {\novtekst the minimizer $\vektor b^{a,*}$ of function $\L_{E\vektor {a}}$ in the set  $\C^{a,*}_2$ satisfies $\vektor b^{a,*}\in   \C^{a,*}_2\subset E \C^*_{\Rtext 1}$ and so we have that }
$\tilde{\vektor b}^*=E^{-1}\vektor b^{a,*}$
is an $\e_1$-minimizer of $\L_{\vektor {a}}$ in the class $\C^*_{\Rtext 1}$. {\novtekst
We denote $\tilde{\vektor b}^*=( \tilde{ b}^*_j )_{j=1}^{\Mtext \janul}.$}
Let ${\Rtext  {\vektor{a}} =E^{-1} {{\tilde{\vektor{a}}}}}$ so that
$$
 {\F}^*(\vektor{a})=\sum_{j=0}^{\Mtext \janul} {\Rtext  a}_j  \varphi_j(x)=u(x).$$
Then, by applying Lemma \ref{lem: properties of minimizers 2} 
{\novtekst we see that  $ v^*= \sum_{j=0}^{\janul} \tilde{ b}^*_j \varphi_j$  satisfies} (\ref{eq: v-chi u}).
Then, choosing $d_j^a=\tilde a_j-b^{a,*}_j, \, j=0, 1, \dots, \janul$, we see that $\tilde v^a=\sum_{j=0}^{\janul } d_j^a \linetilde  \varphi_j$ satisfies (\ref{28.1b}).
This proves Theorem \ref{main sliceing with errors}.
\hfill\Box\medskip

\begin{Remark} \label{rem:25.09 A}
{\dtext Similarly to Remark \ref{Rem: size of J} and Theorem \ref{main sliceing with errors},
we see that if
the collection of
 $g^a|_{B_e(r_0)}$ and $((\lambda_j^{a},\varphi_j^{a} |_{B_e( r_0)}))_{j=0}^J$
 is
  $\delta$-close to  FISD of  a manifold
  ${\dtext (M,g,p)} \in  \overline{{\bM}_{\dtext n}}$ then without loss of generality, we can assume that $J$ satisfies (\ref{08.10.2017}). Indeed, the eigenvalues $\lambda_j$  with index   $j>  2^{n/2}\RobCthree^{n}   J_0(\delta)$
 are not used in the proof of Theorem \ref{main sliceing with errors}.}
 \end{Remark}  
}}

\section{Construction of the approximate interior distance maps.}

\subsection{Volume estimates}

Our next goal is to approximately evaluate
the volume of $M(\alpha)$, see (\ref{Ma set}) with $b=0$. 
\begin{lemma}\label{lemma: main slicing}
{\novtekst
There are uniform constants $\e_0^* >0,\, \SlavaCnine(s)>1,$
{\ztext depending only on  $s$, $n$, $R, D, i_0$ and $r_0$,}
 such that the following holds:

Let $\e_0 \le \e_0^*$.  
  Let  $\e_1,\,\e_2$ be defined by (\ref{epsone defined}) while 
  $\ga_0, \jnull$ be defined by (\ref{pre 2 new condition to e1 e2})
and (\ref{24Bnew})--(\ref {janolla}). 
Assume that  {\dtext we are given} $\big(g^a|_{B_e(r_0)};\, 
\{(\la_j^a, \varphi^a_j|_{B_e(r_0)}) \}_{j=0}^{J}\big)$ 
that is
 $\delta$-close to FISD of
${\dtext (M,g,p)} \in  \overline{\bM_{\dtext n}}.$
Here $J$ satisfies by \eqref{08.10.2017}. 

Let also $i \in \{ L, \dots,  {\dtext L-1+N_1}\},$ 
where $N_1=N_1(\sigma)$ is defined {\dtext as in}
Lemma 
\ref{ga-separation} (ii). {\dtext Assume that 
\beq\label{tau and sigma}
 \sigma \leq {\dtext \tau_0}/2
\eeq  
where $\tau_0$ is} defined in Proposition
\ref{uniform_covering} and let $\a \in \mathcal A^{(i)}$
satisfy (\ref{06.07.1}).

  Then
we can compute an approximate volume, $\vol^a (\M(\alpha))$, of   {\dtext the set $M(\a)$ that satisfies}
\beq \label{2.11.15}
|\vol^a (\M(\alpha))-\vol(\M(\alpha))|\le \SlavaCnine(s) \e_0.
\eeq}
\end{lemma}

\noindent{\bf Proof.} Recall that
\beq \label{07.06.2}
\varphi_0(x)=\vol(M)^{-1/2}, \quad \mathcal F(\varphi_0)=(1, 0, 0, \dots), 
\quad \|\varphi_0\|_{H^s} =1 \,\, \hbox{for} \,\,
s >0.
\eeq
The  interval $I_0=(a_0, b_0)$ in Definition \ref{def:3} contains only
 $\la_0=0$. Thus $ \varphi_0^a|_{B_e(r_0)}$
is a $\delta$-close to $ \varphi_0|_{B_e(r_0)}=\tilde \varphi_0|_{B_e(r_0)}$. 
These allow us to evaluate $\vol^a(M)$ so that $|\vol^a(M) -\vol(M)|< C \delta$.
Using Theorem \ref{main sliceing with errors}
we evaluate the 
{\Mtext Fourier coefficients $( d^a_j)_{j=0}^{\jnull}$  }
of $v^a(x)$ which satisfies (\ref{28.1b}) with
$\tilde u= \varphi_0$.
Let
\beq \label{22.07.1}
\vol^a( M(\a))=\vol^a(M) \, \bigg(\sum_{j=0}^{\jnull} ( d^a_j)^2\bigg)^{1/2}
\eeq
Then, by (\ref{28.1b}),
\beq\label{volume estimate2}
|  \vol^a( M(\a)) -\vol(\M(\a, -2 \derho))   | \le C(\e_0+\delta).
\eeq
Since $|\vol( M(\a))- \vol({ M(\a,-2\derho)}) |< \SlavaCnonumber L \ga$ (cf.\ Lemma \ref{y-lemma 1}), 
(\ref{volume estimate2}) implies
estimate (\ref{2.11.15}),  if $\e_0 \le \e_0^*$ with some uniform constants $\SlavaCnine(s)$ and $\e_0^*$.
Here $\e_0^*$ is defined so that $\delta <\e_0, \, \ga <\ga_0$ for 
$\e_0 <\e_0^*$, see (\ref{29.5}), (\ref{pre 2 new condition to e1 e2}).
\hfill\Box\medskip

Next we  use FISD with errors to approximately find the distances from various
points $x \in M$ to points $z \in B(p, r_0/4)$. The {\dtext main} tool is to
approximately
find the volumes of subdomains of $M$ obtained by 
the slicing procedure.

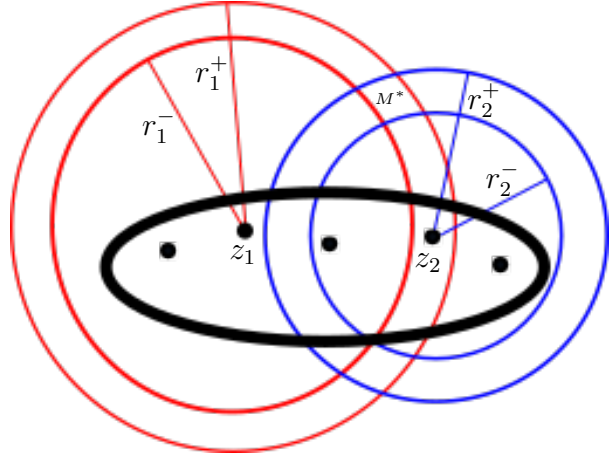
\begin{figure}
\begin {picture}(200,90)(80,0)
  \put(170,0) {\includegraphics[height=6.0cm]{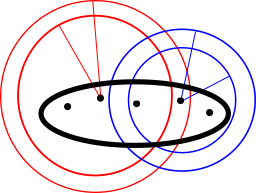} }

\put(253,73){$z_1$}

\put(240,140){$r_1^+$}
\put(220,120){$r_1^-$}

\put(308,132){\tiny $M^{*}$}

\put(323,70){$z_2$}
\put(343,130){$r_2^+$}
\put(350,100){$r_2^-$}

\end {picture}
\caption{
\label{fig_6}
{\it Slicing of the manifold: 
The intersection of ``slices'' $B(z_1,r_1^+)\setminus B(z_1,r_1^-)$ and  $B(z_2,r_2^+)\setminus B(z_2,r_2^-)$, where $z_1,z_2\in B(p,r_0)$,
is the set $M^*=(B(z_1,r_1^+)\setminus B(z_1,r_1^-))\cap (B(z_2,r_2^+)\setminus B(z_2,r_2^-))$. When $r_i^\pm=\beta_i\pm 2\sigma$  and $\beta=(\beta_1,\beta_2,0,0,\dots,0)$
we have $M^*=M^*(\beta)$, see (\ref{21}). We consider analogous indexes $\beta$ containing $L$ non-zero elements and the sets $M^*(\beta)$ for which
 the approximate volume $\vol^a(M^*(\beta))$ is large enough. Using those we choose  the set $\mathcal B$ of  the ``admissible'' indexes $\beta$ and consider the collection $X$ that is the union of the points $x_\beta$ chosen from  the sets $M^*(\beta)$ with $\beta\in \mathcal B$,
 and the points in a {maximal} $\sigma-$net in  $B(p, r_0/2)$, see (\ref{set X}). The set $X$ can be considered as a discrete approximation of the manifold $M$. The intersection of slices is also used to construct a distance function $d_X$ on the discrete set $X$ which makes $(X,d_X)$ an approximation of the manifold $(M,d_M)$. \vspace{-3mm}}}
\end{figure}

For $i \in \{L, \dots,  {\dtext L-1+N_1}(\sigma)\}$ and $\beta\in \mathcal A^{(i)}$,
 $M(\beta)$ are the domains defined in (\ref{Ma set}) with $\a$ replaced by $\beta$.
We consider the intersection of slices,
\beq \label {21}
\M^*_{ (i)}(\beta) = \bigcap_{ \ell\in K_i}
 (B(z_\ell, \beta_\ell+2 \sigma)\setminus B(z_\ell, \beta_\ell-2 \sigma))
\\ \nonumber
=\Big(\bigcap_{  \ell \in K_i}
  B(z_\ell, \beta_\ell+2 \sigma)\Big)\cap  \Big(\bigcup_{  \ell \in K_i}
 B(z_\ell, \beta_\ell-2 \sigma) \Big)^c.
\eeq
Here for $\Omega \subset M,\, \Omega^c =M \setminus \Omega$. 
%
%
Note that
\beq   \label{4.11.4}
& &\vol\Big( (\bigcap_{  \ell \in K_i}
 \Omega_\ell) \bigcap  \widetilde \Omega^c \Big)=
\sum_{  \ell \in K_i}
 \vol(\Omega_\ell \cup \widetilde \Omega) \\ \nonumber
& & - \sum_{\ell \neq \ell'=1}^n \vol(\Omega_\ell \cup \Omega_{\ell'}\cup \widetilde \Omega)
+ \dots +(-1)^{(L+1)} \vol \Big( (\bigcup_{ \ell \in K_i}
 \Omega_\ell ) \cup \widetilde \Omega \Big)
- \vol(\widetilde \Omega).
\eeq
By (\ref{21}),  $M^*_{ (i)}(\beta)$ has form (\ref{4.11.4})
with $$\Omega_\ell= B(z_\ell, \beta_\ell+2 \sigma), \, \widetilde \Omega=
 \bigcup_{  \ell \in K_i} B(z_\ell, \beta_\ell-2 \sigma).$$
For any  $\a_1,\a_2\in \mathcal A^{(i)}$ we have
\ba
\M(\alpha_1) \cup \M(\a_2) = \M(\a_m), \quad
 \hbox {where} \,\, (\a_m)_\ell = \max ((\alpha_1)_\ell, (\a_2)_\ell).
\ea
Therefore all terms in (\ref{4.11.4}) are of form $\vol(M(\a))$ for some 
$\a\in \mathcal A^{(i)}$. Thus,
using Lemma \ref{lemma: main slicing}, we can approximately compute each term
of (\ref{4.11.4}) with error $\SlavaCnine \e_0$. Since there are $2^L+1$ terms in (\ref{4.11.4}),
we
obtain the following result.
\begin{lemma} \label{volumes}
{\Mtext Under the conditions of Lemma \ref{lemma: main slicing},
there exists $\e_4(n, R, D, i_0, r_0) >0$ and $\Slavacten \in (0,1)$,
{\ztext depending only on $n$, $R, D, i_0$ and $r_0$,}
 with the following property:

Let  $0<\e_4 <\e_4(n, R, D, i_0, r_0)$.}   
It is possible to evaluate approximate volumes,
$\,\vol^a(M^*_{ (i)}(\beta))$, of the sets
$M^*_{ (i)}(\beta)$ of form (\ref{21}). Moreover, 
\beq \label{4.11.15}
\bigg| \vol^a(M^*_{ (i)}(\beta))- \vol(M^*_{ (i)}(\beta)) \bigg| \le \e_4,
\quad \hbox{if}\,\, \e_0 \le \Slavacten \e_4.
\eeq
\end{lemma}

\subsection{Distance functions approximation} \label{Distance functions approximation}
A function $r(\cdot)\in C^{0, 1}(B(p, r_0/4))$ is an {\it interior distance function}  
if there is $x\in M$ such that
$
r(z) =r_x(z) = d(x,z),$ for any $ z \in B(p, r_0/4).
$

The interior distance functions determine the interior distance map
\bfo
R_M: \,(M, g) \to L^{\infty}(B(p, r_0/4)), \quad
R_M(x) = r_{x}(\cdot).
\efo
The map $R_M$ or, more precisely, its image
\beq \label{7.11.1}
R_M(M):= \{ r_{x}(\cdot),\,\, x \in M \}
\subset L^\infty(B(p, r_0/4)),
\eeq
may be used to reconstruct $(\M,g)$.  Namely, in   \cite{Ku1},
\cite{KaKuLa}
it was shown  how to reconstruct
$(N,\, g|_N)$, where
$
N=M \setminus B(p, r_0/50),
$
from the knowledge of {\it boundary distance functions}
\beq \label{11.10.2017}
R_N(N)
= \{r^N_x(\cdot) \in L^\infty(\p N);  x\in N\}, 
\quad
r_x^{ N}(z)=d_N(x,z),
\eeq
where $d_N$ is the distance in $N$.
Later, in Section \ref{6.1} we show that 
a Hausdorff approximation $R^*_M$ to $R_M(M)$ makes it possible to
construct an approximation $R^*_N$ to $R_{ N}(N)$.

Thus, our next goal is
to construct a desired approximation $R_M^*$.
 To this end, we use the volume approximations
of the previous subsection.

First, for $z,z'\in B(p,r_0/2)$, we define  an approximate distance $d^a(z, z')$
using the metric $g^a$.
Then Definition \ref{def:3} (iv) together  with convexity of $B(p, r_0)$, 
see (\ref{4.09.2017}),
imply that
\beq \label{3.21.12}
|d^a(z, z') -d(z,z')| \le \SlavaCeleven \sigma, \quad \hbox{if}\,\, \delta < \SlavaCnineteen \sigma.
\eeq
{\dtextRB Recall that above we have used a parameter $\sigma>0$ 
 which  satisfy $\sigma\le \tau_0/2$ 
and we have chosen points $\{z_1,\dots,z_{{\dtext L-1+N_1(\sigma)}}\}\subset B(p, r_0/4)$
such that $\{z_1,\dots,z_{L-1}\}$ is a $\tau_0$-net in $B(p, r_0/4) $. Moreover, the set $\{z_L,\dots,z_{{\dtext L-1+N_1(\sigma)}}\}$
is a maximal $\sigma$-separated set in
$B(p, r_0/4) $, see (\ref{19.07.1}).}

For any $i\in \{ L, \dots,  {\Mtext {\dtext {\dtext L-1+N_1(\sigma)}}}\}$
 and $\beta=(\beta_\ell)_{\ell=1}^{\dtext L-1+N_1(\sigma)}\in \R^{ {\Mtext {\dtext {\dtext L-1+N_1(\sigma)}}}},$ $r_0/8 < \beta_\ell < 2D$,
 cf. (\ref{06.07.1}), we define 
 $\T^{(i)}(\beta)=\beta^{(i)}\in \R^{ {\Mtext {\dtext {\dtext L-1+N_1(\sigma)}}}},$ where
 \beq
\nonumber
&&\beta^{(i)}_\ell=\beta_\ell,\hbox{ if $\ell  \in K_i,$} 
\\
\label{Ai cut off}
& &\beta^{(i)}_\ell=0,\hbox{ if $\ell \not \in K_i.$} 
\eeq
 Then, $\T^{(i)}(\beta) \in \mathcal A^{(i)}.$

Observe that, for any $x \in M \setminus B(p, 3r_0/8+\sigma)$ and $\ell=1, \dots,{\Mtext {\dtext {\dtext L-1+N_1(\sigma)}}}$
there is
$\beta_\ell(x)\in \sigma \Z_+$
such that
 $\beta_\ell(x)-\sigma\leq d_M(x,z_\ell)\leq \beta_\ell(x)+\sigma.$ Therefore,
$B(x, \sigma) \subset B(z_\ell, \beta_\ell(x) +2 \sigma) \setminus B(z_\ell, \beta_\ell(x)-2\sigma)$,
so that, due to \eqref{12.06.13},
$$\vol\left(B(z_\ell, \beta_\ell(x) +2 \sigma) \setminus B(z_\ell, \beta_\ell(x)-2\sigma\right) 
\ge \frac 1{\RobCone} \sigma ^n.$$
Taking into account this inequality together with (\ref{4.11.15}) we require
\beq \label{4.11.20}
\e_4  \le  \frac 1{4\RobCone}\sigma ^n.
\eeq
Thus, for $i\in \{ L, \dots, {\Mtext {\dtext {\dtext L-1+N_1(\sigma)}}}\}$, {\dtext  the 
{\dtext volume and the approximate volume}
of the set
 $M^*_{ (i)}\left(\beta^{(i)}(x)\right)$, $\beta^{(i)}(x)=\T^{(i)}(\beta(x))$ satisfy}
$$
\vol\big(M^*_{ (i)}(\beta^{(i)}(x))\big) \ge 4 \e_4, \quad 
\vol^a\big(M^*_{ (i)}(\beta^{(i)}(x))\big) \ge 3 \e_4,
$$
where we use (\ref{4.11.15}).
{\dtextRB The above considerations motivate the following definition. 
In order to use only finitely many indexes $\beta$, in the following we are going to 
 consider  $\beta=(\beta_i)_{\ell=1}^{\dtext L-1+N_1(\sigma)}$ where $\beta_i \in \sigma\mathbb Z_+$, $\beta_i\leq 2D$.}

\begin{definition} \label{admissible}
Let $\beta =(\beta_\ell)_{\ell=1}^{\dtext L-1+N_1(\sigma)}\in \sigma\mathbb Z_+^{{\Mtext {\dtext {\dtext L-1+N_1(\sigma)}}}}\subset \mathbb R_+^{{\Mtext {\dtext {\dtext L-1+N_1(\sigma)}}}}
$. Such sequence
 $\beta$ is called {\rm admissible}, if  $r_0/8\le \beta_\ell \le 2D$ and 
 for all indexes
 $i\in \{L,\dots, {\Mtext {\dtext {\dtext L-1+N_1(\sigma)}}}\}$,
 the modified index
  $\,\beta^{(i)}=\T^{(i)} (\beta)\in \mathcal A^{(i)}$ satisfies
\beq \label{4.11.19}
\vol^a(M^*_{(i)}(\beta^{(i)})) \ge 3\e_4.
\eeq
We define the set $\mathcal B = 
\{\beta\in \sigma\mathbb Z_+^{  {\Mtext {\dtext L-1+N_1(\sigma)}}};\, \beta \mbox{ is admissible}\}$.
\end{definition}

\begin{lemma} \label{lem: admissible betas}
For any $x \in M\setminus B(p, 3r_0/8+\sigma)$, there exists an admissible
$\beta \in  \sigma\mathbb Z_+^{  {\Mtext {\dtext {\dtext L-1+N_1(\sigma)}}}}$
such that
$
|d(x, z_\ell)-\beta_\ell| \le 2 \sigma, $ for $ \ell \in \{1,2,\dots, {\Mtext {\dtext {\dtext L-1+N_1(\sigma)}}}\}.
$

Conversely, there is $\SlavaCtwelve>0$
{\ztext depending only on $n$, $R, D, i_0$ and $r_0$,}
 such that, if $\beta$ is admissible, then
 there is ${\ftext x=x_\beta} \in M \setminus B(p, 3r_0/8- \SlavaCtwelve \sigma)$ such that,
for all
$\ell\in \{1,2,\dots, {\Mtext {\dtext {\dtext L-1+N_1(\sigma)}}}\}$, we have
\beq\label{DD}
\,|\beta_{\ell}-d(x,z_\ell)|\leq  \SlavaCtwelve\sigma.
\eeq
 \end{lemma}

\noindent{\bf Proof.}
 The first statement follows from considerations before Definition \ref{admissible}.

On the other hand, assume that $\beta\in \mathcal B$. Then
equations (\ref{4.11.15}) and \eqref{4.11.19} guarantee  that, for any $i\in \{L,\dots, {\Mtext {\dtext {\dtext L-1+N_1(\sigma)}}}\}$,
 there is $x_i   \in M^*_{ (i)}(\mathcal T^{(i)}(\beta))$.} {Moreover, we have 
$
|d(x_i, z_\ell) -\beta_\ell| \le 2 \sigma,$  for  $\ell\in \{1, \dots, L-1\}\cup\{ i\}.$
Moreover, in view of (\ref{tau coordinates}), for $j, k \in \{L, \dots,  {\Mtext {\dtext {\dtext L-1+N_1(\sigma)}}}\}$,
\ba
d_M(x_j, x_k) \le \RobCsix |H^L(x_j)-H^L(x_k)| \leq 4 \RobCsix \sqrt L \sigma.
\ea
{\novtekst Defining
\beq \label{distance-appr}
\SlavaCtwelve= 4 \RobCsix \sqrt L+3,
\eeq
and taking $x=x_{\novtekst i_1} $  with {\novtekst arbitrary $i_1$}, we see that
$x \in  M \setminus B(p, 3r_0/8-\SlavaCtwelve \sigma)$ and that  (\ref{DD})
is satisfied.}
\hfill\Box\medskip

{\dtext For the points $\{z_\ell: \ell \in \{1,\dots, {\dtext {\dtext L-1+N_1(\sigma)}}\}\subset B(p, r_0/4)$, let 
$$V_\ell =\{y\in B(p, r_0/4):\ \hbox{$z_\ell$ is the unique closest point to $y$ in $\{
z_{\ell'}\}\,$}\} ,$$
where $\ell=1, \dots,{\dtext {\dtext L-1+N_1(\sigma)}},$  be the corresponding  Voronoi region.}
With any  $\beta \in \mathcal B$ we then associate
a piecewise constant function   {\dtext $r_\beta\in L^\infty(B(p, r_0/4))$ by defining 
$r_\beta(z)= \beta_\ell,$ for $\,z \in V_\ell$.} Clearly,
\beq\label{11.10.2017b}
d_{L^\infty(M)}(r_\beta,\, r_x) \le \SlavaCthirteen \sigma, \quad \SlavaCthirteen=\SlavaCtwelve+2.
\eeq
Let
$$
R^*_{M, >}=\{r_\beta(\cdot):\ \beta\in \mathcal B\ \} \subset  L^\infty(B(p, r_0/4)).
$$
 Choose a {\Mtext maximal} $\sigma-$net $\{x_1, \dots, x_{ N_2(\sigma)}\} \subset B(p, r_0/2)$
by adding to $z_1, \dots, z_{ {\Mtext N_0}(\sigma)}$ a $\sigma-$net 
$z_{{\Mtext N_0}(\sigma)+1}, \dots, z_{N_2(\sigma)}$   in $B(p,r_0/2) \setminus B(p, r_0/4)$.
{\Mtext Again, using Lemma \ref{ga-separation},
we see that  $ N_2(\sigma) \le \MattiCseven \sigma^{-n}$.}
Next we define
\bfo
r_k(z)=d^a(x_k, z_\ell), \quad  \hbox{for}\ z \in V_\ell, \quad k=1, \dots,  N_2(\sigma),\,\,
\ell=1,\dots,  {\dtext L-1+N_1(\sigma)};
\efo
\beq \label{8.07.1 B} \nonumber
& &R^*_{M,<}=\{r_k(\cdot):\ k=1, \dots,  N_2(\sigma)\} \subset L^\infty(B(p, r_0/4)), \\
& & R^*_M =R^*_{M,> } \cup R^*_{M, <}.
\eeq
{\ftext In Figure \ref{fig_6}, we consider the set
\beq\label{set X}
X=\{x_\beta\ :\ \beta\in \mathcal B\}\cup \{x_1, \dots, x_{ N_2(\sigma)}\}\subset M.
\eeq
Thus,}
{\Mtext denoting $\SlavaCfourteen= 2\SlavaCeleven+2\SlavaCthirteen+1$, see (\ref{3.21.12}) and } (\ref{distance-appr}),
 we obtain
\begin{lemma}\label{1st GH}
We have
\beq \label{5.11.4}
d_H(R_M(M), R^*_M) \le \SlavaCfourteen \sigma,
\eeq
where $d_H$ is the Hausdorff distance in $L^\infty(B(p, r_0/4))$.
\end{lemma}

\section{ Proof of Theorem \ref{l:10} and Proposition \ref{l:10-A}}
\subsection{From interior distance functions to boundary distance functions} \label{6.1}

By standard estimates for the differential
of the exponential map, see \cite[Ch.\ 6, Cor. 2.4]{Pe} the diameter of the sphere $\p B(p,r), \, r<r_0,$ is bounded
\beq\label{diameter of sphere}
\diam(\p B(p,r))\leq \pi r\,\cdotp  \frac{\sinh( \sqrt{K} r)}{\sqrt{K} r}\leq \pi r \,\cosh(\frac \pi 2)\leq 10\, r,
\eeq
where we use condition 
{\Rtext  (\ref{26.1}).} Let $N=M\setminus B(p,r_0/25)$.

\begin{lemma} \label{distance_across_sphere}
Let  $x \in M\setminus B(p,r_0/4)$
and $y\in \p N$ and ${z\in \p B(p,r_0/4)}$, let
\beq\label{distance through sphere1}
& &f(y, x, z)=d_N(y,z)+d_M(z,x),\\
& &f(y,x)=\min_{z_1\in \p B(p,r_0/4)}f(x, y; z_1),\nonumber
\eeq
where $d_N$  and $ d_M$ are the distances in $N$ and $M$, respectively.
Then,
\beq\label{distance through sphere}
d_N(y,x)=f(y,x)
\eeq
\end{lemma}

\noindent {\bf Proof.}
Clearly, as $d_M(z,x)\leq d_N(z,x)$ and  a  shortest curve
in $N$ from $y$ to $x$ intersects the sphere $\p B(p,r_0/4)$, we see that  $d_N(y,x)\geq f(y,x)$.

On the other hand let $z'=\hbox{ argmin}_z(f(y, x; z))$ and $\mu([0, f(y,x)]$ be the
corresponding union of the distance minimizing paths from $y$  to $z'$ and from $z'$ to $x$ for
which the minimum in (\ref{distance through sphere1}) is achieved.
Denote $s_1= d_N(y,z')$ and consider $\mu([s_1, f(y, x)]$. We  show next that
$\mu([s_1, f(y, x)] \subset N$.
 If this is not the case, there would exists $ s_1 < s_2 <s_3 <f(y,x)$
such that $\mu(s_1), \mu(s_3) \in \p B(p,  r_0/4)$, $\,\mu(s_2) \in \p B(r_0/25)$
and $\mu[s_3, f(y,x)] \subset M \setminus B(p,r_0/4)$. Then,
\beq \label{2.16.10}
s_1 \ge r_0\left(\frac14-\frac{1}{25}\right),  \,\, s_2 -s_1 \ge  r_0\left(\frac14-\frac{1}{25}\right),
\,\, s_3 -s_2 \ge  r_0\left(\frac14-\frac{1}{25}\right).
\eeq
On the other hand, consider a curve $\mu'([0, l])$ which is parametrised by the arclength
and consists of the radial path from $\mu(s_3)$
to $y' \in \p B(r_0/25)$ followed by a shortest
path along $\p B(r_0/25)$ from $y'$ to $y$. Due to (\ref{diameter of sphere}) and (\ref{2.16.10}),
\ba
l \le r_0\left(\frac{10}{25}+\frac14-\frac{1}{25}  \right) <3 r_0 \left(\frac14-\frac{1}{25} \right) \le s_3.
\ea
Taking the union of the path $\mu'([0, l])$, connecting
 $\mu(s_3)$
to $y'$, and the path
 $\mu(s_3, f(y,   x))$, connecting $y'$ to $x$, we get a contradiction to definition
(\ref{distance through sphere}). Thus, $\mu([s_1, f(y, x)]) \subset N$, i.e.,  $d_N(y, x) \le f(y, x)$.
\hfill\Box\medskip

Next, using the already constructed set $R^*$, see (\ref{8.07.1 B})
together with Lemmata \ref{1st GH} and \ref{distance_across_sphere},
 we construct a set $R^*(N) \subset
L^\infty(\p N)$  which {\Rtext approximates 
\beq \label{R_boundary}
R^{\p N}(N)= \{r_x^{\p N} \in L^\infty(\p N):  x\in N\},
\eeq
 where
$
r_x^{\p N}(z)=d_N(x,z),
\ \hbox{for } z\in \p B(p, r_0/25).
$}

\begin{lemma} \label{boundary_distance}
Let $R^*$  be the set given in (\ref{8.07.1 B}),
which satisfies (\ref{5.11.4}) be given. Then it defines a set $R^*(N) \subset L^\infty(\p N)$
such that
\beq \label{5.11.4 B}
d_H(R^{\p N}(N), \,R^*(N)) \le \RobCthirtyfive \sigma, \quad \RobCthirtyfive= 2\SlavaCfourteen+2\SlavaCeleven+1 .
\eeq
Here $\SlavaCfourteen$ is defined in (\ref{5.11.4}) and $\SlavaCeleven $ is defined in (\ref{3.21.12}).
\end{lemma}
Note that here we assume that $\delta $ satisfies (\ref{29.5}), $\sigma$ satisfies (\ref{4.11.20}) with the
related equations for $\e_4, \e_0,$ etc.
\medskip

\noindent{\bf Proof.} The proof is based on the construction of $R^*(N)$ which satisfies (\ref{5.11.4 B}).

Observe first that it follows from the proof of Lemma \ref{distance_across_sphere} that,
if $x, y \in B(p, r_0/4) \setminus B(p, r_0/25) \subset N$,
then
\bfo
d_N(x, y) \le \frac{r_0}{2} + \frac{8 r_0}{25},
\efo
so that a shortest path in $N$ connecting $x$ and $y$ lies in $B(p,  r_0)$. Thus it is possible, using (\ref{1.21.12}),
to construct an approximation $\tilde r_x^{ \p N}:\p N\to \R$ that satisfies
\beq \label{4.16.10}
\|r_x^{\p N}-\tilde r_x^{ \p N}\|_{L^\infty(\p N)} \le \RobCfortyfive \sigma,
\eeq
with a uniform constant $\RobCfortyfive$, cf.\ (\ref{3.21.12}). Denote
$R^*_<(N)=\{ \tilde r_x^{ \p N};\, x \in B(p, r_0/4) \setminus B(p, r_0/25) \}$, then
\beq \label{23.07.3}
d_H(R^{\p N}({B(p, r_0/4)}),\, R^*_<(N)) \le \RobCfortyfive \sigma,
\eeq
for $\delta< \delta_0$,
cf. construction of $R^*_<$ in subsection \ref{Distance functions approximation}.

Next, let
\bfo
R^*_c=\{r \in R^*:\, \min_{z \in \p N}(r(z)) \geq \frac{r_0}{8}  \}
\efo
For $y,z \in B(p, r_0/4) \setminus B(p, r_0/25)$ denote by $d^a_N(y,z)$ the distance between $y$ and $z$ in the metric $g^a$ along the curves lying in $B(p, r_0/2) \setminus B(p, r_0/25)$.
For each $r \in R^*_c$ we define
\beq \label{23.07.4}
& &\tilde r^{\p N} \in L^\infty(\p N): \,\, \tilde r^{ \p N}(y)=
\inf_{z \in \p B(p, r_0/4)} \left(d^a_N(z, y)+ r(z)  \right);
\\ \nonumber
& & R^*_>(N)=\{\tilde r^{ \p N}(\cdot):\, r \in R^*_c  \}.
\eeq
Then,
  with $R^*(N)= R^*_<(N) \cup R^*_>(N)$, we have that
\beq\label{Cfortyfour}
d_H(R^{\p N}(N),\, R^*(N)) \le (2\RobCfortyfive+\SlavaCfourteen)\sigma=\Cfortyfour\sigma.
\eeq
Here $\SlavaCfourteen \sigma$ error comes from an approximation of $d_N(y, z),$
see
 (\ref{5.11.4}), and
$2\RobCfortyfive \sigma$ error comes from approximating $d_M(z, x)$ and $d_N(y, z)$ in formula (\ref{23.07.3}), see also (\ref{distance through sphere1})-\eqref{distance through sphere} and (\ref{distance_across_sphere}). At last,
we use again that $\delta$ satisfies the  uniformly bound  (\ref{29.5}).
\hfill\Box\medskip
\\\\
\noindent
Recall that the metric tensor $g$  on $B(p,r_0)$ is a representation of a
metric in Riemannian normal coordinates and the $C^{2,\alpha}$-norm of the metric
is uniformly bounded. Using the fundamental equations of the Riemannian geometry, \cite[Ch.\ 2, Prop. 4.1 (3)]{Pe},
we have that the shape operator $S$ of the surface $\p B(p,r), \, r <r_0,$ can be given in the Riemannian
normal coordinates centered at $p$
in terms
of the metric tensor
as $S=g^{-1}\p_\nu g$, where $\nu$  is the unit normal vector of  $\p B(p,r)$. Taking $r=r_0/25$, we see
 that the  $C^{1,\alpha}$-norm of the shape operator $S$ of $\p N$ is uniformly bounded.
Also, by  
{\Rtext  (\ref{26.1})}
the boundary injectivity radius of $(N,g|_N)$ is bounded below by $\frac{24}{25} i_0$. As the sectional curvature
of $M$ and the second fundamental form (that is equivalent to the shape operator) of its submanifold $\p N$
are bounded, the Gauss-Codazzi equations imply that the sectional curvature of $\p N$ is bounded.
As the metric tensor of $M$ is bounded in normal coordinates in $B(p,r_0)$, we see that
the $(n-1)$-dimensional volume of  $\p N=\p B(p,r_0/25)$ is bounded from below by a uniform constant.
Thus by Cheeger's theorem, see \cite[Ch.\ 10, Cor.\ 4.4]{Pe},
the injectivity radius of $\p N$ is bounded from below by a uniform constant.

Summarising the above, the Ricci curvature of $(N,g|_N  )$   is uniformly bounded in $C^{\alpha}$, the
second fundamental form of $\p N$ is uniformly bounded in $C^{1,\alpha}$,
and the diameter and injectivity radii of $N$ and $\p N$, and the boundary injectivity radius of $(N, \p   N)$ are
uniformly bounded.
By \cite{KaKuLa}, using the knowledge of the set, $R^*(N)$ of approximate boundary distance functions,
which are $\RobCthirtyfive \sigma-$Hausdorff close
to  the set, $R^{\p N}(N)$  of the boundary distance functions  of manifold $(N,g|_N)$, one can construct
on the set  $R^*(N)$  a new distance function
$d^*_N:\, R^*(N) \times R^*(N) \to \R_+,$ such that
\beq\label{final estimat}
d_{GH}((N,d_N),(R^*(N),d^*_N))\leq \RobCthirtysix( \RobCthirtyfive \sigma)^{1/36},
\eeq
with a uniform $\RobCthirtysix>0$.

Having constructed ($R^*(N),d^*_N)$ we can now construct an approximate metric space $(M^*, d^*_M)$
which is $ \RobCthirtysix( \RobCthirtyfive \sigma)^{1/36}-$ close to $(M, d_M)$.
Indeed, let $x, y \in N$ and $\mu[0, l],  \,l=d_M(x, y)$
be a shortest between $x$ and $y$. If  $\mu[0, l] \subset N$ then $d_M(x, y))=d_N(x, y)$. If, however,
$\mu[0, l]$ intersects with $B(p, r_0/25)$  then, due to the convexity of $B(p, r_0/25)$,
there are $0<s_1 <s_2<l$
such that
\bfo
\mu[0, s_1] \subset N, \quad \mu[s_1, s_2] \subset B(p, r_0/25), \quad \mu[s_2, l] \subset N.
\efo
Therefore,  similar to Lemma \ref{distance_across_sphere}, we obtain
\begin{corollary} \label{distance-everywhere}
Let $x, y \in N$. Then
\beq \label{24.07.6}
& &d_M(x, y) =
\\
\nonumber & & \min \left(d_N(x, y),\, \min_{z_1, z_2 \in \p B(p, r_0/25)} \left[d_N(x, z_1)+
d_M(z_1, z_2)+d_N(z_2, y)  \right]  \right).
\eeq
\end{corollary}

Next define, for $\tilde r^{\p N}_1, \tilde r_2^{\p N} \in R^*(N)$,
\beq \label{24.07.10}
& &d^*_M(\tilde r_1^{\p N},\, \tilde r_2^{\p  N})= \\
& &\nonumber\min \left(d^*_N(\tilde r_1^{\p N}, \tilde r_2^{\p N}),\,
\min_{z_1, z_2 \in \p B(p, r_0/25)} \left[\tilde r_1^{\p N}(z_1)+ d^a(z_1, z_2)+\tilde r_2^{\p N}(z_2)  \right]  \right)
\eeq
Using (\ref{24.07.6}) together with (\ref{3.21.12}),
(\ref{final estimat}) and
 (\ref{1.21.12}),
we see that 
\beq \label{24.07.8}
\hspace{-4mm}d_{GH}((N, d_M), \,(R^*(N), d^*_M) \le  (2 \RobCthirtysix+1)( \RobCthirtyfive \sigma)^{1/36} \  \hbox{if}\,\,
\RobCfortyfive \sigma \le (\RobCthirtyfive \sigma)^{1/36}\hspace{-1mm}.
\eeq
Here $(N, d_M)$ is  the manifold $N$ with the distance function inherited from $M$ and
 $\delta <\delta_0$, cf. (\ref{23.07.3}).
\\
 Let us define the disjoint union
$M^*= R^*(N) \cup B(p, r_0/25) $. Next we  define a metric $d^*_M$ on this set.
To this end, consider first $\tilde r^{\p N} \in R^*(N),\, y \in B(p, r_0/25)$. Recall,
see the proof of Lemma \ref{boundary_distance},
  that the set $R^*(N)$ is bijective with $R^*_c \cup \left(B(p, r_0/4) \setminus B(p, r_0/25)\right)$.
In the case when $\tilde r^{\p N}$ is obtained from $r \in R^*_c$, we define $d_M^*(\tilde r^{\p N}, y)=r(y)$.
 Moreover, in the case when
 $\tilde r^{\p N}$ is obtained from $x \in B(p, r_0/4) \setminus B(p, r_0/25)$, we define $d_M^*(\tilde r^{\p N}, y)=
 d^a(x, y)$.
At last, if $x, y \in B(p, r_0/25)$, we take $d_M^*(x, y) =d^a(x, y)$.

 It follows from (\ref{24.07.8}) together with equations
(\ref{1.21.12}), (\ref{3.21.12}),
(\ref{5.11.4}) and considerations preceding Lemma \ref{1st GH} that
\beq \label{24.07.7}
d_{GH}((M^*, d^*_M), (M, d_M)) \le (2 \RobCthirtysix+1)( \RobCthirtyfive \sigma)^{1/36}.
\eeq
Summarizing, we obtain
\begin{lemma} \label{final_metric} Let $R^*$  satisfy  (\ref{5.11.4}) and
$M^*= R^*(N) \cup B(p, r_0/25) $ with metric $d^*_M$. Then,
\beq \label{5.16.10}
d_{GH}((M, d_M), (M^*, d_M^*)) \le \RobCfortyseven \sigma^{1/36}, \quad
\RobCfortyseven= (2 \RobCthirtysix+1) \RobCthirtyfive^{1/36}.
\eeq
\end{lemma}

\subsection{Proof of Theorem \ref{l:10} and Proposition \ref{l:10-A}}
{\ktext {\bf Proof of Proposition \ref{l:10-A}.}}
To prove the statement of the 
{\dtextRBX  Proposition}, 
we collect all the previous estimates. 
The aim is to find the relation between the final error $\e$
(i.e. $d_{GH}((M, d_M), (M^*, d_M^*)) \le \e$)
and the initial error $\delta$.
We proceed by following the chain of {\Mtext relations:
\beq \label{23.07.5}
\e \mapsto \sigma \mapsto \e_4 \mapsto \e_0 \mapsto \e_1 \mapsto \ga \mapsto \e_2
\mapsto \origjnull \mapsto \janul \mapsto \delta.
\eeq
{\Mtext To obtain inequality (\ref{13.10.2017b}) from
(\ref{5.16.10}) we set
\ba
\sigma = \left( \frac{\e}{\RobCfortyseven} \right)^{36}
\ea
 and use it in
%
\eqref{4.11.20},
\eqref{4.11.15}, (\ref{5.16.10}) and \eqref{epsone defined} with $\Lambda_s=1$
to determine values of $\e_4$, $\e_0$ and $\e_1$ by setting} 
\beq\label{new-e0e1} 
&&\hspace{-5mm}\e_4 =\frac{ \e^{36n}}{4 \RobCone \RobCfortyseven^{36n}},\   \e_0 = 
\Slavacten \e_4 \le \frac{1}{10},
\hbox{ and}\\
 \nonumber
& &\hspace{-5mm}\e_1 = \RobCforty \e^{72n} \ \mbox{with }
\RobCforty=\frac{ \Slavacten^2 }{160 \RobCone^2 \RobCfortyseven^{72n}}.
\eeq
{\Mtext 
\dtextRB 
To define $\ga$ so that  \eqref{22.06.1}, \eqref{gamma and sigma},
and \eqref{pre 2 new condition to e1 e2}  are valid, we set  
\beq\label{new-e0e2} 
\ga= \RobCfortyone \e_1^{1/b(s)}, \mbox{ with }\; \RobCfortyone = \min\big(\RobCforty^{-1/(2n)}\RobCfortyseven^{-36}, \SlavaCthreeagain, r_0/32 \big).
\eeq
Here we have used that $\sigma = \e_1^{1/(2n)}\RobCforty^{-1/(2n)}\RobCfortyseven^{-36}$ and noticed that $b(s) < 2n$.
}
 From (\ref{epsone defined})  and  \eqref{19.1}
 we get
 \beq \label{19.1 B}
\hspace{-6mm}\e_2\hspace{-1mm}=
  \frac{(\RobCfortyone \e_1^{1/b(s)})^{s/(s-1)}}{\big( \exp\big[ \big(   4L 
\RobCthirty \e_1^{-1}  (\RobCfortyone \e_1^{1/b(s)})^{-2+\theta/2} \; \exp(\RobCfortyone^{-c_{200}} \e_1^{-c_{200}/b(s)}) 
  \big
  )^{1/\beta
  } \big]  \big)^{ \frac{s}{(s-1)}}} \;
 \eeq
{\novtekst{with $\e_1$ given by (\ref{new-e0e1}).}}
 Finally, to choose $\origjnull$ and $\delta$ so that  \eqref{24Bnew}, (\ref{def SlavaCthree(s)}), \eqref{janolla}
 \eqref{29.5} 
 are satisfied,
 we set
  \beq\label{def SlavaCthree(s) again}
\origjnull =\hat \jbasis( \frac{\e_2}8;  \derho, 1)=
\SlavaCthree   \RobCfortyone^{-n}   8^{-n/s}  \e_1^{-n/b(s)}
\e_2^{-\frac ns}, 
\eeq
{\novtekst{with $\e_2$ given by (\ref{19.1 B}),}}
and choose $\delta$ so that
\beq
\nonumber 
\delta&\leq &2^{-n/2}
\Slavacsevenagain \RobCthree^{-n}  (\origjnull)^{-1}\e_2 
 \; = \;
2^{-n/2}
\Slavacsevenagain \RobCthree^{-n}
\SlavaCthree^{-1}   \RobCfortyone^{n}  8^{ n/s}
  \e_1^{n/b(s)}\e_2^{1+\frac ns}\\ 
 &=& {\RobCeigthy \,\e_1^{\RobCthirtynine}}
{\exp\left[-\RobCeightyone \,\e_1^{-\RobCfortyeight} \exp( \RobCeightytwo \,\e_1^{-\RobCfortynine}) \right]},  \label{1.2017 pre}
\eeq
with
\beq
\nonumber &&
C_{34}= \frac{1}{b(s)}\Big(n+ \frac{s+n}{s-1}\Big)
, \;\;
C_{35}= \frac{(s+n)}{(s-1)}\Big( 4L C_{12} \RobCfortyone^{-(2- \frac{\theta}{2})}\Big)^{1/\beta},
\\ \nonumber &&
C_{36}=\frac{\RobCfortyone^{-c_{200}}}{\beta}, \;\;
C_{37}= \frac{1}{\beta}\Big(1+ \frac{1}{b(s)}\Big(2-\frac{\theta}{2}\Big)\Big), \;\;
 C_{38}=\frac{c_{200}}{b(s)},
\\ &&
C_{39}= 2^{-n/2} \Slavacsevenagain \RobCthree^{-n} \SlavaCthree^{-1}    8^{ n/s} \RobCfortyone^{(s+n)/(s-1)+n}. \label{1.2017 preB}
\eeq
We use the inequality $x\le \exp(x)$  to bound from below the right hand side of the  estimate above to obtain,
by calling $  \CCfortythree= \max(C_{34},C_{37},C_{38},1/(2n))$,
\beq\label{CCfortytwo}
\delta \le
{\exp\left[-\exp((C_{39}^{-1}+ C_{35} +  C_{36} )\e_1^{- \CCfortythree}) \right]}.
\eeq
Notice that \eqref{3.21.12} is also satisfied,  by replacing  $C_{39}$ with $ \CCfortytwo=\SlavaCnineteen/(C_{33}^{36} \RobCforty^{1/(2n)})$ and by including $1/(2n)$ in $ \CCfortythree$.
Assuming $0< \delta \le \exp(-e)$,
we get
\beq \label{ineq_delta_e1}
{( \CCfortytwo^{-1}+ C_{35} +  C_{36} )}/{\ln\big(\ln\frac{1}{\delta}\big)} \le
\e_1^{ \CCfortythree},
\eeq
{\dtextRB
Let $\tau_0$ be the uniform constant introduced in Proposition \ref{uniform_covering} and define
 \beq
 \label{44}
 C_{44} = \min\big(1000^{- \CCfortythree},\, \RobCforty^{ \CCfortythree}(\RobCfortyseven^{36} \tau_0/2)^{2n \CCfortythree} \big).
 \eeq
In this way we can set in \eqref{ineq_delta_e1} the two constraints \eqref{tau and sigma} and $\e_1 \le 1/ 1000$ (derived from \eqref{epsone defined} with $\Lambda_s=1$)  and obtain
{\dtextRBX
$$ \delta \le \delta^*,\; \mbox{ with } \delta^* = \min\big( \exp(-e),\; 
\exp[-\exp [C_{44}^{-1}( \CCfortytwo^{-1}+ C_{35} +  C_{36} )]]\big).$$
}
}
Finally by using \eqref{new-e0e1} to rewrite $\e_1$ in \eqref{ineq_delta_e1}, and defining $\RobCninety=1/(72n \CCfortythree)$ and
$ \RobCfortythree = ( \CCfortytwo^{-1}+ C_{35} +  C_{36} )^{\RobCninety}\, \RobCforty^{-1/(72n)}$,
we obtain \eqref{13.10.2017b}.
\hfill\Box\medskip
\\
{\ktext {\bf Proof of Theorem \ref{l:10}.}}
Let $\delta \le \delta^*$ and let the ISD of $M^{(i)}, \, i=1,2$ be $\delta$-close. 
Take the finite collection
\bfo
\mathcal D = \left((B^e(r_0), g^{(1)}),\, \{(\la_j^{(1)},\, \varphi_j^{(1)}) \}_{j=0}^{J} \right), 
\efo
where the index $(1)$ is related to the IDS of $M^{(1)}$.
By construction the data $\mathcal D_0$ are $\delta$-close to the ISD of both $M^{(1)}$ and
$M^{(2)}$. By 
{\dtextRBX Proposition \ref{l:10-A}} 
the metric space $(M^*, d^*_M)$ constructed with these data is
$\e-$close to both $(M^{(i)}, d^{(i)}), i=1,2,$
where $\e$ is given by  the right hand side  of \eqref{13.10.2017b}.
We then conclude by triangular inequality,   
for any $0 <\delta \le \delta^*$,
\beq\label{final 1}
d_{GH}((M^{(1)}, d^{(1)}), \, (M^{(2)}, d^{(2)})) \le 2\e 
\eeq
We now extend this estimate to the case $\delta \in (0, \exp(-e)]$, when $\delta^* < \exp(-e)$.
To this end, observe that the definition of the GH-topology and \eqref{26.1} imply that:
$
d_{GH}((M^{(1)}, d^{(1)}), \, (M^{(2)}, d^{(2)})) \le D,
$
for any $\delta$.
By combining the latter inequality and \eqref{final 1} 
we obtain the inequality \eqref{12.06.3} with
$
\RobCeightyfour=
\max\Big(2\RobCfortythree,
D\big(\ln\big(-\ln\delta^*\big)\big)^{\RobCninety}
\Big)
$.
\hfill\Box\medskip

\noindent
 {\bf  Acknowledgements} RB and ML were partially supported by Academy of Finland, projects
 303754, 284715 and 263235.
YK was partially supported by
 EPSRC grant EP/L01937X and Institute Henri Poincare.
\vfill \newpage
{\ktext \begin{center}
{\bf Table of constants $C_k$, $c_k$, $\tau_0$ and $s$.}
\\
$  $
\\
 Note that all constant depend on
$n, R, D, i_0, r_0$ and variables  in brackets.
\\
$  $
\\
\begin{tabular}{p{2cm}|p{5cm}|p{2cm}|p{5cm}}
Name &  Introduced in / Notes &Name &  Introduced in / Notes
\\
\hline
$\RobCeightyfour$
 &     Thm. \ref{l:10} &
$\RobCninety$ &    Thm. \ref{l:10} \\ \hline
$
\RobCfortythree$ &    Prop. \ref{l:10-A}  &
$
\RobCzero 
$ &     Cor. \ref{Lip-stability} \\ \hline$
\RobCzeroone 
$ &     Cor. \ref{Lip-stability} & $
\RobCone
$ &    \eqref{12.06.13} \\ \hline$
C^{(har)}$
 &    \eqref{26.2A}&
$
C^{(Lip)}
$ &   \eqref{uniform Lip}
\\ \hline$
\RobCthree
$ &   \eqref{26.5}&
$
\RobCfive
$ & Prop. \ref{uniform_covering} 
 \\ \hline$
\RobCsix
$ &  Prop. \ref{uniform_covering}  &
$
\RobCseven
$ & Prop. \ref{uniform_covering}   \\ \hline$
\MattiCseven
$ & Lemma  \ref{ga-separation}& 
$\RobCthirty(\theta)$
& see \eqref{30.06.5}, we use $\theta = 1/2$ \\ \hline$
\RobCeleven
 $ &
\eqref{Sobo1}, we use
 $\a=1-\theta/2$ & $c_{200}$ &\eqref{19.06.1}\\ \hline
 $c_{205}{\ktext ( \theta)}$ & \eqref{19.06.1}, Appendix & $c_{206}{\ktext(\ga,\theta)}$ & \eqref{19.06.1}, Appendix\\ 
 \hline
  $
\slavaConetheta
 $&  \eqref{SlavacnonumberA}
 & $\tau_0$ &Prop. \ref{uniform_covering}
  \\ \hline$
\SlavaCnonumber
  $ &Lemma \ref{y-lemma 1} & 
$
\Slavacnonumber(s)
$ &Lemma \ref{y-lemma 1}  \\ \hline$
b(n)
$ & Lemma \ref{y-lemma 1} & 
$
\Slavacseven(s)
$ & Lemma \ref{lemma: tilde u exists} \\ \hline$
s
$ & \eqref{eq s}& 
$
\SlavaCthree(s)
$ & \eqref{def SlavaCthree(s)} \\ \hline$
\SlavaCone(s;\derho)
$ &Lemma \ref{lemma: tilde u exists} & 
$
\SlavaCsixagain
$ & \eqref{problem} \\ \hline$
\Slavacsevenagain
$ & \eqref{29.5}&
$\SlavaCnine(s)$ 
&Lemma \ref{lemma: main slicing}
  \\ \hline
$\Slavacten$ & Lemma \ref{volumes}
&
$
\SlavaCeleven
$ & \eqref{3.21.12} \\ \hline$
\SlavaCnineteen
$ & \eqref{3.21.12}&  $
\SlavaCtwelve
$ & \eqref{lem: admissible betas} \\ \hline$
\SlavaCthirteen
$ & \eqref{11.10.2017b}& 
$
\SlavaCfourteen
$ &  \eqref{1st GH}\\ \hline$
\SlavaCthreeagain
$ & \eqref{pre 2 new condition to e1 e2}& $
\RobCthirtyfive
$ & Lemma \ref{boundary_distance} \\ \hline$
\RobCfortyfive
$ &\eqref{4.16.10} &$
\Cfortyfour
$ & \eqref{Cfortyfour} \\ \hline$
\RobCthirtysix
$ & \eqref{final estimat}& $
\RobCfortyseven
$ &\eqref{5.16.10}  \\ \hline$
\RobCthirtynine
$ & \eqref{1.2017 preB}& $
\RobCeightyone
$ & \eqref{1.2017 preB} \\ \hline$
\RobCeightytwo
$ & \eqref{1.2017 preB}& $
\RobCfortyeight
$ & \eqref{1.2017 preB} \\ \hline$
\RobCfortynine
$ & \eqref{1.2017 preB}&
$
\RobCeigthy
$ & \eqref{1.2017 preB} \\ \hline$
\RobCforty
$ & \eqref{new-e0e1} & 
$
\RobCfortyone
$ &  \eqref{new-e0e2} 
\\ \hline$
\CCfortytwo
$ &\eqref{CCfortytwo} & 
$
\CCfortythree
$ &  \eqref{ineq_delta_e1}
\\ \hline$
\CCfortyfour
$ &\eqref{44} & 
$
\
$ & \\ \hline\end{tabular}
\end{center}
}


\section{Appendix}
\subsection{\dtextRBX Calculation of $c_{206}(\ga,\theta)$ in Theorem \ref{loc_stability}}
{\ztext 
To prove Theorem \ref{loc_stability} we need to show that the solution $w$ of the wave equation
\beq
\nonumber
P(y,D) w(y)= \tilde{q}(y),\quad y \in \Omega(T)\subset \R^{n+1}
\eeq
can be estimated in the set $\D(z,\gamma,T)\subset
\R^{n+1}$ in \eqref{2.3} that is between two double cones of the spacetime,
i.e. $\Sigma(z,\gamma,T) \subset \D(z,\gamma,T)  \subset \Sigma(z,0,T)$. Here, $\gamma>0$ is the
parameter that indicates how close the set $\D(z,\gamma,T) $  is to the optimal double cone $ \Sigma(z,0,T)$.
Theorem \ref{loc_stability} is proven by applying a proper iterative procedure and the dependency of the
coefficient $c_{206}$ on $\gamma$ is crucial for our considerations.
%
%
%
%
%
\\
The calculation of  $c_{206}$ can be consider as the final step of a long geometric construction.
In order to understand it we summarize the previous steps with related references.\\
In Section 3 of \cite{B} we calculated the parameters of the inequality associated with a (conormally) pseudo-convex function $\psi$ with respect to the wave operator $P(y,D)$.
Then we used this property to calculate the coefficients of the Tataru inequality (recalled in the following section \ref{sub-tataru})
\begin{eqnarray*}
 \|e^{-\epsilon |D_0|^2/2\tau} e^{\tau f} u\|_{1,\tau} \le c_{1,T}\, \tau^{-1/2}\|e^{-\epsilon |D_0|^2/2\tau} e^{\tau f} P(y,D)u\|_{0} + c_{2,T} e^{-\tau R_2^2/4\epsilon} \|e^{\tau f} u\|_{1,\tau}
\end{eqnarray*}
 and to prove the local stability of the unique continuation for the wave operator. \\
In  \cite{BKL} we used the previous result to prove the global stability of the unique continuation for the wave operator. As recalled in the following section \ref{sub-global}, the proof is based on the iteration $N$ times of the a local stability for the 'low temporal frequency' component of the solution $u$ of the wave equation:
\begin{eqnarray*}
\|A(D_0/\omega)b((y-y_j)/r)u_j\|_{H^1} \le c_{155,j}\exp({-c_{132}\mu_j^{\alpha^2}}),\quad \forall\, \omega \le \mu_j^\alpha/(3c_{131}).
\end{eqnarray*}
Moreover in Section 3.1. of \cite{BKL} and Appendix A of \cite{BKL}  we applied the stability result in the domain of influence of a cylinder. 
\\
{\zztext
In the case of the present paper, the mentioned domain of influence is called $\Sigma(z,0,T)$ in Theorem 2 and, according to the iterative procedure, it contains a covering of the set $\Lambda=\D(z,\gamma,T)$. The balls of the covering have radius $2R$, that depends on the distance to the boundary,
on the regularity and pseudo-convexity property of the function $\psi$, and on extra constraints imposed by the Tataru inequality. The local stability step holds for smaller balls with radius $r$, and $r < R$. In Table 1 we summarize these values and in particular we obtain, up to a multiplicative constant,
\beq \label{r formula}
r\sim \gamma^{58}.
\eeq
The $\sim$ symbol is defined precisely in subsection \ref{sub-gamma}.  
By  construction and for \eqref{r formula}, one can calculate the number of local steps of the iteration (see also Table 2)
\beq \label{N formula}
N\sim \gamma^{-58(n+1)}.
\eeq
These two values combined with the calculation of the coefficients for the local and the global stability lead to the following relationship between $c_{206}$ and $\gamma$ 
$$c_{206}  = \zeta_1 \,\exp(  \gamma^{-\zeta_2}), \mbox{  for proper positive numbers } \zeta_1, \zeta_2.$$
We will prove that formula \eqref{N formula} plays a big role for the calculation of $\zeta_2$.\\
In both articles \cite{B,BKL} 
}
we used consistent notations for the geometric quantities and labeled the important coefficient as $c_h$, with a unique $h\ge 100$ in order to be able to follow the construction of the final parameters. One can find them  in those papers by searching for the corresponding index $h$.\\
Here in this Appendix our focus is the dependency of all the parameters (in particular $c_{206}$) on the quantity $\gamma$,
since this reflects the cost of getting close to the cone of dependence. For this reason in the following section \ref{sub-gamma} we quickly introduce the main relationship between $\gamma$ and the used Gevrey function localizers, and in the next sections we recalculate the main coefficients of the above results and we summarize them in the Tables 1 and 2.\\
We will follow the same notation as in \cite{B, BKL}. Unfortunately it was not possible to use an analogous notation in the rest of the present paper. Anyway we will point out the different notations.

\subsubsection{Gevrey functions and dependency on $\gamma$}\label{sub-gamma}
\underline{Assumption}. Let $\alpha \in [1/3,1)$, and let $T, \ell, \gamma$ be defined as in Assumption A5, \cite{BKL}. 
\\
(In the present paper, this corresponds to conditions \eqref{22.06.1} and \eqref{2.3}).\\\\
Gevrey functions are used as smooth localizers in the constructions and 
their main properties are outlined in Section 4 of \cite{B}.
In particular in our calculations we consider the following Gevrey function (see \cite{RD}, Ex 1.4.9 for definition):
\begin{eqnarray*}
\chi_1(t) = \chi(1+t)\chi(1-t), \quad \mbox{with }\chi(s) = \exp(-s^{\frac{\alpha}{\alpha-1}})
\mbox{ for s $>$ 0}, \; \chi(s)=0 \mbox{ for s $\le$ 0}.
\end{eqnarray*}
One can slightly modify the definition such that $\chi_1 = 1$ in a ball $B_1 \subset \R$  (with radius 1), $\chi_1 = 0$ outside the ball $B_2$ (with radius 2), and $0\le \chi_1 \le 1$. Observe that $\chi_1 \in G^{1/\alpha}_0(\R)$ 
{\ztext since}

\begin{eqnarray}\label{chi1}
|D^\kappa \chi_1(v)| \le c_{0X}c^{|\kappa|}_{1X} |\kappa|^{|\kappa|/\alpha},
\quad \mbox{with } c_{0X}=O(1), \,c_{1X}=O\Big(\frac{1}{1-\alpha}\Big).
\end{eqnarray}
Here the symbol $O$ (big-O) means
``comparable up to a an absolute multiplicative constant to''
(i.e. $A = O(B)$ implies ${c_{abs} \le A/B \le C_{abs}}\,$, for some positive numbers $c_{abs}, C_{abs}$). \\
Furthermore,
define  $\chi_\delta(v) := \chi_1(v/\delta)$, $v \in \R^M$.
\\
Hence, $\mathcal{F}_{v \to \zeta} \chi_\delta(v) = \delta^M \mathcal{F}_{v \to \delta\zeta} \chi_1(v)$ for
 $\zeta \in \C$, and calling $c_{2X} = 1/(eMc_{1X})^{\alpha}$ we get
\begin{eqnarray}\label{def_chi2}
\nonumber
|\mathcal{F}_{v \to \zeta} \chi_\delta(v)| &\le& \delta^M c_{0X} \exp(\delta H_{B_2}(\mbox{Im}\zeta)-c_{2X}\delta^{\alpha} |\mbox{Re}\zeta|^{\alpha})\cdot \mbox{Vol(supp}(\chi_1),dv).
\\
\end{eqnarray}
Product estimate: for $v \in B_2(\R^M)$, calling $c_{0X,l},c_{1X,l}$ (resp. $c_{0X,m},c_{1X,m}$) the coefficients in \eqref{chi1} for $\chi_l$ (resp. $\chi_m$),
\begin{eqnarray*}\label{def_chi3}
|D^\kappa \chi_l(v)\chi_m(v)| &\le& c_{0X,l}c_{0X,m} \max\{c_{1X,l},c_{1X,m} \} (\max\{c_{1X,l},c_{1X,m} \})^{|\kappa|}|\kappa|^{|\kappa|/\alpha}.
\end{eqnarray*}
%
We start by linking $\gamma$  and the coefficient $(1-\alpha)$, since both quantities tend to zero.\\\\
\underline{Assumption}. We assume $\theta = 1/2$. Next, from now on the symbol\\
$\sim$ means ``comparable up a multiplicative coefficient independent of $\gamma$ or $(1-\alpha)$\,\,to''. \\
(i.e. $A \sim B$ implies $c \le A/B \le C$, with $c,C$ {\wtext independent} 
 of  $\gamma$ or $(1-\alpha)$). The multiplicative coefficient is in general a uniform geometric constant, in the sense specified at the beginning of the paper. 
\\\\
We will call $c_{205}$ the resulting multiplicative coefficient for $c_{206}$.
\\
$\theta$ is the exponent appearing in the global stability of the unique continuation (see  Theorem \ref{loc_stability} of this paper and the following section \ref{sub-global}), while $1/\alpha$ is the order of the used Gevrey functions.
According to \cite{BKL} (end of page 6469), by construction these two values are related in the following way:
\begin{equation}\label{alphaN}
\alpha^N =\theta \mbox{, that implies } \alpha^{1/r^{(n+1)}} =  \frac{1}{2} \quad \Rightarrow (1- \alpha) \sim r^{n+1} \mbox{, as } \alpha \to 1,
\end{equation}
where $N=c_{170} \sim \gamma^{-58(n+1)}$ 
is reported in the following Table 2 and  $r  \sim \gamma^{58}$ is in Table 1. 
Consequently, for $\chi_1$ and $c_{1X}$ defined above, we get
\begin{eqnarray}\label{c1X}
c_{1X} \sim \frac{1}{1-\alpha} \sim \frac{1}{\gamma^{58(n+1)}},\qquad
|\chi_1'|_{C^0(\Omega_0)} \sim c_{1X}, \qquad |\chi_1''|_{C^0(\Omega_0)} \sim c_{1X}^2. 
\end{eqnarray}
\subsubsection{Tataru inequality and Table 1}\label{sub-tataru}

We consider the wave operator in $\R^{n+1}$,
\begin{eqnarray}\label{wave_op}
P(y,D) = -D_0^2+ \sum_{j,k=1}^n g^{jk}(x)  D_j D_k +  \sum_{j=1}^n h^j(x) D_j + q(x),
\end{eqnarray}
where  $y=(t,x)\in \R \times \R^{n}$ are the time-space variables, $D_0 = -i \partial_{t}$, $D_j = -i \partial_{x_j}$.
The  coefficients $g^{jk} \in C^1(\R^n)$ are real and independent of time, and  $[g^{jk}]$ is a symmetric positive-definite matrix.
The  coefficients $h^j, q\in C^{0}(\R^n)$
are complex valued and independent of time. 
Call $\xi=(\xi_0,\tilde\xi)$ the Fourier dual variable of $y=(t,x)$.
{\ztext  In the next theorem we use the exponential pseudodifferential operator}
$e^{-\epsilon |D_0|^2/2\tau} v= \mathcal{F}^{-1}_{\xi_0 \to t}e^{-\epsilon \xi_0^2/2\tau}\mathcal{F}_{t'\to\xi_0}v$, with $\mathcal{F}$ and $\mathcal{F}^{-1}$ representing respectively the Fourier transform and its inverse. 
Let us also define
\begin{eqnarray}\label{def_f}
f(y)= \sum_{|\upsilon|\le 2} (\partial^\upsilon {\phi})(y_0)\, (y-y_0)^\upsilon / \upsilon! - \sigma |y-y_0|^2.
\end{eqnarray}
In following theorem (called Theorem 2.1 in \cite{B}) we recall the Carleman-type estimate by Tataru, named `Tararu inequality'.
\begin{theorem}
\label{th_carleman}(\cite{B}, Theorem 2.1; or \cite{BKL}, Theorem 2.3.) Let $\Omega$ be an open subset of $\R\times\R^n$.
Let $P(y,D)$ be the wave operator \eqref{wave_op}, with $g^{jk}(x) \in C^1(\Omega)$, $h^j, q \in C^{0}(\Omega)$.
Let $y_0 \in \Omega$ and $\psi \in C^{2,\rho}(\Omega)$ be real valued, for some fixed $\rho\in(0,1)$, such that $\psi'(y_0)\neq 0$ and $S=\{y;\psi(y)=0\}$ being an  oriented hypersurface non-characteristic
in $y_0$. \\
Consequently there is $\lambda>1$ such that $\phi(y)=\exp(\lambda \psi)$ is a conormally strongly pseudoconvex  function with respect to $P$ at $y_0$.
\\
Then there is a real valued quadratic polynomial $f$ defined in \eqref{def_f} with proper $\sigma > 0$,
and a ball $B_{R_2}(y_0)$ such that $f(y) < \phi(y)$ when $y \in B_{R_2}-\{y_0\}$ and $f(y_0)=\phi(y_0)$; and $f$ being a conormally strongly pseudoconvex function with respect to $P$ in $B_{R_2}$.
This implies that there exist $ \,\epsilon_0,\, \tau_0,\,c_{1,T},\,c_{2,T}, R$, such that, for each small enough  $\epsilon<\epsilon_0$ and large enough $\tau > \tau_0$, we have
\begin{eqnarray*}
 \|e^{-\epsilon |D_0|^2/2\tau} e^{\tau f} u\|_{1,\tau} \le c_{1,T}\, \tau^{-1/2}\|e^{-\epsilon |D_0|^2/2\tau} e^{\tau f} P(y,D)u\|_{0} + c_{2,T} e^{-\tau R_2^2/4\epsilon} \|e^{\tau f} u\|_{1,\tau}.
\end{eqnarray*}
Here $u \in H^1_{loc}(\Omega),$ with $P(y,D)u \in L^2(\Omega)$ and supp$(u) \subset B_{R}(y_0)$.
\end{theorem}

\noindent
\underline{Assumption}:  We now consider the  `hyperbolic function' 
\begin{eqnarray}\label{psi1}
\psi(t,x;T,z) = (T - d_g(x,z))^2 - t^2
\end{eqnarray}
introduced in Definition 3.1 of \cite{BKL}, and its level set $\psi(y) - \gamma^2 = 0$.
%
\\\\
Starting from  a 
{\ztext  general}
$\psi$, in section 3 of \cite{B}, page 180, we
have already calculated all the geometric constants associated either with the related pseudo-convexity estimates of $\psi$ or with the Tataru inequality. 
They are summarized in Table 1, page 191 of \cite{B}, and are copied in Table A.3. of \cite{BKL} 
(with few modifications explained in the related Appendix A).
Then in section A.1.1. of \cite{BKL}, page 6487, we have recalculated these quantities for the particular case of  the `hyperbolic function' $\psi$ in \eqref{psi1}. 
\\
Our aim here is to start from Table A.3. and section A.1.1. of \cite{BKL}  in order to 
 find the $\gamma-$dependency of those coefficients.\\
The following new Table 1 must be read from the top to the bottom, since it starts with the basic inequalities and continues with more complicated expressions.\\
As said, we assume $\psi$ as in \eqref{psi1}, and calculate all coefficients in the Tataru inequality.
The first two values $C_l =\min|\psi'(y)|$ and $p_1=\min p(y, \psi')$ are defined at page 6484 of \cite{BKL},  Section A.1., Assumption b).
Their limit value is calculated in \cite{BKL}, formula (A.7): i.e. 
$C_l = 2 \gamma_I b_0^{-1/2}$,  $p_1=4\gamma_I^2$.
Since $\gamma_I=\gamma/\sqrt{2}$ (see Lemma A.3.a, page 6489), and $b_0$ is defined as a constant (see formula (3.1), page 6452), then 
the $\gamma-$dependency of the two coefficients is respectively $\gamma$ and $\gamma^2$, as shown in the table. 
The third value of Table 1 is $dist\{\partial \Omega_0,  \Omega_a\}$ (alias $dist(\Lambda, \partial \Omega_0)$) and behaves like $\gamma^2$, thanks to the estimates (A.12) and (A.11) in \cite{BKL}.
On the other hand, the following coefficients until $C_3$ in Table 1  are independent from $\gamma$, because of formula (A.6) and (A.8) in \cite{BKL}.\\
The next  values in Table 1 are obtained by substituting the upper values:\\
$M_P, M_1,..., R_1$, defined in section 3.1 of \cite{B};\\
$c_T \sim \lambda^3$ (replacing $4n|\lambda \psi|_{max,\Omega_0}$), see (A.2) and Remark A.1 in \cite{BKL};\\
$\tau_0, c_{1,T}, c_{2,T}, c_{133}$ defined in section 3.2 of \cite{B};\\
$r, \delta, R$, defined in section 3.3 of \cite{B}, here we have renamed $r_0$ by $r$.\\
These 3 coefficients are used to prove Proposition 2.5 of \cite{B}, which is related to the result of local stability for the unique continuation.\\
$c_{111}$ of \cite{B} is not used here.\\
Note that $\sigma, r, \delta, R, \tau_0, R_1, R_2$ have nothing to do with quantities with the same name used in the rest of this paper (outside from the Appendix). 
 
\subsubsection{Global stability coefficients and Table 2}\label{sub-global}

{\zztext
This section can be seen as an overview of the proof of Theorem \ref{loc_stability}, with the final estimate for $c_{206}$.\\
We introduce the main steps and we always follow the notation of \cite{B,BKL} to better follow the calculations.
}
\\\\
%
\noindent
\underline{Assumption}:
Define a net of center points $(t_k,z_k)$ for the translated hyperbolic functions:
$$\psi(y;T_k,z_k,t_k) = (T_k-d_g(x,z_k))^2-(t-t_k)^2.$$
Let $\Upsilon=W(z,T,\ell)$ be the initial cylinder  (called $\Gamma$ in \eqref{2.3}) 
 and let $\Sigma(z,\ell,T)$ be the related domain of influence (in the paper called $\Sigma(z,0,T)$ according to \eqref{2.3b}).
\\
We choose the domains for the covering 
$$\Omega_{0,k} \subset \{y; y\in [-T_k +t_k, T_k+ t_k ]\times \R^n;\, \psi(y;T_k,z_k,t_k) \ge \gamma_k^2/2, \,T_k\ge d_g(x,z_k)\}$$ and 
$$\Lambda_{k} \subset \{y; y\in [-T_k +t_k, T_k+ t_k ]\times \R^n;\, \psi(y;T_k,z_k,t_k) \ge \gamma_k^2, \,T_k\ge d_g(x,z_k)\}.$$ Let $\gamma_k \ge \gamma$, for all k.
\\
The construction is similar to the one in Figure 1, page 6470 of \cite{BKL}. \\
The parameters $(t_k,z_k,T_k,\gamma_k)$ should be chosen such that the $x-$projection of $\Omega_{0,k}$ is  contained  in the domain $0 < d_g(z_k,x) \le \frac{7}{8}i_0$, that is within the injectivity radius $i_0$, in order to guarantee the $C^{2,\rho}$-regularity of $\psi(y;T_k,z_k,t_k)$.
Moreover 
the union $\Lambda= \bigcup_{k=1}^K \Lambda_{k}$ should cover a subset of the domain of influence $\Sigma(z,\ell,T)$.
For example, let $\Lambda = S(z,\ell,T,\gamma)$, (alias $\D(z,\gamma,T)$  in \eqref{2.3}).
\\\\
The above construction together with the assumption on the Gevrey-regularity of the localizers 
let us apply Theorem 1.2 in \cite{BKL}. The details of Assumptions A2-A3 can be checked in the paper, while $P$ is defined in \eqref{wave_op}.

\begin{theorem}\label{global2}(\cite{BKL}, Theorem 1.2)
Under the conditions of Assumptions A2-A3,
define the open set
$\Omega_1 = \bigcup_{k=1}^K \Omega_{0,k} \backslash \overline{\Upsilon}$ containing $\Lambda$.
Then for every $0<\theta < 1$ we have
\begin{eqnarray*}
 \|u\|_{L^{2}(\Lambda)} \le c_{161} \frac{\|u\|_{H^{1}(\Omega_1)}}{\Big(\ln\Big(1+\frac{\|u\|_{H^{1}(\Omega_1)}}{\|Pu\|_{L^{2}(\Omega_1)}}\Big)\Big)^{\theta}}.
\end{eqnarray*}
Moreover, for any $m \in (0,1]$ we get
\begin{eqnarray*}
\|u\|_{H^{1-m}(\Lambda)} \le c_{161}^m \frac{\|u\|_{H^1(\Omega_1)}}{\Big(\ln \Big(1 + \frac{\|u\|_{H^1(\Omega_1)}}{\|Pu\|_{L^2(\Omega_1)}}\Big)\Big)^{m \theta}}\,.
\end{eqnarray*}
The constant $c_{161}$ is calculated in the proof.
\end{theorem}
Up to a uniform multiplicative constant (and according to Remark 3.8. of \cite{BKL}), we can identify the constant $c_{161}$ with our final constant $c_{206}$,  even if the first one is defined for a bounded domain of the Euclidean space and the second one is defined for a compact manifold $(M,g)$. Indeed by assumption, in each chart of $M$ holds the inequality 
$$a_0 I\leq [g_{jk}(x)]_{j,k=1}^n\leq b_0\, I,\quad \hbox{and}\quad \|g_{jk}\|_{C^4(M)}\leq b_3, \; a_0<1<b_0,$$
which let one approximate all spatial subdomains to an Euclidean ball.
%
\\
\\
%
Theorem 1.2 is a generalization of Theorem 1.1 in \cite{BKL} 
{\ztext for a more complex domain, but with a similar final estimate where the inverse-log term has a different multiplicative constant.}
For each  $\Omega_{0,k}$ Theorem 1.1 in \cite{BKL} holds with constant $c_{160}$ in place of $c_{161}$.
The number 
{\ztext $K$ of the used sets}
 $\Omega_{0,k}$ is by construction proportional to the number of charts covering the domain. This number depends on 
{\ztext the bounds for the diameter, the injectivity radius and the harmonic radius  of $M$, called respectively $D$, $i_0$ and  $r^{(har)}$ in the notation of the paper.}
Hence we can also write $c_{206} \sim c_{161} \sim c_{160}$.
\\\\
%
The technique used to prove the above Theorem consists in iterating the local stability result, but considering the low temporal frequencies separately from the high temporal frequencies.
The pseudodifferential operator $A(D_0/\omega)$ defined below is used to localize the low temporal frequencies of the solution $u$, where the estimate is more complicated.\\\\
\underline{Assumption}:
We consider a pseudo-differential operator $A(D_0)$  with symbol $a(\xi_0) \in G^{1/\alpha}_0(\R)$, $0\le a\le 1$,
supported in $|\xi_0|\le 2$ and equal to one in $|\xi_0|\le 1$.
Hence we can write $A(\beta |D_0|/\omega)v = \mathcal{F}^{-1}_{\xi_0 \to t}a(\beta |\xi_0|/\omega)\mathcal{F}_{t'\to\xi_0}v$.
We fix $a$ as in \eqref{chi1}.
\\\\
Another complication comes from the fact that the local stability result holds just in small balls $B_r(y_j)$, centered in $y_j$ with radius $r$.
{\ztext 
It is important for our estimates that the balls $B_{r}(y_j)$, with $j=1,\dots, N$,  cover the set   $\Lambda$. 
We will choose the center points $y_j$ in the set ${\cal E}$, so that the union of the balls is contained in the domain of influence of the cylinder $\Upsilon$, i.e. $\bigcup_{k=1}^N B_{r}(y_j) \subseteq  \bigcup_{k=1}^K \Omega_{0,k} \subset \Sigma(z,\ell,T)$.
}
\\
 Furthermore there are particular  conditions on the support of $u$ to be fulfilled, 
{\ztext  also affecting the set ${\cal E}$ and the iteration.}
\\
Hence in the final domain $\bigcup_{k=1}^K \Omega_{0,k}$ the local stability result must be applied several times to a sequence $\{u_j\}_{j=2}^ N$ of proper cut-offs of the solution $u$.
Let $u_j$ be defined as:
\begin{eqnarray}\label{uj}
u_j = \prod_{k=1}^{j-1}\big(1-b_k\big)u, \quad b_{k}:=b\Big(\frac{2(y-y_k)}{r}\Big).
\end{eqnarray}

Then we can introduce the following Theorem 2.7 in \cite{BKL}, formulating a local stability estimate (of the unique continuation for the wave operator) of inverse exponential type for the low temporal frequencies of $u_j$.
\\
The exact construction of the radii $r$ and $R$ is in Proposition 2.5 of \cite{B}, as intersection of several geometric and analytic constraints.
The $\gamma$-dependency of $r$ and $R$ is shown in Table 1.
In particular we get $r \sim \gamma^{58}$.
The number of balls used in the iteration is $N=c_{170} \sim \gamma^{-58(n+1)}$, as shown in Table 2.
The constant $c_{170}$ is defined in formula (2.5) of \cite{BKL}.
\\
Notice that at each step we reduce the support of the temporal localizer $A(D_0)$, by defining the term $\mu_j = c_{156} \mu_{j-1}^{\alpha}$.\\
{\ztext We will show that 
$c_{155,N} \sim {\gamma^{- \zeta_4}} c_{155,1}\sim {\gamma^{- \zeta_5}},$ 
and that $c_{161}  \sim N \gamma^{- \zeta_6}  c_{156}^{-\frac{\alpha}{(1-\alpha)}}$,
for proper positive numbers $\zeta_4, \zeta_5, \zeta_6$.
}
\\
The details of Assumptions A1-A2-A4 and of the set ${\cal E}$ can be checked in the paper.

\begin{theorem}\label{iter}(\cite{BKL}, Theorem 2.7)
Under the Assumptions A1-A2-A4,
let $y_k \in {\cal E}$
and let
$b\in G^{1/\alpha}_0(\R^{n+1})$ be a Gevrey functions of class $1/\alpha$ with compact support, such that $0<\alpha<1$.
\\
Then, there exist constants $R,r$ with $R \ge 2r >0$, and $c_{159}>1$
such that
 for $\mu > c_{159}$ there are coefficients $c_{151}, c_{152}, c_{154},c_{155},c_{156}, \beta, N$ for which,\\
if
\begin{eqnarray}\label{ipo1}
\|u\|_{H^1(\Omega_1)}= 1,\quad \|Pu\|_{L^2(\Omega_1)}
< 1
,\quad \|A\big(\frac{D_0}{\beta \mu}\big)l(y)Pu\|_{L^2} \le  \exp(-\mu^{\alpha}),
\end{eqnarray}
then calling  $\mu_1=\mu$ and $\mu_j = c_{156} \mu_{j-1}^{\alpha}$ for $2 \le j\le N$, we have $\mu_j \ge 1$ and
\begin{eqnarray}\label{ipo2}
\|u_j\|_{H^1(B_{2R}(y_j))}\le c_{152},\quad \|Pu_j\|_{L^2(B_{2R}(y_j))} \le c_{153},
\end{eqnarray}
\begin{eqnarray}\label{ipo3}
\|A(D_0/\mu_j)b((y-y_j)/R)Pu_j\|_0 \le c_{154,j}\exp({-\mu_j^\alpha}),
\end{eqnarray}
and consequently
\begin{eqnarray}\label{ipo4}
\|A(D_0/\omega)b((y-y_j)/r)u_j\|_{H^1} \le c_{155,j}\exp({-c_{132}\mu_j^{\alpha^2}}),\quad \forall\, \omega \le \mu_j^\alpha/(3c_{131}).
\end{eqnarray}
The radii $r$ and $R$ are defined in Table A.3, while the coefficients $c_{k}$ are calculated in the proof of the Theorem.
\end{theorem}

In the following Table 2 we show the $\gamma$ dependency of the coefficients used in the proof of Theorem 1.2 and Theorem 2.7 in \cite{BKL}.
The coefficients $c_h$ of the local stability are defined in \cite{B} and recalled also in the proof of Lemma 2.6. of \cite{BKL}, page 6459.
As said, the index $h$ is unique and here we briefly remind the definition of $c_h$ and the relationship with other coefficients and with Table 1.
\\
In \eqref{c1X} we obtained the $\gamma$ dependency for  $c_{1X}$, i.e.  $c_{1X}  \sim  {\gamma^{-58(n+1)}}$.
It follows that $c_{2X} \sim 1/(c_{1X})^{\alpha} $, where $c_{2X}$ is the coefficient in \eqref{def_chi2} (it was called $c_{102}$ in \cite{B}).\\
Therefore for simplicity 
{\ztext we give below}
the values in Table 2 in terms of their $c_{1X}$ or  $\gamma$ dependency. \\
In order to calculate the rest we need to refine some estimates.\\
First of all we recall and improve the coefficients in Lemma 2.1, \cite{BKL}, for the $L^2$ and $H^m$ norms:
\beq\label{c107-108}
& &c_{107}=  c_3 \Big(\frac{8}{\beta_1}\Gamma\Big(\frac{1}{\alpha}\Big)\frac{1}{\alpha (c_{117})^{\frac{1}{\alpha}}}\Big)^{\frac{1}{2}}\frac{1}{(\alpha c_{106})^{\frac{1}{\alpha}}},\\
& & \nonumber
c_{108}=  c_{107}(1+|D_x^m f|_{C^0}) + c_{107}\frac{(1+m)^{\frac{(m+1)}{\alpha}}}{(\alpha c_{106})^{\frac{m}{\alpha}}}  \\ & & \|A(\beta_1 D_0/\mu)f (1-A(D_0/\mu))v\|_1 \le c_{108} e^{-c_{106} \mu^\alpha} \|v_{\mbox{supp}(f)}\|_m\,. \nonumber
\eeq
Next, following Remark 2.8 (4) in \cite{BKL}, we split each smooth Gevrey localizer in time and space:
$$b\big(\frac{y-y_0}{R}\big)=b\big(\frac{t-t_0}{R}\big)b\big(\frac{x-x_0}{R}\big),$$
with $b(t) = \chi_1(t) \in G^{1/\alpha}_0(\R)$ (as in \eqref{chi1}) and $b(x) \in C^{2}_0(\R^n)$.
Consequently the functions $f_1(y),f_2(y),f_3(y)$ (see formula (2.21) in \cite{BKL}) can generally be written as:
$f_*(y)=f_*(t)f_*(x)$, with $f_*(t) = D_0^2 b_{j-1}(t)+D_0 b_{j-1}(t)+b_{j-1}(t)$ and
$f_*(x) = D_rD_s b_{j-1}(x)+D_r b_{j-1}(x)+b_{j-1}(x)$, for $b_{j-1}(t) := b(2(t-t_{j-1})/r)$.
Let $v = b((y-y_{j-1})/r)u_{j-1}$, then
\begin{eqnarray*}
\|A\Big(\frac{3D_0}{\nu}\Big)f_*(t)\Big(1-A\Big(\frac{D_0}{\nu}\Big)\Big)f_*(x)v\|_1
\le \|A\Big(\frac{3D_0}{\nu}\Big)(D_0 f_*(t))\Big(1-A\Big(\frac{D_0}{\nu}\Big)\Big)f_*(x)v\|_0\\
+ \|A\Big(\frac{3D_0}{\nu}\Big)f_*(t)\Big(1-A\Big(\frac{D_0}{\nu}\Big)\Big)(D_0+D_x+1)(f_*(x)v)\|_0
\le c_{108} c_{152} \exp(-c_{106}\nu^\alpha)
\end{eqnarray*}
with $c_{108}$ calculated as in \eqref{c107-108} with $\beta_1=3, m=3, c_3=(r/2) c_{0X}$. Moreover, we can 
recalculate the terms at page 6466 of \cite{BKL}:
\begin {eqnarray*}
c_{162,j} &=& 2 c_{162,j-1}+ c_{153}c_{164}+c_{155,j-1} |-P_2 b_{j-1} + h^s(x)D_{x_s}b_{j-1}|_{C^0} \\
&&+c_{107} c_{152}\big(1+ n^2|g^{kr}|_{C^0}  + |h^s|_{C^0} \big)
 + c_{155,{j-1}} |2  D_0 b_{j-1}|_{C^1}
 + c_{152} c_{108}  \\
&&+ c_{155,j-1}  |D_0(2  D_0 b_{j-1})|_{C^0}
+ c_{152}c_{107}
\\
&&+  c_{155,j-1} |2n g^{kr}  D_k b_{j-1}|_{C^1}
+c_{152} c_{108} n^2|g^{kr}|_{C^1} \\
&&+ c_{155,j-1} |D_r(2  g^{kr} D_k b_{j-1})|_{C^0}
+c_{107} c_{152}n^2|g^{kr}|_{C^1}
\end{eqnarray*}
\begin {eqnarray*}
c_{162,j}&\sim& 2 c_{162,j-1} + \Big( \frac{N^2 c_{1X}^2}{r^2} \Big)\frac{c_{1X}^{3/2}}{r^{1/2}} + c_{155,{j-1}}(1+ |g^{kr}|_{C^1}+|h^{s}|_{C^0})\Big(\frac{|b'|_0}{r}+\frac{|b''|_0}{r^2}+\frac{|b'|^2_0}{r^2}\Big)+\\
&&+ \Big( \frac{N c_{1X}}{r}\Big) c_{108} (1+ |g^{kr}|_{C^1}+|h^{s}|_{C^0})
\;\sim\; c_{162,j-1} + c_{155,j-1}\frac{c_{1X}^2}{r^2} \\
c_{154,j}&=&c_{162,j} + c_{153} \tilde c_{107} \sim c_{162,j} + \frac{N^2 c_{1X}^3}{r^2}R^n  \sim c_{162,j} \sim  c_{155,j-1}\frac{c_{1X}^2}{r^2} \\
c_{116} &\sim& \gamma^4 c_{154,j}^2 \Big(\frac{N c_{1X}}{\gamma^{48}}\Big)^4.
\end{eqnarray*}
By applying Lemma 2.6 in \cite{BKL} with $c_U=c_{152}$, $c_P=c_{153}$, $c_A=c_{154,j}$, one obtains:
\begin{eqnarray*}
c_{155,j}=c_{150}(c_{152},c_{153},c_{154,j}) \sim c_{1X}^3\frac{\sqrt{c_{116}}}{\gamma^{48}} \sim \frac{N^2 c_{1X}^5}{\gamma^{46+58\cdot 2}}c_{155,j-1},\\
c_{156} =\min \Big( \frac{1}{18\beta c_{131}},c_{132}^{1/\alpha},c_{165}^{1/\alpha}, \frac{c_{106}^{1/\alpha}}{3c_{131}}\Big)
= \frac{c_{106}^{1/\alpha}}{3c_{131}}  \sim \gamma^{56\alpha + 58(n+1)(\alpha+1)+28}.
\end{eqnarray*}
Now we can obtain the $\gamma$ dependency of $c_{160}$ in Theorem 1.1 of \cite{BKL}.\\
{\zztext
Recalling \eqref{alphaN} (i.e. $\alpha^N = \theta = 1/2$ and $(1-\alpha) \sim \gamma^{58(n+1)}$),
}
we get $c_{159} = c_{156}^{-\frac{1}{\alpha^{N-1}(1-\alpha)}} > 1$ and therefore:
\begin{eqnarray*}
c_{159}\sim \Big(\frac{1}{\gamma^{56\alpha + 58(n+1)(\alpha+1)+28}}\Big)^{\frac{1}{2\gamma^{58(n+1)}}}
=\exp\Big(\frac{-[56\alpha + 58(n+1)(\alpha+1)+28]}{2\gamma^{58(n+1)}}\ln (\gamma)\Big),
\end{eqnarray*}
\begin{eqnarray*}
c_{158}&=& N c_{155,N}+ 3N c_{131} c_{152} \Big(1+ \frac{|b'|_{C^0}}{r}\Big)c_{156}^{-\alpha/(1-\alpha)} \sim N c_{131} c_{152}c_{159}^{1/2},\\
c_{160}&=&\big(\ln (1+ e^{c_{159}})\big)^{1/2}+2^{1/2} c_{158} \sim c_{158}.
\end{eqnarray*}
Hence, $c_{160}$ of Theorem 1.1 in \cite{BKL} (and analogously $c_{161}$ in Th. 1.2.) fulfills the estimate
\begin{eqnarray*}
c_{160} \le \exp\Big(\frac{1}{\gamma^{c_{200}}}\Big), 
\mbox{ with  } c_{200} = 58(n+1) + 2.
\end{eqnarray*}
We know that $c_{206} \sim c_{160}$.
{\ztext We denote by}
$c_{205}$ the uniform multiplicative constant that depends on 
{\ztext the uniform geometric parameters  $T,i_0,D,r_0,R,n$, named according to the notation of the rest of the paper. 
}
The number $c_{205}$ also depends on $\theta$, that for simplicity has been fixed here equal to 1/2. 
{\ztext The above inequality gives an estimate for $c_{160}\sim c_{206},$ and thus we}
can conclude that $c_{206}{\ktext(\ga,\theta)}=c_{205}{\ktext ( \theta)}\exp\big({\gamma^{-c_{200}}}\big)$.
%
\\\\
{\bf Remark.}
Please notice that there was a misprint in the paper \cite{BKL}
both in the statements of Theorem 3.3. and in Corollary 3.9. 
However this misprint did not affect the calculations of the present paper (or the results in \cite{BKL}).
\\
Namely in Theorem 3.3., we have the following
erratum (in the denominator of the final inequality):
\begin{eqnarray}
\|u\|_{L^{2}(\Lambda)} \le 
c_{163} \frac{\|u\|_{H^1((\Omega_1)}}
{\Big(\ln \Big(e + \frac{\|u\|_{H^1((\Omega_1)}}{\|f\|_{L^2(\Omega_1)}}\Big)\Big)^{\theta}}.
\nonumber
\end{eqnarray}
And this is the corresponding corrigendum (replace $e$ with $1$):
\begin{eqnarray}
\|u\|_{L^{2}(\Lambda)} \le 
c_{163} \frac{\|u\|_{H^1((\Omega_1)}}
{\Big(\ln \Big(1 + \frac{\|u\|_{H^1((\Omega_1)}}{\|f\|_{L^2(\Omega_1)}}\Big)\Big)^{\theta}}.
\nonumber
\end{eqnarray}
In Corollary 3.9., we have erratum in (3.27), and corrigendum:
\begin{eqnarray}
\|w\|_{L^{2}(\Omega_2\backslash W_1)} \le 
c_{166} \frac{\|w\|_{H^1((\Omega_1\backslash W_0)}}
{\Big(\ln \Big(1 + \frac{\|w\|_{H^1((\Omega_1\backslash W_0)}}{C' \|w |_{W_1}\|_{H^1(W_1)}}\Big)\Big)^{\theta}}.
\nonumber
\end{eqnarray}

\bigskip
\noindent
{\bf Table 1 and Table 2.} We next present the two tables that summerize the previous calculations. They show the $\gamma$ dependency of the parameters. The name of the constants there is unique. 
The order of the parameters in Table 1 is always increasing in complexity, that is the parameters down may depend on the upper ones.
In general the same principle is followed also in Table 2, even if the relationships are more complex. For simplicity 
the values in Table 2 are expressed in terms of their $c_{1X}$ or  $\gamma$ dependency, where we recall that  $c_{1X}  \sim  {\gamma^{-58(n+1)}}$.

\newpage
\begin{center}
{\bf Table 1}
\begin{tabular}{p{2cm}|p{0.3cm}|p{12cm}}
Name &  & Order with respect to $\gamma$ 
\\
\hline
$C_l$ & $\sim$ & $\gamma$ \qquad (\cite{BKL}, formula (A.7))\\
\hline
$p_1$ & $\sim$ & $\gamma^2$  \qquad (\cite{BKL}, formula (A.7))\\
\hline
$dist\{\partial \Omega_0,  \Omega_a\}$ & $\sim$ & $\gamma^2$ \qquad (\cite{BKL}, formula (A.12))\\
\hline
$|\psi'|_{C^k}$ & $\sim$ & $1$ \qquad (\cite{BKL}, formula (A.8))\\
\hline
$d_g(x,z)$ & $\in$ & $[\ell, T-\gamma]$ \quad (in $\Gamma\backslash$ cylinder)\\
\hline
$|\partial_k d_g|$ & $\sim$ & $1$ \qquad (\cite{BKL}, formula (A.6))\\
\hline
$C_3$ & $\ge$ & $1$ \\
\hline
$M_P$ & $\le$ & 1\\
\hline
$M_1$ & $\ge$ & $ \frac{1}{(p_1)^2} = \frac{1}{\gamma^4} $ \\
\hline
$M_2$ & $\ge$ & $M_1= \frac{1}{\gamma^4}$\\
\hline
$\lambda$ & $\ge$ & $\max\{M_1,e,\frac{1}{C_l^2}\}= \frac{1}{\gamma^4}$ \\
\hline
$\phi_0$ & $\ge$ & $e^{-1}$ \\
\hline
$\phi_M$ & $\le$ & $e$ \\
\hline
$R_1$ & $\le$ & $\min\{1, \gamma^2, \frac{1}{\lambda}\}=\gamma^4$ \\
\hline
$c_T$ & $\sim$ & $\lambda^3=\frac{1}{\gamma^{12}}$ \quad  (\cite{BKL}, formula (A.2)  and Remark A.1)\\
\hline
$c_{100}$\quad & $\ge$ & $1$\\
\hline
$\epsilon_0$ & $\le$ & $\frac{1}{(\lambda (1 + \lambda )+ c_T)}=\frac{1}{\lambda^3}=\gamma^{12}$\\
\hline
$R_2$ & $\le$ &
$\min\Big \{R_1,\, \frac{ C_l}{ (1 + \lambda   + c_T/\lambda )},\,
\frac{\lambda^2  C_l^2}{c_T}$,
$\Big(\frac{1}{c_T^2  M_1(1+\lambda^2)}\Big)^{\frac{1}{4}},\, \frac{\epsilon_0}{\sqrt{2M_2}},\,\frac{\lambda } {c_T
\big(1+\lambda^2  +\lambda^2 (1 + \lambda) \big)}\Big\}$,
\\
 &  &
$=\min\{\gamma^4, \gamma^9,\gamma^6, \gamma^9, \gamma^{14}, \gamma^{20}\}=\gamma^{20}$\\
\hline
$\sigma$ & $\ge$  & $c_T R_2=\gamma^{8}$  \qquad  (\cite{BKL}, formula (A.2) and Remark A.1)\\
\hline
$\tau_0$ & $\ge$ & $  M_1  \Big(
\big(\lambda^2 +  c_T R_2\big)^2  + |h|^2_{C^{0}(\Omega_0)}(1+\big(\lambda  + c_T R_2^{2}\big)^2)+ |q|^2_{C^{0}(\Omega_0)} \Big) = \frac{1}{\gamma^{20}}$\\
\hline
$R$ & $\le$ & $R_2 = \gamma^{20}$ \\
\hline
$\delta$ & $\le$ & $c_T R_2^{3}=\gamma^{48}$\\
\hline
$r$ & $\le$ & $\frac{ \lambda^2 C_l^2 R_2^{3} }{\big( \lambda+c_{T} R_2^2 \big)} = \gamma^{58}$\\
\hline
$c_{1,T}$ & $\ge$ & $ \sqrt{ \Big(\frac{ M_1}{\tau_0} + \frac{1}{\lambda }\Big)}=\gamma^2$\\
\hline
$c_{2,T}$ & $\ge$ &
$
\sqrt{M_2}(1+\frac{|\chi_1'|_{C^0(\Omega_0)}}{\tau_0 R}) + \frac{c_{1,T}}{\sqrt{\tau_0}}\,c_{133}=\frac{1}{\gamma^{2}}+\frac{1}{\gamma^{8}}(|\chi_1''|_{C^0(\Omega_0)}+\frac{|\chi_1'|_{C^0(\Omega_0)}}{\gamma^4} ) \sim \frac{c_{1X}^2}{\gamma^8}
$
\\
\hline
$c_{133}$ & $\ge$ & $ \frac{|\chi_1''|_{C^0(\Omega_0)}}{\tau_0 R^2} + \frac{|\chi_1'|_{C^0(\Omega_0)}}{R}(1+\lambda  + c_T R_2^{2}+\frac{|h|_{L^\infty{(\Omega_0)}}}{\tau_0})\,=\frac{1}{\gamma^{20}}(|\chi_1''|_{C^0(\Omega_0)}+\frac{|\chi_1'|_{C^0(\Omega_0)}}{\gamma^4} )$
\\
\hline
\end{tabular}
\end{center}

\newpage
\begin{center}
{\bf Table 2}
\begin{tabular}{p{1cm}|p{6.5cm}|p{1cm}|p{6cm}}
Name &   Value & Name & Value \\
\hline
$c_{2X}$ & $=c_{102}=\frac{1}{(e c_{1X})^{\alpha}}$ & $c_{119}$ & $\delta c_{1X}\sim \gamma^{48}c_{1X}$\\
\hline
$c_{118}$ & $1 + |\phi'|_0(1+R_2) +5n|\phi''|_{0,\rho}R_2^{\rho+1} + |\phi''|_0(1+R_2^2) +  \sigma(2+R_2^2)\sim \frac{1}{\gamma^8}$ & $c_{114}$ & $c_{1,T}^2 |g|^2_{C^1}|\chi_1|^2_{C^2}(1+|\varphi'|^4_{C^0}/\delta^4+|\varphi''|^2_{C^0}/\delta^2)
 \sim \frac{c_{1X}^4}{\gamma^{12+48\cdot 4}}$\\
\hline
 $c_{115}$ &   $c_{2,T}^2(|\varphi'|^2_{C^0} +1)(3^3 e^{-3}/\delta^3)(1+|\chi_1'|^2_{C^0}/\delta^2)
\sim \frac{c_{1X}^6}{\gamma^{8\cdot 3+48\cdot 5}}$ & $c_{121}$ & $\frac{c_{1X}}{\delta}$\\
\hline
 $c_{122}$ & $\frac{c_{1X}^2}{\gamma^{44}}$ & $c_{123}$ & $\sim \frac{\gamma^{56 \cdot \alpha}}{c_{1X}^\alpha}$\\
\hline
 $c_{128}$ &   $\frac {1}{3^{\alpha}2}c_{123} \sim c_{123}$ & $c_{110}$  & $c_{122} \big(\frac{8\Gamma(1/\alpha)}{3[\alpha c_{123}^{1/\alpha}(\alpha c_{128})^{1/\alpha}]}\big)^{1/2}
\sim \frac{ c_{1X}^3}{\gamma^{44+56}}$\\
\hline
 $c_{109}$ &  $\min(\sqrt{\epsilon\,\delta/36}, c_{128}/2, 1)
 \sim \frac{\gamma^{56 \cdot \alpha}}{c_{1X}^\alpha}$ & $c_{130}$  & $ \frac{3 c_{109}}{4 \delta}\Big(\frac{1}{16 }\Big)^{5}
 \quad \sim \frac{\gamma^{56 \cdot \alpha - 48}}{c_{1X}^\alpha}$\\
\hline
 $c_{131}$ &  $\mbox{max}(16^6\sqrt{2}, \frac{16^6 3^{\alpha-1}\sqrt{2\epsilon_0 \delta}}{c_{123}}, \frac{16^6 \sqrt{\epsilon_0 \delta}}{3\sqrt{2}})  \sim\frac{c_{1X}^{\alpha}}{\gamma^{56\cdot \alpha - 30}}$ & $c_{135}$ & $r^{\alpha} c_{2X}\frac{1}{4 \cdot 3^{\alpha}} \sim \frac{\gamma^{58\cdot \alpha}}{c_{1X}^{\alpha}}$\\
\hline
 $c_{137}$ &  $\min(\frac{1}{2}\big(c_{102}\delta^{\alpha}\frac{(c_{130})^{\alpha}}{(\sqrt{2})^{\alpha}} +\delta \frac{c_{130}}{2 \sqrt{2}}\big), \frac{1}{2} c_{102}\delta^{\alpha}(\frac{1}{2\sqrt{2}}c_{130})^{\alpha})
 \sim \frac{\gamma^{48\cdot \alpha}}{c_{1X}^{\alpha}}c_{130}^{\alpha}$ & $c_{132}$ & $\min(c_{135},c_{137}) \quad
\sim \frac{\gamma^{56 \cdot \alpha \cdot \alpha}}{c_{1X}^{\alpha\cdot \alpha }}\frac{1}{c_{1X}^{\alpha}}$\\
\hline
 $c_{170}$ &  $N \sim \frac{1}{\gamma^{58(n+1)}}$ &  $c_{117}$ & $(r/2)^\alpha \frac{1}{(ec_{1X})^\alpha} \sim \frac{r^\alpha}{c_{1X}^\alpha}$\\
\hline
 $\tilde c_{117}$  & $c_{2X} R^{\alpha} = (e c_{1X})^{-\alpha}  R^{\alpha}$ & $\beta $ &  $2 + (\frac{4}{\tilde c_{117}})^{1/\alpha}  \sim \frac{c_{1X}}{R} \quad (2.11)$\\
\hline
 $\tilde c_{106}$  &   $  \frac{1}{\beta^\alpha} \sim \frac{R^\alpha}{c_{1X}^\alpha}$ & $\tilde c_{107}$ & $R^{n+1} c_{0X} \Big(\frac{8}{\beta}\Gamma\Big(\frac{1}{\alpha}\Big)\frac{1}{\alpha (\tilde c_{117})^{1/\alpha}}\Big)^{1/2}\frac{1}{(\alpha \tilde c_{106})^{\frac{1}{\alpha}}} \sim R^{n}c_{1X}$\\
\hline
 $c_{154,1}$  &   $1 + \tilde c_{107} \sim R^{n}c_{1X}$ & $c_{155,1}$ & $\max(c_{134}, c_{136}) = \max(c_{1X}^{2.5} \gamma^{58(n-\frac{3}{2})} , \frac{c_{1X}^6}{\gamma^{180}})=\frac{c_{1X}^6}{\gamma^{180}}$\\
\hline
 $c_{153}$ &   $1+2N\big(1+ n^2|g^{kr}|_{C^0} + |h^s|_{C^0}\big)\big(\frac{|b'|_{C^0}}{r} + \frac{|b''|_{C^0}}{r^2}+ (N-1)\frac{|b'|^2_{C^0}}{r^2}\big) \sim \frac{N^2 c_{1X}^2}{r^2}$ & 
$c_{152}$ & $2\big(1+ N \frac{|b'|_{C^0}}{r}\big) \sim \frac{Nc_{1X}}{r}$\\
\hline
 $c_{162,1}$ &   $1$ & 
 $c_{156}$
 & 
$ 
\sim \frac{c_{106}^{1/\alpha}}{3c_{131}}  
\sim \gamma^{56\alpha + 58(n+1)(\alpha+1)+28}
$
 \\
\hline
 $c_{165}$ &   $c_{117}\beta^\alpha/(3^\alpha 4) \sim \frac{r^\alpha}{R^\alpha} \sim \gamma^{38\alpha}$ & $c_{164}$ & $\frac{r}{2} c_{0X} \big(\frac{8}{3}\Gamma\Big(\frac{1}{\alpha}\Big)\frac{ec_{1X}}{\alpha^{1/\alpha} (r/2)} \big)^{1/2} \frac{e c_{1X} (3^\alpha 4)^{\frac{1}{\alpha}} }{(\alpha^{\frac{1}{\alpha}} (r/2))} \sim \frac{c_{1X}^{3/2}}{r^{1/2}}$\\
\hline
 $c_{107}$ &   $c_{164} \sim \frac{c_{1X}^{3/2}}{r^{1/2}}$ & $c_{108}$ & $\Big(c_{107} + c_{107} \frac{4^{4/\alpha}}{(\alpha c_{106})^{3/\alpha}}\Big)
\Big(1+\frac{|b'|_0}{r}+\frac{|b''|_0}{r^2}+\frac{|b'''|_0}{r^3}\Big) \Big(1 + \frac{|b'|_0}{r}\Big)
\sim \frac{c_{1X}^\frac{17}{2}}{r^\frac{15}{2}}$\\
\hline
\end{tabular}
\end{center}


\end{document}